\documentclass{amsart}
\usepackage{graphicx}
\usepackage{amssymb}
\usepackage[all]{xy}

\input epsf
\newtheorem{theorem}{Theorem}[section]
\newtheorem{cor}[theorem]{Corollary}
\newtheorem{lemma}[theorem]{Lemma}
\newtheorem{prop}[theorem]{Proposition}
\newtheorem{example}[theorem]{Example}
\theoremstyle{definition}

\newtheorem{rem}[theorem]{Remark} 

\theoremstyle{plain}
\numberwithin{equation}{section}

\newcommand{\C}{\mathbb{C}}

\newcommand{\Hom}{{\operatorname{Hom}}}
\newcommand{\N}{{\mathbb{N}}}
\newcommand{\TL}{Temperley--Lieb}
\newcommand{\TLD}{D}
\newcommand{\type}{\Delta}

\providecommand{\noglossaryignore}[1]{}
\newcommand{\globalglossaryentry}[3]{\makebox[1.5in][l]{\tt $\backslash${#1}} 
\makebox[1.1in][l]{{$#2$}} \makebox[2.5in][l]{{#3}}\newline} 
\newcommand{\newcommandabbreviation}[3]{\newcommand{#1}{#2}%
\noglossaryignore{\globalglossaryentry{#3}{#2}{}}}
\newcommand{\renewcommandabbreviation}[3]{\renewcommand{#1}{#2}%
\noglossaryignore{\globalglossaryentry{#3}{#2}{}}}
\newcommand{\newcommandmacro}[4]{\newcommand{#1}{#2}%
\noglossaryignore{\globalglossaryentry{#3}{#2}{#4}}}
\newcommand{\gge}[3]{\noglossaryignore{\globalglossaryentry{#1}{#2}{#3}}}
\newcommand{\myaddress}%
{\parbox{3in}{\footnotesize \begin{center} 
Mathematics Department, City University, \\  
Northampton Square, London EC1V 0HB, UK.\end{center}}}
\newcommand{\twlrm}{}    
\newcounter{minidef}[section]

\newcounter{minicapt}

\newtheorem{pr}[theorem]{Proposition} 
  
\newtheorem{lem}[theorem]{Lemma} 
\theoremstyle{definition}
\newtheorem{de}[theorem]{Definition}     
\noglossaryignore{GREEK ETC.\newline}
\newcommandabbreviation{\e}{\epsilon}{e}        
\newcommandabbreviation{\lam}{\lambda}{lam}  
\newcommandabbreviation{\la}{\langle}{la}        
\newcommandabbreviation{\ran}{\rangle}{ran}
\newcommandabbreviation{\ha}{\#}{ha}             
\newcommandabbreviation{\rmap}{\rightarrow}{rmap}
\newcommandabbreviation{\aaa}{\alpha}{aaa}        
\newcommandabbreviation{\ab}{\alpha,\beta}{ab}
\newcommandabbreviation{\aab}{a(\ab )}{aab}       
\noglossaryignore{\newline RINGS\newline}
\newcommandabbreviation{\HH}{H \!\!\! I}{HH}               

\renewcommandabbreviation{\Re}{\mathbb R}{Re}
\newcommandabbreviation{\Q}{\mathbb Q }{Q}
\renewcommandabbreviation{\H}{\mathbb H }{H}
\noglossaryignore{\newline SYMMETRIC GROUP\newline}
\def\Sym(#1){\Sigma(#1)}                   
\gge{Sym(-)}{\Sym(-)}{symmetric group on - objects}
\def\Sy(#1){\Sigma_{#1}}                   
\gge{Sy(-)}{\Sy(-)}{symmetric group irreducible -}
\def\sym(#1){\mbox{\LARGE s}(#1)}        
\gge{sym(-)}{\sym(-)}{symmetric group on - objects (variant)}
\def\sy(#1){\mbox{\LARGE s}({#1})}        
\gge{sy(-)}{\sy(-)}{symmetric group irreducible - (variant)}
\newcommandmacro{\cs}{\C \, \sy(n)}{cs}{symmetric group algebra over $\C$}
\noglossaryignore{\newline PARTITIONS/SETS\newline}
\newcommand{\Nset}[1]{\underline{#1}}
\gge{Nset\{-\}}{\Nset{-}}{set of natural numbers to -} 
\def\nset(#1){ \{ #1 \}_{ \underline{n} }} 
\gge{nset(-)}{\nset(-)}{a set $-\times\Nset$}
\def\ul(#1){_{\underline{#1}}}             
\gge{ul(-)}{{}\ul(-)}{subscript underline -}
\def\Ee(#1){{\bf E}_{#1}}                  
\gge{Ee(-)}{\Ee(-)}{set of equivalence relations on set -}
\def\Eee(#1){{\bf E}_{\{ #1 \}_{\underline{n}}}}   
\gge{Eee(-)}{\Eee(-)}{ditto for nset}
\def\Een(#1,#2){{\bf E}_{\{ #1 \}_{\underline{#2}}}}   
\def\Ssn(#1,#2){{\bf S}_{\{ #1 \}_{\underline{#2}}}}   
\def\Ss(#1){{\bf S}_{#1}}                  
\def\Sss(#1){{\bf S}_{\{ #1 \}_{\underline{n}}}}   
\def\bbc(#1){((\beta_1)(\beta_2)...(\beta_{#1}))}      
\newcommandmacro{\Ln}{{\Gamma}^{n}}{Ln}{large index set}
\newcommandmacro{\LnQ}{{\Gamma}^{n}_Q}{LnQ}{index set}
\newcommandmacro{\Zz}{\zeta}{Zz}{`shape' function}
\noglossaryignore{\newline PARTITION ALGEBRA\newline}
\def\ka(#1){\kappa_{#1}}                   
\def\Sm(#1){\Sigma_{#1}}                   
\newcommandmacro{\com}{\bullet}{com}{bullet composition}
\newcommandmacro{\enm}{\; e^n(\! m\! ) \;}{enm}{product of idempotents}
\def\Ai(#1){ A^{ #1 \cdot } }              
\def\Aij(#1,#2){ A^{ #1  #2 } }            
\newcommandmacro{\One}{\mbox{\bf $1 \!\!\! 1$}}{One}{algebra unit 1}
\newcommandmacro{\Bp}{B_p}{Bp}{partition basis}
\def\Bb(#1){B_p[#1]}                       
\def\Pp(#1){P_n[#1]}                       
\def\Ps(#1){P_n[#1] \! /}                  
\newcommandmacro{\Ph}{\hat{P}}{Ph}{P hat  algebra}
\def\Is(#1){\sim^{#1}}                     
\noglossaryignore{\newline MODULES\newline}
\def\Wm(#1){{\cal S}_{#1}}                 
\gge{Wm(-)}{\Wm(-)}{Weyl module with index -}
\def\wm(#1,#2){{}_{#1}{\cal S}_{#2}}       
\gge{wm(-1,-)}{\wm(-1,-)}{Weyl module with index -}
\def\Ind(#1,#2,#3){\mbox{Ind}_{#1}^{#2}#3} 
\gge{Ind(-1,-2,-)}{\Ind(-1,-2,-)}{induction}
\def\Res(#1,#2,#3){\mbox{Res}_{#1}^{#2}#3} 
\gge{Res(-1,-2,-)}{\Res(-1,-2,-)}{restriction}
\newcommandabbreviation{\weyl}{standard}{weyl}

\newcommandabbreviation{\head}{\mbox{head }}{head}
\newcommandabbreviation{\Weyl}{Weyl}{Weyl}
\def\SS(#1){{\cal S}_{#1}}                 
\gge{SS(-)}{\SS(-)}{Specht/Weyl module index -}
\def\LL(#1){{\cal L}_{#1}}                 
\gge{LL(-)}{\LL(-)}{simple module index -}
\noglossaryignore{\newline FUNCTORS/MAPS\newline}
\newcommandmacro{\Gg}{{\cal G}}{Gg}{G Functor}
\newcommandmacro{\Fg}{{\cal F}}{Fg}{F Functor}
\newcommandmacro{\ra}{\rightarrow}{ra}{}
\def\ses(#1,#2,#3){0\ra #1 \ra #2 \ra #3 \ra 0}   
\gge{ses(1,2,-)}{\ses(1,2,-)}{\hspace{.5in} short exact sequence}
\def\starr(#1){ \stackrel{ #1 }{\longrightarrow} }
\gge{starr(-)}{\starr(-)}{}
\newcommandmacro{\doublerightarrow}{\; -\!\!\! -\!\!\!\!\!\! \gg \;}
{doublerightarrow}{}
\noglossaryignore{\newline PARTITION ALGEBRA MAPS\newline}
\newcommandmacro{\smap}{s}{smap}{`inclusion' map}
\newcommandmacro{\tmap}{t}{tmap}{$ P_n -> S_n$}
\newcommandmacro{\pmap}{\psi}{pmap}{$ S_n -> P_n $}
\noglossaryignore{\newline MISC.\newline}
\def\Amap(#1){{\cal A}_{#1}}               
\gge{Amap(-)}{\Amap(-)}{}
\def\Rr(#1){R_{#1}}                        
\gge{Rr(-)}{\Rr(-)}{restriction of E}
\def\Cr(#1){C_{#1}}                        
\gge{Cr(-)}{\Cr(-)}{restriction of E to N}
\newcommandmacro{\Tm}{{\cal T}}{Tm}{Transfer Matrix}
\def\On(#1){{\cal I}_{#1}}
\gge{On(-)}{\On(-)}{}
\newcommandmacro{\UU}{\underline{\sqcup}}{UU}{}  
\newcommandmacro{\UUU}{\sqcup}{UUU}{}  
\newcommandmacro{\Vq}{V_Q^{\otimes n}}{Vq}{Potts config. space}
\def\bs(#1,#2){\mbox{{\Large $\ast$}}^{#1}_{#2}}  
\gge{bs(-,-)}{\bs(-,-)}{general plumbing multiplier}
\newcommand{\ignore}[1]{}
\gge{ignore\{-\}}{\ignore{-}}{ignore argument!}
\def\choo(#1,#2){ \left( \begin{array}{c} #1 \\ #2 \end{array} \right) } 
\gge{choo(-1,-)}{\choo(-1,-)}{choose}
\newcommand{\Qed}{$\Box$}
\gge{Qed}{\mbox{\Qed}}{QED}
\def\staq(#1){\stackrel{#1}{=}}            
\gge{staq(-)}{\staq(-)}{}
\def\stam(#1){\stackrel{#1}{\rightarrow}}  
\gge{stam(-)}{\stam(-)}{}
\def\mat{ \left( \begin{array} }    
\def\tam{ \end{array}  \right) }
\gge{mat/tam}{...}{matrix delimiters}
\newcommand{\beq}{\begin{equation} }
\def\eql(#1){ \begin{equation} \label{#1} 
}
\newcommand{\eq}{\end{equation} }
\def\eqal(#1){\begin{eqnarray} \label{#1} }
\def\eqa{\end{eqnarray} }
\def\lab(#1){\label{#1}
}
\def\prl(#1){ \begin{pr} \label{#1} 
}
\def\del(#1){ \begin{de} \label{#1} 
}
\gge{smeq\{-\}}{...}{small equation}
\gge{fneq\{-\}}{...}{very small equation}
\noglossaryignore{\newline HECKE/BLOB\newline}
\newcommandmacro{\Hnq}{H_n(q)}{Hnq}{ * freestanding symbol}
\newcommandmacro{\Hn}{H_n}{Hn}{      *-mod etc.}
\newcommandmacro{\A}{{\cal A}}{A}{}
\newcommandmacro{\Cwts}{C}{Cwts}{}
\newcommandmacro{\CA}{{\cal A}}{CA}{}

\newcommandmacro{\calA}{{\cal A}}{calA}{}
\newcommandmacro{\modi}{\mbox{Mod} }{modi}{was mod not modi!}
\newcommandmacro{\Wgen}{{\Bbb S}}{Wgen}{}
\def\ol(#1){\overline{#1}}
\newcommandmacro{\st}{\mbox{St}}{st}{}
\def\CMult(#1,#2){(#1:#2)}
\def\CM(#1,#2){( #1 : #2 )}
\def\FMult#1,#2{(#1:#2)}
\def\CF#1,#2{(#1:#2)}

\newcommandmacro{\Top}{\mbox{Top}}{Top}{}
\newcommandmacro{\Soc}{\mbox{Soc}}{Soc}{}
\newcommandmacro{\Head}{\mbox{Head}}{Head}{}
\newcommandmacro{\Filt}{{\cal F}}{Filt}{}
\newcommandmacro{\Mod}{\mbox{mod}}{Mod}{}
\newcommandmacro{\Resi}{\mbox{Res }}{Resi}{was without i!}
\newcommandmacro{\Indi}{\mbox{Ind }}{Indi}{was without i!}
\def\RR(#1,#2){R^{#1}_{#2}}   
\def\TT(#1,#2){T^{#1}_{#2}}   

\def\Hom{\mbox{Hom}}

\def\Chi{\chi}

\newcommandmacro{\Ann}{\mbox{Ann}}{Ann}{}
\newcommandmacro{\Cen}{\mbox{Cen}}{Cen}{}
\newcommandmacro{\End}{\mbox{End}}{End}{}
\newcommandabbreviation{\semisimple}{semisimple}{semisimple}
\newcommandabbreviation{\Bratteli}{Bratteli}{Bratteli}
\newcommandabbreviation{\JBC}{Jones Basic Construction}{JBC}
\newcommandabbreviation{\pa}{partition algebra}{pa}
\newcommandabbreviation{\TM}{transfer matrix}{TM}
\newcommandabbreviation{\PM}{Potts model}{PM}
\newcommandabbreviation{\QSC}{quantum spin chain}{QSC}
\newcommandabbreviation{\Hamiltonian}{Hamiltonian}{Hamiltonian}
\newcommandabbreviation{\YS}{Young symmetrizer}{YS}

\newcommand{\ind}{{\mbox{Ind}}}
\newcommand{\res}{{\mbox{Res}}}
\newcommand{\mas}{ \left\{ \begin{array}{c}  }
\newcommand{\sam}{  \end{array}     \right\} }

\newcommand{\cmpx}{CoxMartinParkerXi06}
\newcommand{\dlabringel}{DlabRingel92}
\newcommand{\YounG}{{\mathcal Y}}
\newcommand{\YounGG}{{\mathcal Y}^+}

\newcommand{\CC}{\mathcal{C}}
\newcommand{\edge}[1]{(\phi_{#1})'}
\newcommand{\modules}{\mbox{-mod}}
\newcommand{\lat}{{\mathcal L}}

\newcommand{\Subsection}[1]{\subsection{#1} \ \\}

\begin{document}
\title[Counting Catalan Sets]
{Pascal Arrays: Counting Catalan Sets}
\author{R. J. Marsh}%
\address{Department of Pure Mathematics,
         University of Leeds,
         Leeds,
         LS2 9JT}%
\email{marsh@maths.leeds.ac.uk}%
\author{P. P. Martin}%
\address{Department of Mathematics,
         City University, London,
         Northampton Square,
         London,
         EC1V OHB}%
\email{p.p.martin@city.ac.uk}%
\thanks{2000 {\em Mathematics Subject Classification 05A15, 05C38, (05A18, 17B37, 20C08)}. \\
The first author was supported in part by a Leverhulme Fellowship.}%
\subjclass{}%
\date{December 19, 2006}
\begin{abstract}
Motivated by representation theory
we exhibit an interior structure to Catalan sequences and 
many generalisations thereof.
Certain of these coincide with well-known (but heretofore isolated)
structures. The remainder are new. 

\end{abstract}
\maketitle



{  \parskip=1pt
\tableofcontents }
\section{Catalan Sequences}
\Subsection{Introduction}
The Catalan numbers form a sequence which begins
$$
1,1,2,5,14,42,...,C(n) = \frac{1}{n+1}{\binom{2n}{n}},...,
$$
They occur in a wide variety of distinct combinatorial contexts 
\cite{BloteNightingale82,Fishburn85,Kreweras72,sloane1,Stanley99,Stanley01},
and have many generalisations 
\cite{FominZelevinsky03b,Graham95,Green98,MartinSaleur94a,Matsumoto06}.
In particular, $C(n)$ is equal to the number of non-crossing pair
partitions of $2n$ objects \cite{Kreweras72}. 

The Temperley-Lieb algebra $TL_n$ \cite{TemperleyLieb71}
has a natural basis given by such partitions 
\cite{BloteNightingale82,Martin91}.
Consideration of non-crossing pair partitions
from a Temperley-Lieb representation theoretic
perspective reveals an interior structure
to the Catalan integer sequence,
and an enumeration of each set of non-crossing pair
partitions.
It is interesting to ask whether and how this structure also exists
in the other combinatorial contexts mentioned above. 

The first aim of this paper is to explain this structure
and to give interesting examples. The second aim is to pass this
structure wholesale into generalisations,
both combinatorially and algebraically.
The generalisations appropriate in this context can be considered to be 
classified by arbitrary rooted directed graphs
(with the rooted semi-infinite chain $(A_{\infty},0)$ as the original).

Key motivating examples for our investigation were the Fomin-Zelevinsky
cluster algebras~\cite{FominZelevinsky02}. The cluster algebras of
finite type are classified by the Dynkin diagrams~\cite{FominZelevinsky03b},
and the number of
clusters in a cluster algebra of type $A_{n-1}$ is given by $C(n)$.
We show in particular 
that such clusters possess the structure mentioned above
(see Section~\ref{clusters}). Numerical evidence suggests that clusters of type
$B$ can also be put into this framework. However, as it stands, clusters
of type $D$ do not fit into this picture, suggesting that the
generalisation of the Catalan numbers we consider here (obtained by
considering walks on rooted directed graphs) 
does not fully contain that
arising for cluster algebras~\cite{FominZelevinsky03a}.
On the other hand, our generalisation brings other classical sequences
such as the Bell numbers \cite{Liu68} into the same framework.




 
\Subsection{Representations and Towers of algebras}
The  structure in 
Catalan combinatorics we refer to above consists of two features 
of paths on graphs, which we shall call {\em decomposition} and
{\em edge maps}. 
These will be fully explained in Section~\ref{s:pascalarrays}.
In order to do this, and to explain how they are  
connected to representations of algebras, we need first to recall some
representation theory.
In particular, we recall the formula for the dimension
of a  
finite-dimensional 
algebra in terms of the dimensions of its simple and projective modules.
We also discuss Bratteli diagrams for simple modules in
a tower of algebras with a global limit, and their
close relationship with paths on a related graph.


All our algebras will be finite-dimensional algebras with $1$ over a
field $k$ (except where otherwise stated).
Let $A \subset B$ be an identity-preserving injection of algebras,
with left modules ${}_A M$ and ${}_B N$ respectively. Then 
the left-adjointness of induction to restriction 
\cite{CurtisReiner90}
implies the Frobenius reciprocity:
\eql(FR)
\Hom_B ( {}_B \ind {}_A M, {}_B N) 
\cong \; \Hom_A ( {}_A M , {}_A \res {}_B N  ) . 
\eq
Consider the special case where $A=k$, $M=k$ and $N$ is simple. Then
${}_B \ind {}_A k $ is isomorphic to the left regular module ${}_B B$.
We also have that $\Hom_B ( {}_B B , {}_B N )$ is generated by the
maps from summands of ${}_B B$ with
${}_B N$ in their head --- that is, the copies of the indecomposable
projective module $P_N$ covering ${}_B N$.
By \eqref{FR} we see that there are $\dim \Hom(k,\res N)=\dim {}_B N$
summands $P_N$. Thus
\eql(sumPL)
\dim B \; = \; \sum_{\lambda} \dim P_{L_{\lambda}} \dim L_{\lambda}
\eq
where $\{ L_{\lambda} \}_{\lambda}$ is the set of simple
$B$--modules. 
In particular if $B$ is semisimple then 
\eql(sumLL)
\dim B \; = \; \sum_{\lambda} (\dim L_{\lambda} )^2 . 
\eq
(Combinatorially, 
in case of the group algebra of the symmetric group, 
this is at the heart of the Robinson--Schensted correspondence 
\cite{Knuth98,Schensted61}.)
Thus if a combinatorial set can be equipped with the property of basis
for $B$, it is a matter of representation theory to express its
cardinality as a sum of squares. 
Of course there is no constructive procedure for equipping a set in
this way, in general. 
The aim here is to
present cases in which such an expression can be realised, via a
bijection like the Robinson--Schensted correspondence
(and more specifically via the Robinson--Schensted-like correspondence of
\cite{MartinRollet98}).


We
have not yet used the fact
that combinatorial sets often occur in {\em sequences}.  The
second aspect of representation theory which we want to employ is the
relationship between the algebras in a \emph{tower} of algebras
\cite{GoodmanDelaharpeJones89,CoxMartinParkerXi06} giving these
sequences.

Suppose that $\{A_n\}_{n=0}^{\infty}=A_{\bullet}$ is a tower of algebras (with
identity-preserving inclusions). Let $\Lambda_n$ be an index set for
the simple $A_n$-modules, for each $n$. The \emph{Bratteli diagram} of
$A_{\bullet}$ has vertices given by the simple modules of the $A_n$
arranged in layers indexed by $n$. There are $m$ arrows from
the $A_n$-module $L_n(\mu)$ to the $A_{n+1}$-module $L_{n+1}(\lambda)$
whenever $L_n(\mu)$ appears with multiplicity $m$ as a
composition factor in ${}_{A_{n}}\res_{A_{n+1}}L_{n+1}(\lambda)$.
(Bratteli diagrams appear in \cite{Bratteli72,Effros80}
--- a thorough discussion of the semisimple case may be found in 
\cite{GoodmanDelaharpeJones89}.)

It is clear that the multiplicity of an $A_{n-1}$-module $L_{n-1}({\nu})$ as
a composition factor in ${}_{A_{n-1}}\res_{A_{n+1}}L_{n+1}({\lambda})$ is
given by the number of paths from $L_{n-1}({\nu})$ to $L_{n+1}({\lambda})$.
If $A_0\cong k$, so $\Lambda_0=\{0\}$ for some element $0$, it
follows by induction on $n$ that the dimension of an $A_n$-module
$L_n({\lambda})$ is given by the number of paths from $L_0(0)$ to
$L_{n}(\lambda)$ in the Bratteli diagram.


Suppose now, for a moment, that the following axioms hold:

(i) $\Lambda_{n} \hookrightarrow \Lambda_{n+2}$ for all $n$, so that
the sequences of odd and even index have limits $\Lambda'$ and
$\Lambda$ respectively. That is, for every simple module of $A_{n}$
there is a simple module of $A_{n+2}$ with the same label. 

(ii) There is a graph $G$ (subsequently referred to as the \emph{Rollet graph},
cf. \cite{MartinRollet98}) with vertices $\Lambda\sqcup\Lambda'$ such
that for any $\lambda\in \Lambda_{n+1}$, $\mu\in\Lambda_n$, there is
an arrow from $L_{n+1}(\lambda)$ to $L_{n}(\mu)$ in the Bratteli
diagram if and only if there is an arrow from $\lambda$ to $\mu$ in
$G$, and in this case, the arrows have the same multiplicity, so that
the multiplicity of $L_n(\mu)$ in ${}_{A_n}\res_{A_{n+1}} L_{n+1}(\lambda)$
is given by the number of arrows in $G$ from $\mu$ to $\lambda$.
(Note that there is a natural connection with the concept of 
{\em principal graph} as in \cite[\S4.1]{GoodmanDelaharpeJones89},
from the context of finite dimensional von~Neumann algebras.)


(iii) $A_0 = k$. (Set $\Lambda_0 = \{ 0 \}$.)

An heuristic explanation for the significance of axioms of
this kind in {\em statistical mechanics} and in {\em invariant theory} is given in
Section~\ref{physic};
a more constructive and general axiom set is given in 
Section~\ref{new axioms}.


It  follows that
\prl(f1) 
(i)
The walks on the Rollet graph $G$ from $0$ to $\lambda$ of length $n$ are a
basis of $L_n({\lambda})$.

(ii)
If $A_{n}$ is semisimple, the walks on the Rollet graph $G$ from $0$ to $0$
of length $2n$ are a basis for $A_{n}$ itself. 
\Qed
\end{pr}

We remark that while we have a basis for $L_n(\lambda)$ in (i), the construction does not
tell us the action of $A_n$ on $L_n(\lambda)$.

\smallskip

The idea is that if $A_{n}$ has an interesting combinatorial basis (in
a sense to be elucidated in examples below) then such a
representation-theoretic decomposition will be a way to understand it.

We do not restrict ourselves to consideration of simple modules as in the
above discussion; our more general axiom set allows a much wider collection of
examples. Particularly nice examples arise from towers of quasihereditary
algebras (cf.~\cite{\cmpx}) where the modules considered are the standard
modules over each algebra in the tower, and we will discuss some of these.
It may also be the case that there are nice examples of towers of cellular
algebras, with the modules taken to be the cell modules, 
although we do not
consider that case in this article (other than quasihereditary examples).
(Cellular algebras {\em per se} were introduced in \cite{GrahamLehrer96}.)


\Subsection{Structure of the paper}
The paper is structured as follows. In Section~\ref{s:pascalarrays},
we set up the appropriate
graph-theoretic notation and define the notion of a Catalan
$(G,v_0)$-sequence of sets for a rooted graph $(G,v_0)$. These are sequences
of sets possessing the additional structure referred to above, in the sense
that they arise from an underlying array of sets corresponding to walks on
$G$ starting at $v_0$. In Section~\ref{s:Exa1},
we give a collection of examples of such
sets for the rooted graph $(A_{\infty},0)$, including our motivating example
of basis diagrams of the Temperley-Lieb algebra. Such examples will have
cardinalities given by the sequence $(C(n))$ of Catalan numbers.

In Section~\ref{s:Pascalalgebras} we consider towers of algebras, together
with algebra 
modules at each level. We show that, provided a certain set of axioms
holds, this gives rise to a Bratelli diagram and a Catalan $(G,v_0)$-sequence
for a certain underlying rooted graph $(G,v_0)$. As a way of constructing
examples of towers of algebras satisfying our axioms, we show that a set of
axioms close to that in~\cite{\cmpx} is sufficient to ensure that our axioms
hold. Examples arising in this way will always be quasihereditary, with the
standard modules forming the set of modules at each layer. 
We discuss the example of the Temperley-Lieb algebra in this context.

In Section~\ref{s:blob} we give our first example of a Catalan sequence
corresponding to a rooted graph different from  $(A_{\infty},0)$, 
i.e. blob diagrams. In this case
the corresponding rooted graph is $(A_{\infty}^{\infty},0)$.
In Section~\ref{s:furcated_graphs},
we give a collection of natural examples of Catalan
sequences corresponding to a wide collection of graphs 
--- infinite rooted trees whose
branching properties are governed by the distance to the root. We also
show how these can arise from towers of algebras.

In Section~\ref{s:reentrant}, we consider examples 
for non-tree graphs. Our examples
arise from partition algebras,
Brauer algebras and Hecke algebras. 
The corresponding graphs in our framework include the
double Young graph and the Young graph, respectively, showing that the
underlying graph can be quite complex. The corresponding Catalan sequences
have cardinalities given by the sequences of Bell numbers 
(counting partitions)
and the sequence counting pair partitions. 

In Section~\ref{s:truncation} we indicate how walks on a truncated graph
can arise, by considering certain quotients of the 
`doubly critical' blob algebra \cite{MartinSaleur94a}, and
in Section~\ref{s:dnalgebra} we give an example of a tower of algebras
for the rooted graph $(D_{\infty},0)$ (see Figure~\ref{f:Dinfinity}(b)).
These algebras first appeared in~\cite[\S 7]{Green98}.

\Subsection{Towers in statistical mechanics\label{physic}}
Towers of algebras $A_{\bullet}$ of the kind described above 
occur `in nature' as the transfer matrix algebras for sequences of
lattice statistical mechanical models approaching the thermodynamic limit
(indexed by lattice size --- see \cite{Martin91} for details).
In this setting different simple modules of $A_n$ can be associated to
different correlation functions, that is, to different observables on the
system. Approaching the thermodynamic limit there has to be 
(on physical grounds) a stable notion of observables, essentially
independent of lattice size.  This implies a relationship between
simple modules for different lattice sizes --- i.e. between the simples
responsible for the {\em same observable} in different algebras in
the sequence. At the level of representation theory this is manifest,
in practice, via functors embedding the category of modules for the smaller
system into that of the larger. 
(There are several explicit examples in the literature 
--- see \cite{Martin91,Martin2000} and references therein.)
This is the heuristic explanation for why we expect a non-empty
set of examples satisfying axioms of the kind described above. 
That is, the inclusion $\Lambda_{n-2} \hookrightarrow \Lambda_{n}$
takes the label corresponding to a given observable in level $n-2$ to
the label corresponding to the same observable in level $n$. We remark
that the gap of two here is not forced: see Section~\ref{new axioms}
for consideration of a more general set up.

\Subsection{Geometrical considerations in representation theory}
The endowment of the index set for simple modules of an algebra
with a geometrical structure
(as in the weight space and weight lattice for Lie algebras) has been very
useful in representation theory, {\em where it is possible}. 
For example the alcove geometry which
describes the representation theory of quantum groups starts with an
embedding of the index set in a space, on which a reflection group
then acts.
The number of instances of this {\em precise } abstract setting is manifestly
limited \cite{Humphreys90}. 
Thus it is of interest to consider generalisations. 
A very mild generalisation of the classical use of alcove geometry
\cite{Jantzen87} is discussed in \cite{MartinWoodcock03}, 
and another in \cite{MartinWoodcock2000}, 
but there are many towers of algebras known for which,
if geometry is to be applied, a more radical generalisation will be needed. 
Interpolating between the two examples above, it can be seen that
while there may be no obvious metric space in which to embed a given
index set, there may be a natural underlying graph. 
From a representation theory perspective this paper may be seen in part as 
an exploration of the uses of this idea.


\section{Pascal arrays and Catalan sequences} \label{s:pascalarrays}
\Subsection{Paths on Graphs}
A {\em graph} $G$ is defined to be a set of vertices $V=V_G$
together with a set of edges $E=E_G$ and a pair of maps $i$ and $f$ 
(initial and final) from $E$ to $V$. 
Thus an element $e \in E$ is a directed edge from $i(e)$ to $f(e)$.
\\
If $G$ is a graph then the \emph{opposite graph} to $G$ is the graph
$G^{op}$ obtained by interchanging the functions $i$ and $f$. 

We recall that a graph is said to be \emph{loop-free} if it
contains no edge $e$ such that $i(e)=f(e)$; it is said to be
\emph{simply-laced} if it contains no pair $e,e'$ of edges such that
$i(e)=i(e')$ and $f(e)=f(e')$, and it is said to be
\emph{undirected} if for each edge $e$ there is an edge $e'$ such that
$i(e')=f(e)$ and $f(e')=i(e)$. 

Thus a simply-laced undirected graph $G$ may be defined to be a set of
vertices $V$ together with a set of edges $E$ consisting of subsets of $V$ of
cardinality $2$. 

A {\em rooted graph} $(G,v_0)$ is a graph $G$ together
with a distinguished vertex $v_0$.
If $v$ is a vertex in a graph, its adjacent vertices are defined to
be those linked to $v$ by an edge starting at $v$. 
We assume the {\em valency} of each
vertex (the number of adjacent vertices) to be finite.

We use Dynkin diagram naming conventions for
appropriate simply laced graphs 
--- see, for example, Figures~\ref{f:Ainfinity}(a) and~\ref{f:Dinfinity}(b).
\begin{figure}[htbp]
$$
\hspace{.5in}
\xy
\xymatrix{
0 \ar@{-}[r] & 1 \ar@{-}[r] & 2 \ar@{-}[r] & 3 \ar@{--}[r]
& {\ }
}
\endxy
\hspace{1.32in}
\xy
\xymatrix{
0 \ar@{-}[dr] \\
& 1 \ar@{-}[r] & 2 \ar@{-}[r] & 3 \ar@{--}[r] & {\ } \\
0' \ar@{-}[ur]
}
\endxy
$$
\caption{(a) The graph $A_{\infty}$. \hspace{1.5in} (b) The graph $D_{\infty}$}
\label{f:Ainfinity}
\label{f:Dinfinity}
\end{figure}


A {\em path} in graph $G$ is defined to be a sequence $\mathbf{e}=
(e_1,e_2,\ldots ,e_l)$ of edges of $G$ such that
$f(e_i)=i(e_{i+1})$ for $i=1,2,\ldots ,l$. 
The {\em length} of $\mathbf{e}$
is defined to be the integer $l\geq 0$.
For $v,w\in V$, let $P_G(n;v,w)$ denote the set of paths in $G$ starting at
$v$ and ending at $w$ with length $n$. If $v=v_0$ (the distinguished
vertex) then we may omit it. We denote the cardinality of
$P_G(n;v,w)$ by $N_{G,v}(n;w)$, and drop the subscripts $G,v$ if they are
not needed.

\del(def array)
For fixed $v\in V$, we regard the collection of sets $P_G(n;v,w)$ as an array
$Y_{G,v}=(Y_{G,v}(n;w))_{n\in\mathbb{N},w\in V}$ indexed by the length
$n$ and the vertex $w$. We regard
$(Y_{G,v}(n;w))_{w \in V}$
for fixed $n$ as forming the $n$th {\em layer} of the array.
\end{de}

\begin{figure}[htbp]
\[
\includegraphics[width=13.1cm]{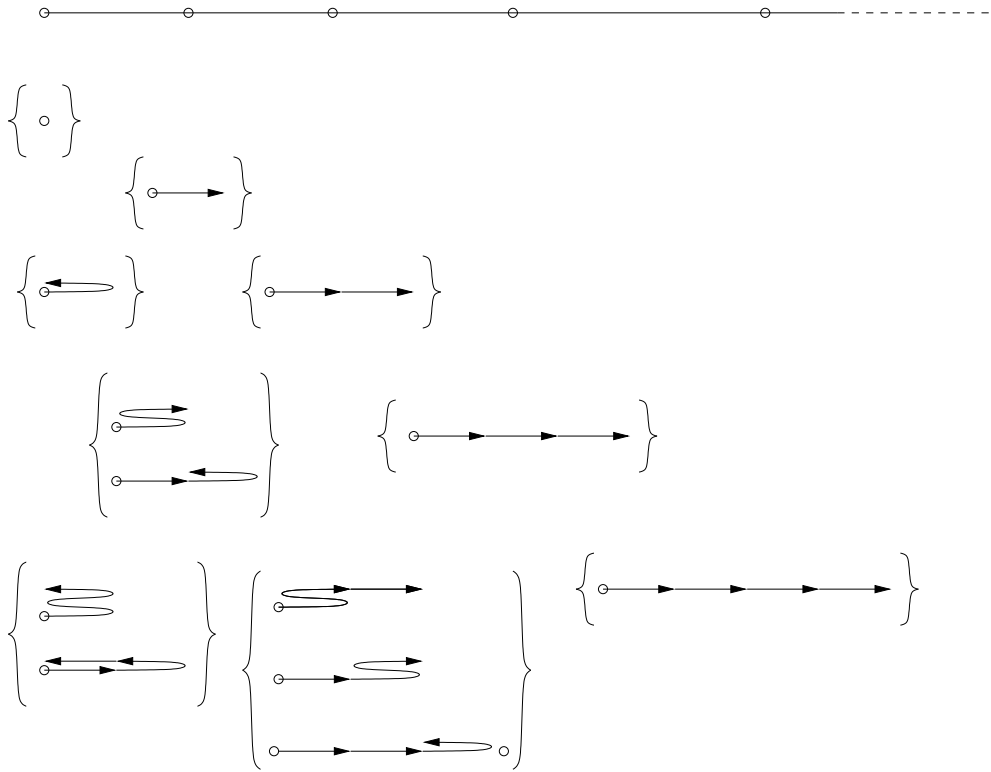}
\]
\caption{The array $Y_{A_{\infty},0}$ of sets of paths
from 0 on $A_{\infty}$, truncated below.  
The accompanying
copy of the graph is positioned so as 
to associate vertices to the columns in the array.
The column position in the array thus 
corresponds to the walk endpoint.}
\label{f:Aarray}
\end{figure}

\begin{figure}
\[
\includegraphics{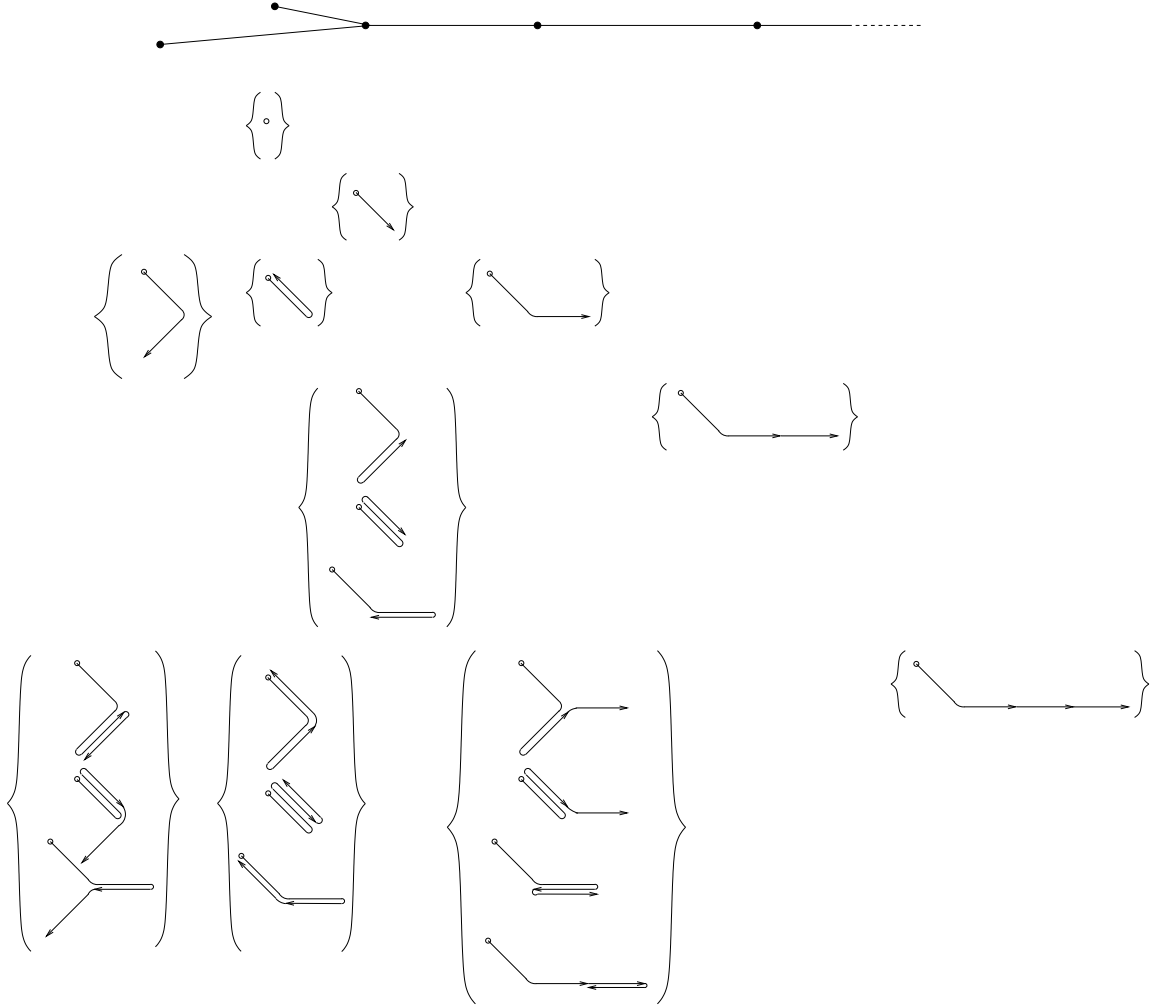}
\]
\caption{\label{D-seq walk}
The Pascal array 
$Y_{D_{\infty},0}$.
As before we include a copy of the graph to label the columns.}
\end{figure}

\begin{example} \label{e:Ainfinity} \rm
Let $G$ be the graph $A_{\infty}$ or  $D_{\infty}$
--- see Figure~\ref{f:Ainfinity}. 
We consider the rooted graph $(G,0)$.

In this case the top of the array of paths is given in
Figure~\ref{f:Aarray}
(respectively Figure~\ref{D-seq walk}).

Note that paths in $A_{\infty}$ are displayed as finite graphs which are
to be read from bottom to top; thus (for example) the graph in the bottom
left of the picture corresponds to the path which visits the vertices
$0,1,2,1$ and then $0$ (in that order).
\end{example}


\Subsection{Arrays of paths and edge maps}
Let $(G,v_0)$ be a rooted graph, with vertices $V$ and edges $E$.
For $e \in E$ let 
$$
\phi_{e} : Y_{G,v_0}(n-1;i(e)) \hookrightarrow Y_{G,v_0}(n;f(e))
$$
be the natural embedding of the shorter path in the longer one,
i.e. $\phi_e(e_1,e_2,\ldots ,e_l)=(e_1,e_2,\ldots ,e_l,e)$.
We call these maps \emph{edge maps}.
Then by construction

\prl(edgemap1)
\begin{equation} \label{inductionstep}
Y_{G,v_0}(n;v) = 
\bigsqcup_{ e \in E ; \; f(e)=v} 
   \phi_{e} ( Y_{G,v_0}(n-1;i(e))  ) 
\end{equation}
And
 any total order on $E$ passes to a total order
on the union of sets  in the $n^{th}$ layer of $Y_{G,v_0}$ 
(i.e. $\cup_{v\in V} Y_{G,v_0}(n;v)$) 
via the lexicographic ordering of sequences on $E$.
\Qed
\end{pr}


\del(p:array)
An array of sets $Y$ is said to be {\em equivalent} to $Y_{G,v_0}=
(Y_{G,v_0}(n;w))_{n\in\mathbb{N},w\in V}$,
and called a \emph{Pascal $(G,v_0)-$array},
if it has 
\newline
(1) the same index sets ($V$ and $\mathbb{N}$), 
\newline
(2) pointwise bijections
between the sets  
in layer $n=0$ of $Y$ and of $Y_{G,v_0}$,
and 
\newline
(3)
a set of inclusions corresponding to the edge maps $\phi_{e}$ which
satisfy equalities corresponding to equation~\eqref{inductionstep},
i.e. for each edge $e$, a map
$$
\edge{e} : Y(n-1;i(e)) \hookrightarrow Y(n;f(e))
$$
satisfying:
\begin{equation} \label{inductionstepgeneral}
Y(n;v) = 
\bigsqcup_{ e \in E ; \; f(e)=v} 
   \edge{e} ( Y(n-1;i(e))  ).
\end{equation}
\end{de}
We will refer to the maps $\edge{e}$ in (3) as edge maps also.

\prl(p:equivalence)
Let $Y$ be an array of sets equivalent to $Y_{G,v_0}$. Then the entire
array of sets is in pointwise explicit bijection with $Y_{G,v_0}$.
\end{pr}
{\bf Proof:} 
This follows from requirement (3) of Definition \ref{p:array}
together with the fact (requirement (2)) that the sets in layer $0$ all
have cardinality $0$ or $1$.
Each path in  $Y_{G,v_0}$ is given by a sequence of edge maps.
Its image is given by the sequence of corresponding maps from (3). 
$\Box$


{\bf Example:} For any undirected graph $G$ 
the array 
$\{ P_{G^{op}}(n;w,v_0) \}_{n\in\mathbb{N},w\in V}$ is equivalent to $Y_{G,v_0}$.


Suppose that $Y$ is a Pascal $(G,v_0)$-array. 
Let $N(n;v)$
denote the cardinality of $Y(n;v)$ for all $n\in\mathbb{N}$ and
$v\in V$. Note that necessarily we have $N(n;v)=|P_G(n;v_0,v)|$ for
all $n\in\mathbb{N}$ and $v\in V$. We have an array of cardinalities
corresponding to the array $Y$ of sets satisfying:
\begin{eqnarray} \label{inductionstepcardinality}
N(0;v) & = & 1 \\
N(n;v) & = & 
\sum_{x \in G, \; e \in E} 
N(n-1;i(e)).   \nonumber
\end{eqnarray}

\begin{figure}[htbp]
\includegraphics{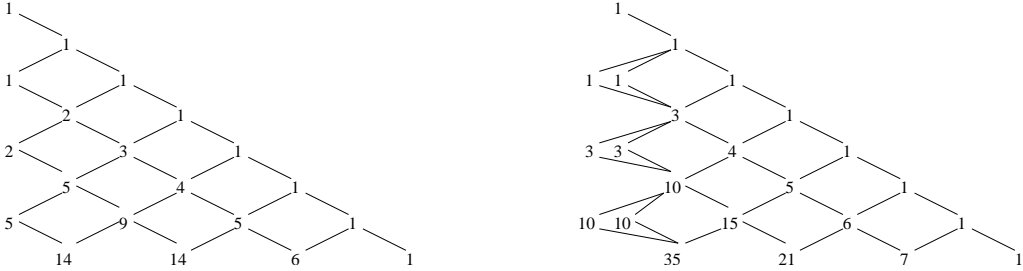}
\caption{The Catalan triangle (cardinality array for
  Figure~\ref{f:Aarray});
and the cardinality array for Figure~\ref{D-seq walk}.}
\label{f:catalantriangle}
\end{figure}

{\bf Example:}
Let $G$ be the graph $A_{\infty}$ considered in the previous section. Then
equations~\ref{inductionstepcardinality} become:
\begin{eqnarray*}
N(0;v)&=&1, \\
N(n;v)&=&N(n-1;v-1)+N(n-1;v+1),
\end{eqnarray*}
where we define $N(n;-1)=0$ for all $n\in\mathbb{N}$.
The solution to this recurrence is the array given by the Catalan triangle,
the top of which is shown in Figure~\ref{f:catalantriangle};
see for example~\cite[page 796]{abramowitzstegun1}.

{{\rem
The above combinatorics can be encoded using the
adjacency matrix $M(G)=
(M_{v,w})_{v,w\in V}$ of the graph $G$, where $M_{v,w}$ is the number of
edges $e$ with $i(e)=v$ and $f(e)=w$. The $v_0,w$ entry in $M(G)^n$ is
equal to $N_{G,v_0}(n;w)$. 
}}
{{\rem
The entries in  
the Catalan triangle are also sometimes known as generalised
Catalan numbers \cite{Ng06}. 
Note that, for $n-v$ even,
$$N(n,v)=\binom{n}{\frac{1}{2}(n-v)}-
\binom{n}{\frac{1}{2}(n-v)-1},$$
with the convention that
$\binom{n}{-1}=0$. (See e.g.~\cite{Guy00,Ng06}).
}}



\Subsection{Closed path decomposition and Catalan sequences}
Let $*:P_G(n;v,v')\times P_G(m;v',v'')\rightarrow P_G(n+m;v',v'')$
denote the composition of paths.
Note that the triple ($V_G$, $\cup_n P_G(n;-,-)$, $*$), is a category
with units given by paths of length zero (it is known as the
\emph{free category} on $G$).
The following is immediate. 

\prl(decomp1)
In general, we can express the set of paths from $v_0$ back to $v_0$ 
in $2n$ steps in the following way:
\begin{equation} 
\label{N=sumN2}
P_G(2n;v_0,v_0)= * ( \cup_{v\in G} P_G(n;v_0,v) \times P_G(n;v,v_0) )
\end{equation}
That is, $*$ defines an bijection between
from $\cup_{v\in G} P_G(n;v_0,v) \times P_G(n;v,v_0)$ to
$P_G(2n;v_0,v_0)$.\Qed
\end{pr}
In this setting the inverse {\em decomposition} map to $*$ is obvious.
In this paper we will be concerned with sets in bijection with 
$P_G(2n;v_0,v_0)$ (for various choices of $(G,v_0)$) for which the 
decomposition map is not so obvious (even the target set is
unobvious),
but which reveals important features of these sets (and associated
algebras). 


If $G$ is undirected, it follows from Proposition~\ref{decomp1} that
\begin{equation} \label{e:sumofsquares}
N_{G,v_0}(2n;v_0)=\sum_{x\in V}N_{G,v_0}(n;x)^2.
\end{equation}

We call the sequence $ \{ N(2n;v_0) \}_{n}$ the 
{\em Catalan  $(G,v_0)$-numbers}.

\del(d: cat seq)
A \emph{Catalan $(G,v_0)$-sequence} $\CC$ is any sequence of sets
$(\CC(n))_{n\in\mathbb{N}}$
in explicit natural bijection $\Phi$ with
$\{ \cup_{v \in G} Y(n;v) \times  Y'(n;v) \}_{n\in\mathbb{N}}$
where $Y$ is a Pascal $(G,v_0)$-array and $Y'$ is a Pascal
$(G^{op},v_0)$-array.
In this case we say that the Pascal array $Y$ {\em underlies} the Catalan
sequence (via $\Phi$). 
\end{de}

We remark that if $G$ is undirected then $Y'$ can be taken to be a
Pascal $(G,v_0)$-array in the above: we make much use of this in the
sequel. In this case, the sequence of 
cardinalities of such a sequence is given by the Catalan $(G,v_0)$-numbers.

We will see later that, beside the many well-known Catalan sequences
(of `type-A') this formalism brings many entirely distinct
combinatorial sequences, such as the sets of partitions of sets of
degree $n$, into the same framework. 

Note that 
\prl(p:C equivalence)
The pointwise bijections from Proposition~\ref{p:equivalence}
together with the bijections $\Phi$ give an explicit pointwise
bijection between any two Catalan $(G,v_0)$-sequences. \Qed
\end{pr}

The reason for this nomenclature is the following example. 

{\bf Example:}
For $(G,v_0)=(A_{\infty},0)$ we obtain that
$$
N_{A_{\infty},0}(2n;v_0)=\sum_{v=2n,2n-2,\ldots ,1\text{ or 0}}
\left( \binom{n}{\frac{1}{2}(n-v)}-
\binom{n}{\frac{1}{2}(n-v)-1}\right)^2=C(n).
$$
That is, $\{ P_{A_{\infty},0}(2n;0,0) \}_n $  is the
{\em ordinary} Catalan sequence.
Equation~\ref{e:sumofsquares} states
that the $n$-th term is the sum of the squares
in layer $n$ of the triangle in Figure~\ref{f:catalantriangle}.
This corresponds 
via equation~\ref{sumLL} 
to a well-known result from representation theory
in this case \cite{Martin91}. We will review this shortly
(in Section~\ref{TLA})
before going
on to consider a number of generalisations.  


\section{Examples of $(A_{\infty},0)$-arrays} \label{s:Exa1}

In this section we recall some important Catalan sequences.
By discovering decomposition maps for these (in the sense of
Proposition~\ref{decomp1}), and edge maps for the target sets 
(in the sense of Proposition~\ref{edgemap1} and Definition~\ref{p:array}),
we give the corresponding Pascal $G$-arrays
with $G$ the rooted graph $(A_{\infty},0)$.
Later we will give
examples of the same kind for other rooted graphs.

We start with the core example from representation theory.

\Subsection{Temperley-Lieb diagrams\label{s:TLdiagrams}}
Recall that a pair partition of $2n$ objects is a partition of the objects
into parts each containing precisely two elements.
One may draw a picture of such a pair partition by placing the objects
on some horizontal line in the plane, and drawing lines between them
in the plane below this horizontal. 
A pair partition is {\em non-crossing} if 
 each pair may be connected by a line drawn on the
lower half plane
{\em simultaneously} with no lines crossing.

A famous Catalan sequence~\cite{Stanley99} has $n^{th}$ term
given by the set of non-crossing pair partitions of $2n$ objects.
The $n=3$ cases may be represented as follows. 

\includegraphics{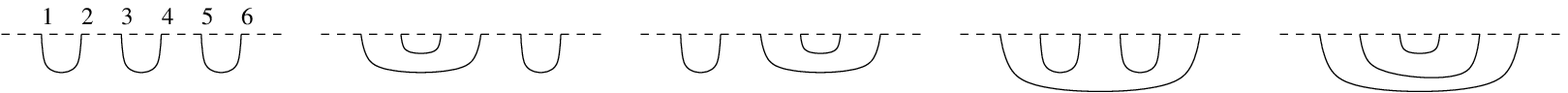}


The question here is: Is there an underlying Pascal $G$-array?
This is the same as to ask:
Is there a decomposition corresponding to the path decomposition
and a set of edge maps?
(The answer is yes. And in this case we should understand it in terms of
representation theory.)

Note that the distinction between the 
embedding in the 
lower half plane and any other
plane interval, such as a disk, is immaterial in this 
non-crossing pair partition
construction. 
(The non-crossing property is preserved by conformal transformations,
for example.)
Let us think now of a representative non-crossing line
realisation of such a pair partition as indeed being drawn on a disk. 

Suppose we number the objects $1,2,...$ clockwise from some chosen
point.
The set $\TLD(m,2n-m)$ of {\em ordinary (\TL) diagrams} is simply 
a representative set of these realisations with a
notional aggregation of the objects into `northern' subset
$\{1,2,..,m\}$ of adjacent objects and a southern subset $\{m+1,..,2n \}$
of adjacent objects. 
Often the disk is also replaced by a rectangle, with the northern
(resp. southern) subset lying on the northern (resp. southern) edge.
\[
\includegraphics{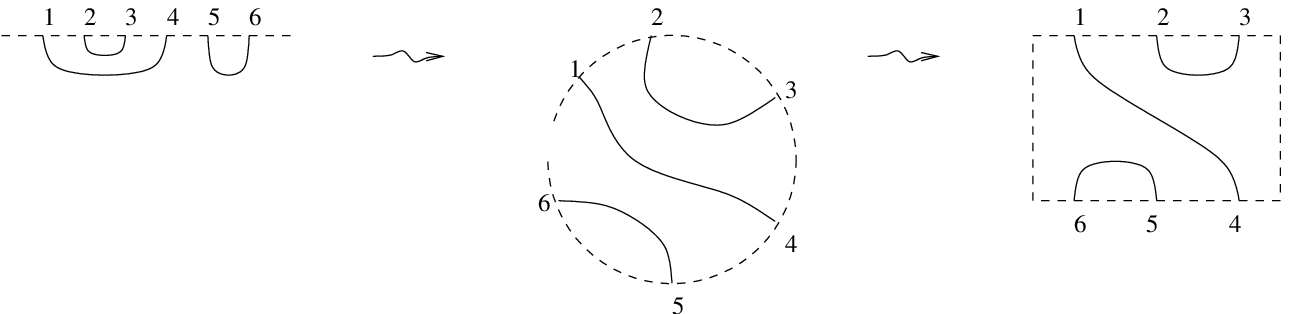}
\]


Note that for $n>0$ (and $m+n$ even) there is a natural bijection 
\[
\phi^u : \TLD(m,n) \rightarrow \TLD(m+1,n-1)
\]
via the common underlying set.

The {\em propagating} lines in a \TL\ diagram are the
lines {\em between} the northern and southern subsets. Note that
the number of these is {\em not} fixed by $\phi^u$. 

Let $\TLD_l(m,n)$ denote the subset of $\TLD(m,n)$
consisting of elements with $l$ propagating lines.

\begin{figure}[htbp]
\[
\includegraphics{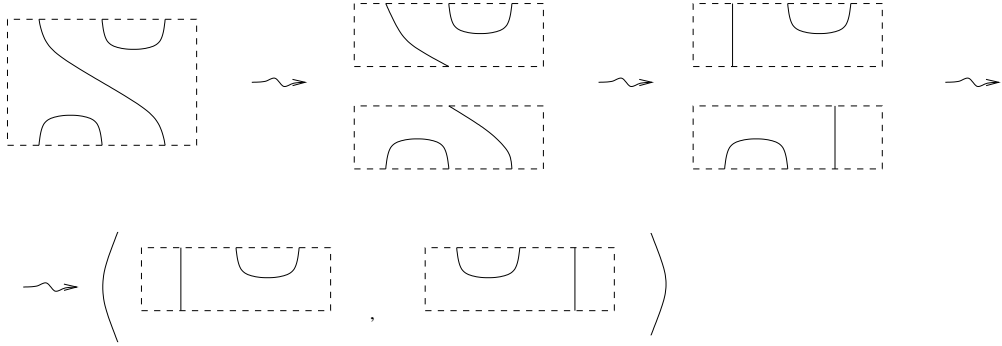}
\]
\caption{Splitting a $D(3,3)$ diagram into two. 
\label{f:TLhalfdiagram}}
\end{figure}

The basic operations on \TL\ diagrams we need to consider are:

(1) Bra-ket decomposition: It is always possible to cut a diagram from east
to west such that each propagating line is cut once, and no other line
is cut. This process produces a well defined pair of half--diagrams,
which itself may be thought of 
(in case $m=n$)
as lying in $\TLD(n,l) \times \TLD(l,n)$ for some $l$.
Indeed
\begin{equation}
\label{ppm01}
\TLD_l(n,n) \cong \TLD_l(n,l) \times \TLD_l(l,n)
\end{equation}
where the map from left to right 
is given by cutting the propogating lines; the inverse map 
is to join the propagating lines in the obvious order.
See Figure~\ref{f:TLhalfdiagram} for an example.


{\bf Remark.} The reason for using Dirac's bra-ket terminology will become
apparent in Section~\ref{TLA}. 


(2) Edge maps:

(2a) Given an upper half diagram in $\TLD_l(n,l)$ as above we may
construct one in $\TLD_{l+1}(n+1,l+1)$ by adding a propagating line on the
right. We denote this map by $\phi^1$. For example:
\[
\phi^1: \raisebox{-22pt}{\includegraphics{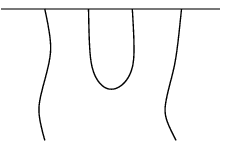}} 
\ \ \ \mapsto\ \ \ \raisebox{-22pt}{\includegraphics{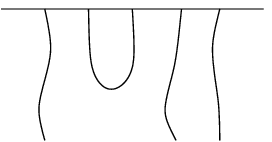}}
\]
(2b) Given an upper half diagram in $\TLD_l(n,l)$ with $l>0$ we may
   construct another 
half diagram
via the natural injection $\phi^u$ 
into $\TLD_{l-1}(n+1,l-1)$, i.e.
``bending over'' the rightmost propagating line onto a new vertex on the
horizontal line. 
For example:
\[
\phi^u : \raisebox{-22pt}{\includegraphics{./xfig/halfdi1}}
\ \ \ \mapsto\ \ \ \raisebox{-22pt}{\includegraphics{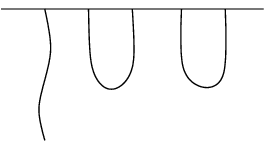}}
\]


\begin{prop}
The array of sets $Y_{TL}=((\TLD_l(n,l))_{l\in A_{\infty}})_n$ 
(see Figure~\ref{f:TL1})
is a Pascal $(A_{\infty},0)$--array with edge maps 
$\edge{i,i+1}$
given by
$\phi^1$ and edge maps $\edge{i+1,i}$ given by $\phi^u$.
This array underlies the sequence $( D(n,n) )_n $
via equation~(\ref{ppm01}). 
\end{prop}

{\bf Proof:} For the first part it suffices to note that 
$$
\TLD_l(n,l) 
\; = \; \phi^u( \TLD_{l+1}(n-1,l+1) ) \; \bigsqcup \; \phi^1(
\TLD_{l-1}(n-1,l-1) ).
$$

For the second part, note that
an 
explicit construction (and {\em en passant} a counting) of the Catalan 
sets in the sequence $(\TLD(n,n))_n$ follows from equation (\ref{ppm01}),
and from the fact that $\TLD(n,n) = \bigcup_l \TLD_l(n,n)$. We obtain:
\eql(ep1)
\TLD(n,n) = \bigsqcup_l \TLD_l(n,l) \times D_l(l,n).
\eq
$\Box$

Note that in Figure~\ref{f:TL1} we have drawn edges indicating when the
edge maps are being applied (all such edges should be regarded as being
oriented down the page). We sometimes use such edges for clarity;
in such cases the brackets separating out the sets in the array are no
longer necessary.


{\oddsidemargin -1cm
\newcommand{\aar}{\ar@{-}}         
\newcommand{\headroom}{80}       
\newcommand{\raised}{-\headroom} 
\newcommand{\pp}[2]{\begin{picture}(#1,\headroom)(0,\raised)
    \put(-20,0){#2} \end{picture}}
\newcommand{\ppp}[4]{\begin{picture}(#1,#2)(0,-#3) \put(-20,0){#4}
    \end{picture}}       

\newcommand{\putch}{\put(-.3,-26)}
\newcommand{\poutch}{\put(-.3,-20)}

\begin{figure}
\xymatrix%
@R=35pt@C=15pt@M=0pt%
{
{\begin{picture}(0,10)(10,-5) 
\put(0,-40){ {\Large $\emptyset$ } 
} 
\end{picture}}
\\
\aar[dr]+U
\\
 & 
\ppp{.2}{0}{55}{
{\put(12,-60){
\put(0, 0){\dashbox{1.0}(20,0){}}
\thicklines
\put(10, 0){\line(0,-1){20}}
}}
} 
\\
&\aar[dl]+U_{\phi^u} 
\aar[drr]+U^{\phi^1}
& 
\\
\ppp{15}{0}{45}{
{\put(20,-55){
\put(0, 0){\dashbox{1.0}(30,0){}}
\thicklines
\put(15.0, 0){\oval(10,10)[b]}
}}
}
&&&
\ppp{20}{0}{45}{
{\put(20,-55){
\put(0, 0){\dashbox{1.0}(30,0){}}
\thicklines
\put(10, 0){\line(0,-1){20}}
\put(20, 0){\line(0,-1){20}}
}}
}
\\
\aar[dr]+U &&&
\aar[dll]+U\aar[drrr]+U
\\
&
\ppp{20}{30}{20}{
{\put(20,0){
\put(0, 0){\dashbox{1.0}(40,0){}}
\thicklines
\put(15.0, 0){\oval(10,10)[b]}
\put(30, 0){\line(0,-1){20}}
}}
\putch{
{\put(20,0){
\put(0, 0){\dashbox{1.0}(40,0){}}
\thicklines
\put(10, 0){\line(0,-1){20}}
\put(25.0, 0){\oval(10,10)[b]}
}}}
}
\newcommand{\pppshift}[1]{{\ppp{20}{.3}{-15}{##1}}}
&&&&&
\ppp{20}{23}{5}{
{\put(20,0){
\put(0, 0){\dashbox{1.0}(40,0){}}
\thicklines
\put(10, 0){\line(0,-1){20}}
\put(20, 0){\line(0,-1){20}}
\put(30, 0){\line(0,-1){20}}
}}
}
\\
&\aar[dl]+U\aar[drr]+U
&&&&&\aar[dlll]+U\aar[drr]+U
\\
\ppp{50}{100}{90}{
{\put(20,0){
\put(0, 0){\dashbox{1.0}(50,0){}}
\thicklines
\put(15.0, 0){\oval(10,10)[b]}
\put(35.0, 0){\oval(10,10)[b]}
}}
\poutch{
{\put(20,0){
\put(0, 0){\dashbox{1.0}(50,0){}}
\thicklines
\put(25.0, 0){\oval(30,18.0)[b]}
\put(25.0, 0){\oval(10,10)[b]}
}}
}}
&&&
\ppp{50}{100}{90}{
{\put(20,0){
\put(0, 0){\dashbox{1.0}(50,0){}}
\thicklines
\put(15.0, 0){\oval(10,10)[b]}
\put(30, 0){\line(0,-1){20}}
\put(40, 0){\line(0,-1){20}}
}}
\putch{
{\put(20,0){
\put(0, 0){\dashbox{1.0}(50,0){}}
\thicklines
\put(10, 0){\line(0,-1){20}}
\put(25.0, 0){\oval(10,10)[b]}
\put(40, 0){\line(0,-1){20}}
}}
\putch{
{\put(20,0){
\put(0, 0){\dashbox{1.0}(50,0){}}
\thicklines
\put(10, 0){\line(0,-1){20}}
\put(20, 0){\line(0,-1){20}}
\put(35.0, 0){\oval(10,10)[b]}
}}
}}}
&&&&&
\ppp{20}{100}{85}{
{\put(20,0){
\put(0, 0){\dashbox{1.0}(50,0){}}
\thicklines
\put(10, 0){\line(0,-1){20}}
\put(20, 0){\line(0,-1){20}}
\put(30, 0){\line(0,-1){20}}
\put(40, 0){\line(0,-1){20}}
}}
}
}
\caption{\label{f:TL1} The start of the Pascal array $Y_{TL}$}
\end{figure}
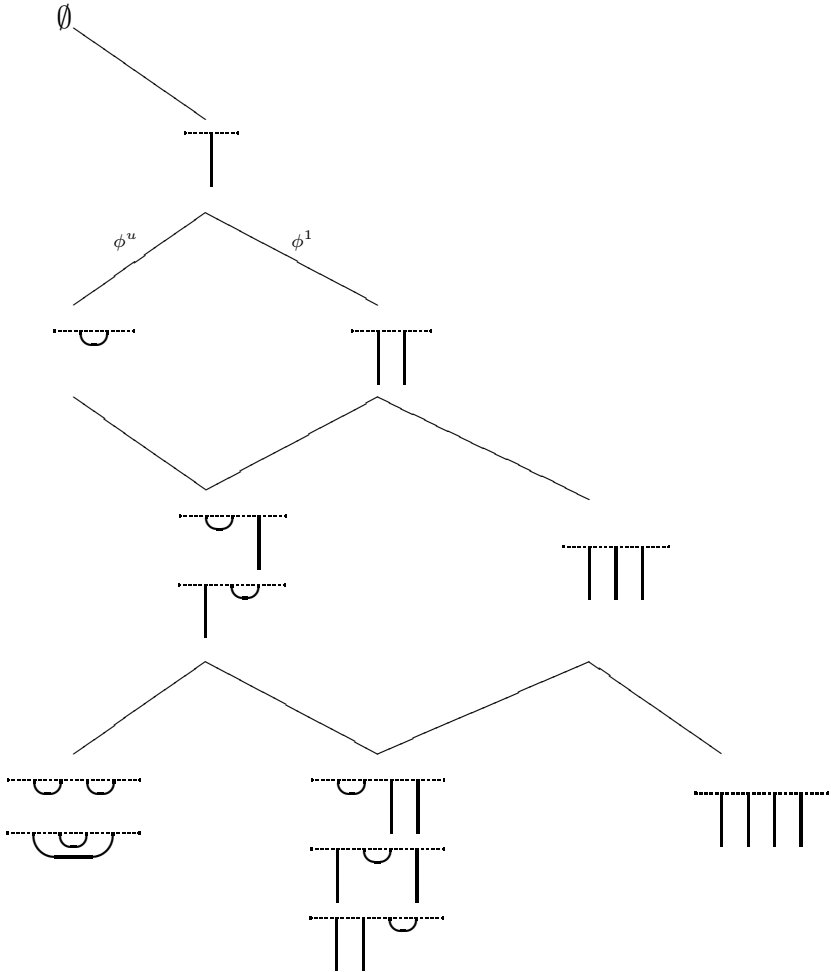
}


We shall see in Section~\ref{TLA} that $\TLD(n,n)$ is a basis for 
the Temperley--Lieb algebra $TL_n(q)$
(independently of a choice of groundfield $k$ and parameter $q$), 
and that $\TLD_l(n,l)$ is a basis for a special kind of $TL_n$-module
denoted $\Delta_n(l)$ (i.e. a standard module, see e.g.~\cite{\dlabringel}).

On the one hand (\ref{ep1}) is (via the array equivalence) a special case of 
equation (\ref{N=sumN2}), and on the other hand it implies that
\begin{equation} \label{dimTL}
\dim(TL_n) = \sum_{\lambda} (\dim(\Delta_n(l)))^2.
\end{equation}
Since for generic $q\in k=\mathbb{C}$,
$\Delta_n(l)$ is isomorphic to the corresponding simple module
$L_n(l)$, and $TL_n$ is semsimple (see \cite{Martin91}), 
this provides a example of the deep
algebraic result 
in equation~\ref{sumLL}. 

\Subsection{Parentheses\label{ss:Cb}}
Consider the set $\CC_b(n)$ consisting of matched bracket sequences containing
$n$ pairs of brackets, all of the same type. For example:
\[ 
\CC_b(3) = \{ ()()(), (())(), ()(()) , (()()), ((())) \}.
\]

Let us represent a walk in $Y_{A_{\infty},0}$ 
(an element in one of the sets in Figure~\ref{f:Aarray})
by a sequence constructed from the symbols $\{1,2\}$,
where a 1 means step along an edge from $i$ to $i+1$ (some $i$)
and a 2 means step from $i+1$ to $i$. 
Thus all paths start with a 1, and the running total of 2's in any
truncation never exceeds the running total of 1's 
(we also call this a {\em standard} sequence
\cite{StantonWhite86}). 
If we simply replace 1 with ``('' and 2 with ``)'' we see that these
sequences pass to properly nested, but not necessarily closed bracket
sequences. That is, no close-bracket appears which does not close an
open-bracket.
Write $Y_b(n;l)$ for the set of sequences of $n$ individual brackets
in which there are $l$ more open-brackets than close-brackets.
The Pascal array corresponding to $\CC_b$ comes from the following
operations:
\\
(1)
Decomposition: An element $x \in \CC_b(n)$ has $2n$ individual
brackets. Break this sequence into two parts, each of $n$ brackets.
The left-hand sequence is of the form of an element of $Y_b(n,l)$ for
some $l$. The right-hand sequence is not of such form, but the reverse
sequence is. 
\\
(2) 
Edge maps:  
For $i\geq 0$, the map
$\edge{i,i+1}$ is given by adding an open bracket and the map
$\edge{i+1,i}$ is given by adding a close bracket.

\begin{figure}
$$
\xymatrix@R=5pt@C=5pt{ \{ \emptyset \} \\ 
 & { \mas ( \sam } 
 \\
{ \mas ()  \sam }     && { \mas (( \sam }
 \\
 &  { \mas ()( \\ (()  \sam  }    && { \mas ((( \sam }  
 \\ 
{ \mas ()() \\ (()) \sam }
           && {\mas ()(( \\ (()( \\ ((()  \sam} && {\mas  (((( \sam}
\\
 & {\ldots}
}
$$
\caption{\label{Abra} Pacal array $Y_b$. Here $\emptyset$ denotes the
empty sequence.}
\end{figure}
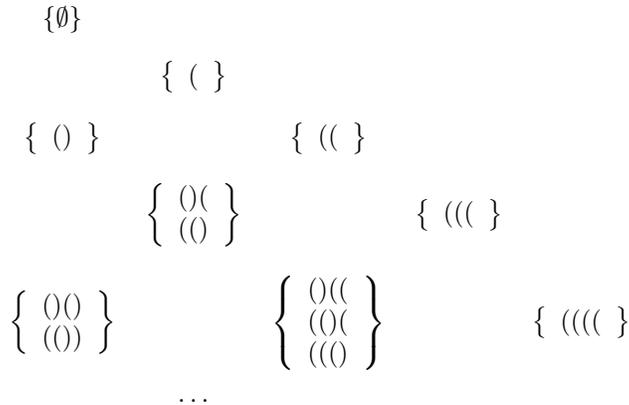

Thus we have: 
\prl(p:Y_b)
The array $Y_b$  (Figure~\ref{Abra})
is a Pascal $(A_{\infty},0)$--array  
(cf. Figure~\ref{f:Aarray}) 
underlying the sequence $\CC_b$.
\Qed
\end{pr}

\Subsection{Rooted planar trees\label{ss:Ct}}
Another famous Catalan sequence \cite{Stanley99}
has $n^{th}$ set the set $\CC_t(n)$ of
rooted planar trees with $n$ edges.
For example, $\CC_t(3)$ is given as follows:
\[
\CC_t(3) : \;\;\;\qquad 
\raisebox{-.41in}{
\includegraphics[width=8cm]{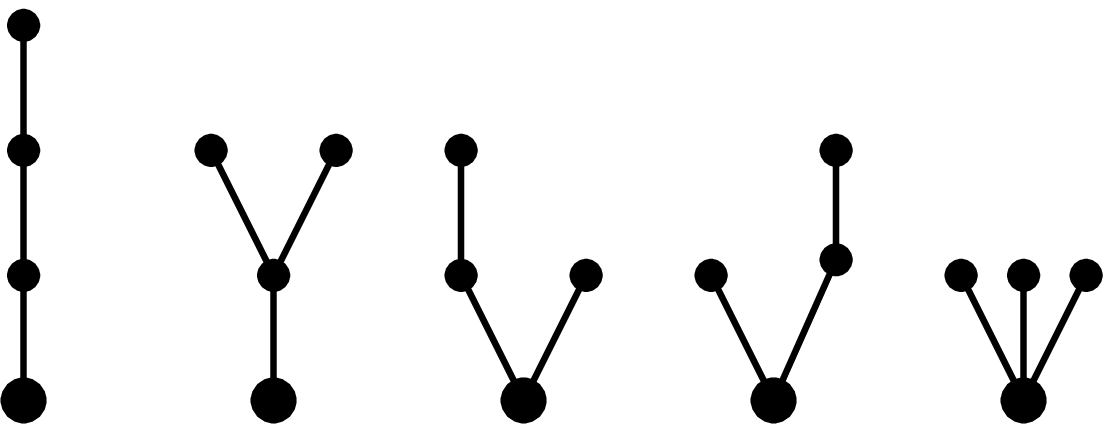}
}
\]
More precisely, a {\em concrete} rooted planar tree is a tree 
(a connected acyclic graph of vertices with at most single edges
between them), together with an embedding in the upper half-plane 
that takes the root vertex (and no other vertex) to the half-plane boundary,
and has straight edges and no edge crossings. 
(Note that every tree has such an embedding.)
Two such embedded trees are {\em equivalent} if 
they are related by an isotopy (note that this 
excludes moves which pass one branch over another). 
A rooted planar tree is an equivalence class of such embedded trees.


Algorithmically,  
two such embedded trees are distinct if they are distinguished by the following
procedure. 
\newline
Consider the planar figure constructed by drawing around the
outside of the embedded tree, 
starting from the root and drawing (say) anticlockwise 
(cf. the contour traversals in \cite{ChassaingDurhuus06}). 
On completion this figure follows the shape of the tree, except that
each edge is replaced by an outward/return pair of edges: 

\[
\raisebox{-.01in}{
\includegraphics[width=7cm]{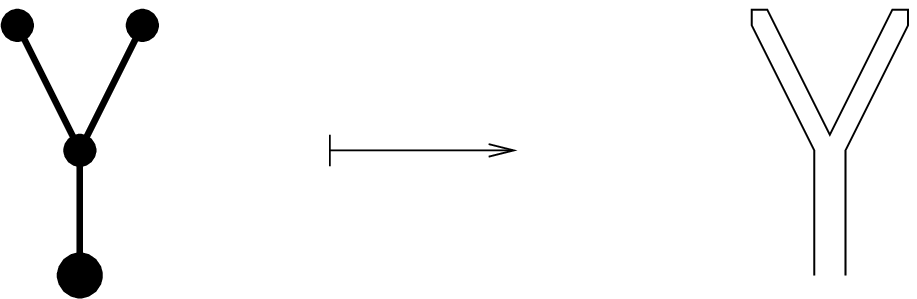}
}
\]


This figure may be called a {\em boundary tree}. Drawing
clockwise produces the same boundary tree up to orientation. 
Traversing the boundary tree anticlockwise, each step is either moving
away from or towards the root. Thus the embedded tree defines a
sequence on any pair of symbols representing these two step types.
(Our example could be the sequence 112122, say.)

\prl(lugosi)
Two such trees are distinct if their sequences are distinct. \Qed
\end{pr}

(Remark: choosing the symbols $($ and $)$ 
gives a bijection between $\CC_t(n)$ and $\CC_b(n)$.)


Next we show that 
the boundary tree realisation leads us to a decomposition analogous
to \ref{s:TLdiagrams}(1) and \ref{ss:Cb}(1). 

Every rooted planar tree has an embedding in which the root lies at
the origin $(0,0)$, some sequence of vertices and edges starting
at the root lie on the positive $y$-axis (and no others), and all vertices
lie in the closure of the positive quadrant.

The edges and vertices lying on the $y$-axis are called the 
{\em trunk} in such an embedding. Each vertex $v$ of the trunk defines a
tree with root $v$ induced by the child vertices of $v$ obtained by
moving along edges not in the trunk; we refer to these as
the \emph{branches} of the tree in such an embedding.
The intersection of the boundary tree of such a tree with the positive
quadrant is called a {\em (right) half-tree}.
Note that, of the outward and return pair of 
boundary tree edges associated to an edge in the trunk,
only the outward edge is retained in the right half-tree
(while both are retained for edges not in the trunk).   
There is a corresponding notion of left half-trees. 

The set of all right half-trees derived from trees with $n$ edges in
total and $l$ edges in the trunk is denoted  $Y_t(n;l)=Y^R_t(n;l)$.
The sets $Y^L_t(n;l)$ of left half-trees are defined similarly. 

\begin{figure}
\includegraphics{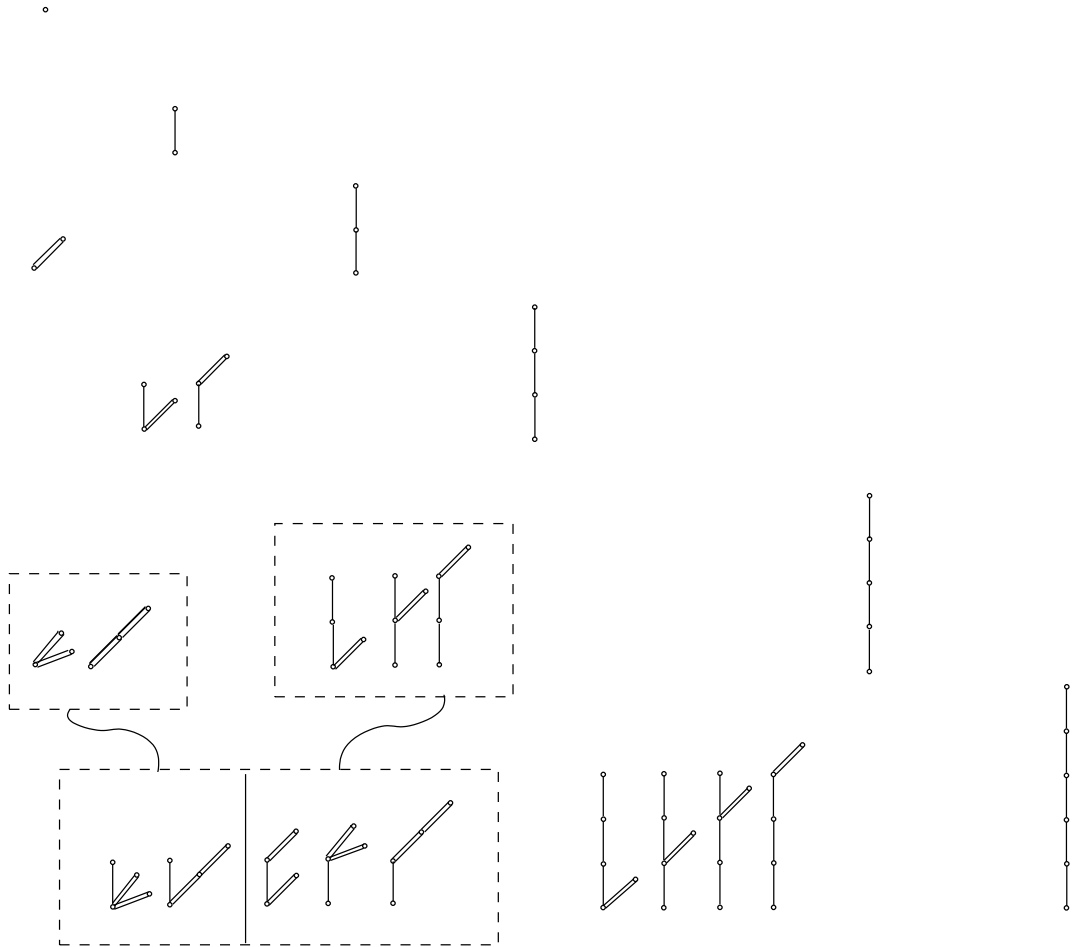}
\caption{\label{halftreeA} The Pascal $(A_{\infty},0)$--array 
$Y_t$, underlying $\CC_t$.
Note that the precise angle of edges is irrelevant.}
\end{figure}

The Pascal $(A_{\infty},0)$--array corresponding to $\CC_t$ may be
constructed as follows (cf. Figure~\ref{halftreeA}).

(1) Bra-ket decomposition: 
Consider an $n$-edge planar tree drawn as described above.
Note that each boundary tree contains $2n$ edges. We cut such a boundary tree
figure into two halves by cutting after the $n^{th}$ edge: 

\[
\raisebox{-.01in}{
\includegraphics[width=7cm]{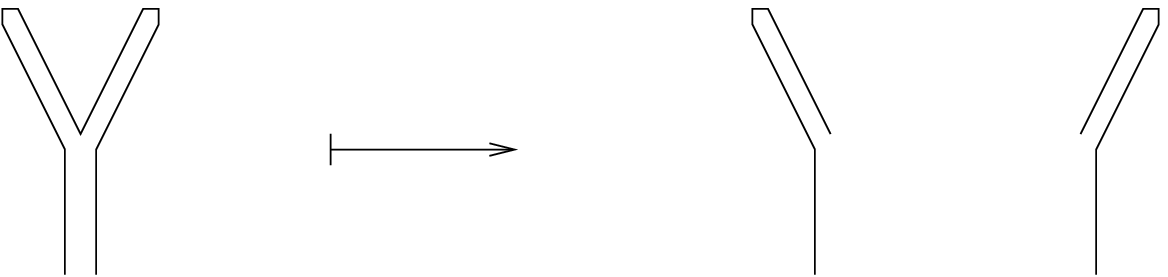}
}
\]

{\lem{ \label{ethoc} This defines a map
\[
\CC_t(n) \rightarrow \bigsqcup_{l=0}^n Y^L_t(n;l) \times Y^R_t(n;l) . 
\]
}}


(2) Edge maps: \\
Let $\phi_t^+$ denote the map from $Y_t(n;l)$ to $Y_t(n+1;l+1)$ which adds
a single edge to the trunk of a half-tree and leaves the branches unchanged.
Let $\phi_t^-$ denote the map from $Y_t(n;l)$ to $Y_t(n+1;l-1)$ which
adds an extra edge to the last edge of the trunk (and rotates it clockwise
off the $y$-axis, together with its branch, without going through another
edge). The pair of edges thus created is regarded as part of the uppermost
branch of the new half-tree. (See Figure~\ref{halftreeA}).

A short argument shows that $\phi_t^+$ and $\phi_t^-$ are both injective
and that
\[
Y_t(n;l) = \phi^+_t(Y_t(n-1;l-1)) \bigsqcup \phi^-_t(Y_t(n-1;l+1))
\]


{\em Proof of Lemma~\ref{ethoc}:}
If we cut the boundary tree of a tree $T$ in $\CC_t(n)$ in half as
described above, the right-hand half lies in $Y_t(n;l)$ for some $l$
(regarding edges in the boundary tree as being incident with the same
vertices as the corresponding vertices in the tree itself).
This is because a step in the path to a higher layer corresponds to
$\phi_t^+$ and a step in the path to a lower layer corresponds to $\phi_t^-$.
Similarly, the left-hand half is the reflection of a tree in $Y_t(n;l')$ in
its trunk, for some $l'$. Since the edges of multiplicity one in each
half come from the same set of edges of $T$, we must have $l=l'$.
\Qed

Given a pair of half-trees in $Y_t(n;l)$, we can splice the first
with the reflection of the second in its trunk to obtain a tree with
all edges of multiplicity two which is the boundary tree of
a tree in $\CC_t(n)$. We denote this map by $s$.
It is clear that this operation is inverse to the bra-ket
extraction.
We have confirmed that:
\prl(p:Y_t)
The array $Y_t(n;l)$ is a Pascal $(A_{\infty},0)$-array
underlying $\CC_t$. \Qed
\end{pr}


It follows from Proposition~\ref{p:equivalence} via a comparison of the edge
maps in the Pascal $(A_{\infty},0)$-arrays $Y_t$ and $Y_b$ that there is a
bijection between $Y_t(n;l)$ and $Y_b(n;l)$
obtained by writing an open-bracket (respectively, close-bracket) 
for each step away from (respectively, towards) the root, as one traverses 
a half-tree anticlockwise from its root passing along each edge exactly once.

A comparison of the procedures for bra-ket extraction for
$\CC_t(n)$ and $\CC_b(n)$ tells us that there is a bijection between
$\CC_t(n)$ and $\CC_b(n)$ given by the same rule applied to the entirety of
the boundary tree of a tree.

A similar argument shows there is a bijection between
$D(n,n)$ and $\CC_t(n)$ obtained by
drawing a vertex in each connected component of the complement of a
diagram in $D(n,n)$ in its bounding rectangle and an edge connecting any two
vertices in components separated by a single line. Again the vertex in the
region touching the western edge of the rectangle is taken to be the root.

\Subsection{Interval orders}
(For a general introduction to interval orders
see \cite{Fishburn85}.)

Consider the set 
$\Chi_n$ 
of $n$--tuples of closed intervals of length 1 on
the real line with distinct end-points. 
It will be convenient for each interval $I$ to be 
considered as a function
$I:[0,1]\rightarrow \mathbb{R}$.
We may label (totally order) the intervals in a tuple in $\Chi_n$ by
starting point: $I^1(0) < I^2(0) < \ldots < I^n(0)$. 
We may
partially order the intervals in a tuple in $\Chi_n$ by stipulating
that 
$I^i<I^j$ if
$I^i(1)<I^j(0)$. The \emph{unit interval orders} of degree $n$ are the
nonisomorphic partial orders 
of $\{ I^1, I^2, \ldots, I^n \}$ 
that can arise in this way, denoted by $\CC_i(n)$. We remark that our
restriction that the end-points are distinct does not change this set
as a small enough displacement of an intervals will not change the
isomorphism class.
 
We have \cite{WineFreund57} that $|\CC_i(n)| = C(n)$.
(See also~\cite[p98]{Fishburn85}).
For example,
\[
\CC_i(3) = \{ (I<J<K), (I<J,K), (I,J<K), (I<K), () \}
\]


Suppose that $0\leq l\leq n$ and $n-2l$ is even.
Consider the set 
$\Chi_n[l]$ 
of tuples of $l$ distinct points and $(n-l)/2$ unit intervals
such that the points and end-points of intervals never coincide
and there is never a point $a$ and an interval $I$ such that
$a<I(0)$.

Partially order the elements in a tuple by ordering the intervals as above.
In addition, if $I$ is an interval and $a$ a point, define $I<a$ if $I(1)<a$.
Thus two points are always incomparable, and a point can never be less
than an interval.
Let $Y_i(n;l)$ denote the set of distinct partial orders arising in this
way. Each such partial order comes with a distinguished subset (given by
the points). We call such labelled partial orders
\emph{unit interval-point} orders.


(1) Bra-ket decomposition: Given a unit interval order of degree $n$,
we consider the lowest $n$ end-points.
If both end-points of an interval lie in this set we keep the
interval, otherwise we just take the endpoint that does lie in the
set. Note that it is
not possible for such a point $a$ to be such that $a<I(0)$ for one of
the intervals $I$ kept in the above procedure,
since then $a$ would be the start of an interval whose
size must be greater than $1$. The corresponding labelled
partial order is thus an element of $Y_i(n;l)$ for some $l$. 
The negations of the remaining end points give rise to
a second element of $Y_i(n;l)$ in the same way. It can be shown that
this procedure is well-defined; the proof involves a description of the
different interval representations of an interval order.

(2) Edge maps.  \\
Let $\phi_i^+$ denote the map from $Y_i(n;l)$ to $Y_i(n+1;l+1)$
given by adding a new labelled element greater in the ordering than all
unlabelled elements and incomparable with all other labelled elements.
Since this corresponds to adding a new point to a representative of the
unit interval-point order greater than all the intervals in it, it is clear
that the new labelled partial order is again a unit interval-point order.

Define the \emph{height} of a labelled element to be the number of
(necessarily unlabelled) elements less than it.
Let $\phi_i^-$ denote the map from $Y_i(n;l)$ to $Y_i(n+1;l-1)$
given by removing the label from any labelled point of minimal height.
It can be checked that the new labelled partial order is indeed an element
of $Y_i(n+1;l-1)$ and is independent of the choice of labelled point.

We remark that two unit interval-point orders $X_1$ and $X_2$ with the same
number of labelled points and orderings $<_1$ and $<_2$ can be combined to
make a new order. The labelled points in $X_1,X_2$ can be ordered in
increasing height. The new order has set $X_1\cup X_2$,
where we identify the labelled points in $X_1$ following this order with the
labelled points in $X_2$ ordered in decreasing height.

We then set $x<y$ in the new order if either:

(a) $x,y\in X_1$ and $x<_1 y$; \\
(b) $x,y\in X_2$ and $y<_2 x$; or \\
(c) $x\in X_1$, $y\in X_2$, and $x,y$ are both unlabelled.

It can be seen that this is a unit interval order using a result of
Roberts~\cite{Roberts69} (see~\cite{BogartWest99} for a short proof)
that states that a finite partial order is an interval order if and only if
it has an interval representation without nested intervals.


The following proposition can be shown via an analysis of interval
representations of interval orders.

\prl(p:intervalorder)
The sets of unit interval-point orders $Y_i(n;l)$ described above form an
$(A_{\infty},0)$-array of sets underlying the sequence of sets $\CC_i(n)$.
\Qed
\end{pr}

\begin{figure}
\[
\includegraphics{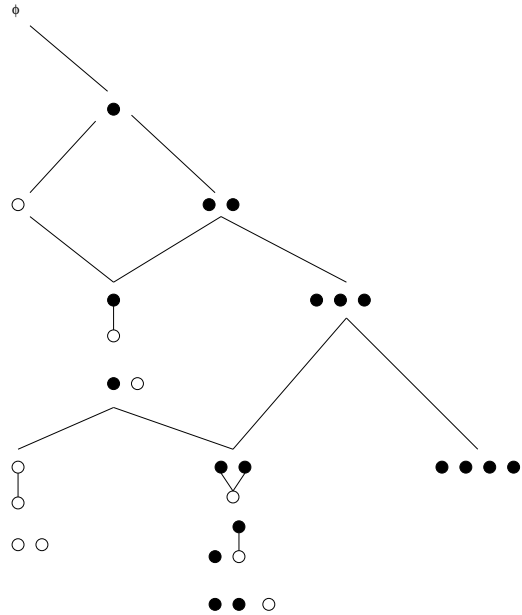}
\]
\caption{\label{f:intervalarray}
The start of the array of unit interval-point orders;
a filled-in circle denotes
a distinguished element, corresponding to a point in a representation.
Empty circles correspond to intervals.}
\end{figure}

The start of the array is displayed in Figure~\ref{f:intervalarray}.
Labelled points are shown as filled-in circles.

It follows from Proposition~\ref{p:equivalence} via a comparison of the edge
maps in the Pascal $A_{\infty}$-arrays $Y_i$ and $Y_b$ that there is a
bijection between $Y_i(n;l)$ and $Y_b(n;l)$
obtained by writing an open-bracket (respectively, close-bracket) 
for each point or lower end-point of an interval (respectively, upper
end-point of an interval) reading along the real line.

\begin{figure}[htbp]
\[
\includegraphics{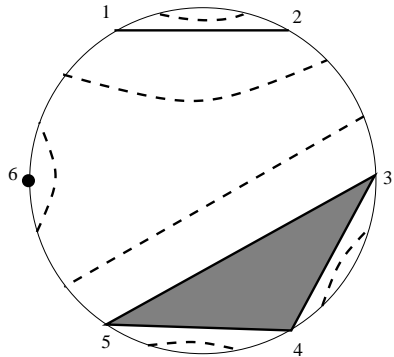}
\]
\caption{Example of correspondence between noncrossing pair partitions and
noncrossing partitions. The dotted lines represent the noncrossing pair partition corresponding to the noncrossing partition $\{\{1,2\},\{3,4,5\},\{6\}\}$.}
\label{f:noncrossing}
\end{figure}

A comparison of the procedures for bra-ket extraction for
$\CC_i(n)$ and $\CC_b(n)$ tells us that there is a bijection between
$\CC_i(n)$ and $\CC_b(n)$ given by the same rule applied to an interval
order (where the point case does not arise).
\Subsection{Noncrossing partitions\label{noncrossing}}
Let $n\in\mathbb{N}$. A \emph{noncrossing partition} of $n$ is a partition
of $n$ such that, if the numbers $1,2,\ldots ,n$ are equally spaced
(in an anticlockwise order) around
a circle in the plane then the convex hulls of the parts do not intersect
\cite{Kreweras72}
(see~\cite{Mccammond06} for a recent survey).
There is a bijection with noncrossing pair partitions of $2n$
obtained by drawing around the outside of the boundaries of the convex hulls
(see \cite[\S6.2.1]{Martin91}, or Figure~\ref{f:noncrossing} for an
example). Our usual set of constructions follows
from Section~\ref{s:TLdiagrams} via this bijection.


\Subsection{Clusters} \label{clusters}
In order to study the (dual) canonical basis of a quantized enveloping
algebra of a simple Lie algebra, Fomin and Zelevinsky~\cite{FominZelevinsky02}
have introduced the notion of a \emph{cluster algebra}. Cluster algebras
of finite type have been classified~\cite{FominZelevinsky03b}, and are
parametrized (in families) by the Dynkin diagrams. The generators of a
cluster algebra (apart from the coefficients) are organised into distinguished
subsets all of the same
size, known as \emph{clusters}. In the finite type case, the generators can be
parametrized by the \emph{almost positive} roots $\Phi_{\geq -1}$
(i.e. the positive roots together with the negative simple roots) in the
root system $\Phi$ of the corresponding Dynkin diagram, and the clusters
can be described using root system combinatorics; the subsets of
$\Phi_{\geq -1}$ corresponding to clusters are known as
\emph{root clusters}.

Root clusters can be regarded as a special case of the $\Gamma$-clusters
of~\cite{MRZ03} associated to each orientation $\Gamma$ of the
Dynkin diagram; here $\Gamma$ is taken to be the alternating orientation,
in which each vertex is a sink or source.

Since $\Gamma$-clusters are known to form a sequence of sets with
Catalan cardinalities in
type $A_n$ (this can be seen from~\cite[3.8]{FominZelevinsky03a}
and~\cite[4.3,4.5]{BMRRT06}) it is natural to ask whether
they fit into the framework of this paper (as we have already mentioned
in the introduction). In this section, we show that in the simplest
case, i.e. root clusters associated to the linear orientation of the
Dynkin diagram of type $A_n$, this can be done.
We describe a Pascal array for clusters that was constructed using
the bijections given in Reading~\cite{Reading05} between Coxeter-sortable
elements and clusters and between Coxeter-sortable elements and noncrossing
partitions. See Section~\ref{noncrossing} for information concerning
the Pascal array underlying the sequence of sets of noncrossing partitions.



Reading~\cite{Reading05} gives an elementary description of
$\Gamma$-clusters which we use here. (Note that Reading generalises
$\Gamma$-clusters to include the non-crystallographic case also).
In Reading's approach, the clusters depend on a choice of Coxeter element
$c$ in the Weyl group $W$ of $\Phi$, and are therefore known as $c$-clusters.
This is a natural generalisation of the initial combinatorics of clusters
given in~\cite{FominZelevinsky03a}; see also~\cite{FominZelevinsky03b}.
We regard $c$-clusters as subsets of $\Phi_{\geq -1}$ as in Fomin-Zelevinsky's
work.

Let $s_1,s_2,\ldots ,s_n$ be the simple reflections in $W$ corresponding
to the simple roots $\alpha_1,\alpha_2,\ldots ,\alpha_n$ in the root
system $\Phi$ of type $A_n$. Note that
$$\Phi=\{\pm(\alpha_i+\alpha_{i+1}+\cdots+\alpha_j)\,:\,1\leq i\leq j\leq n\}.$$
For positive roots $\alpha=\alpha_i+\alpha_{i+1}+\cdots +\alpha_j$
and $\alpha'=\alpha_{i'}+\alpha_{i'+1}+\cdots +\alpha_{j'}$ we write
$\alpha\subseteq \alpha'$ if $i\leq i'\leq j'\leq j$.

For each $i=1,2,\ldots ,n$ let $\sigma_i$ be the involution of
$\Phi_{\geq -1}$ given by the formula:
$$\sigma_i(\alpha):=\left\{\begin{array}{cc}
-\alpha, & \text{if $\alpha=\pm \alpha_i$}, \\
\alpha, & \text{if $\alpha=-\alpha_j$, $j\not=i$}, \\
s_i(\alpha), & \text{otherwise}.
\end{array}\right.$$


Define a family of binary relations $||_c$, for $c$ a Coxeter element,
on the set $\Phi_{\geq -1}$, unique satisfying:

(I) For any $-\alpha_i$ and any $\alpha\in \Phi^+$,
$-\alpha_i||_c \alpha$ if and only if $\alpha_i$ does not occur in
the expansion of $\alpha$ in terms of the simple roots.

(II) For any $\alpha,\beta\in\Phi_{\geq -1}$ and any initial letter
$s$ (equiv. final letter) of $c$, $\alpha||_c\beta$ if and only if
$\sigma(s)(\alpha)||_{scs} \sigma(s)(\beta)$.

Then a $c$-cluster is a maximal subset of $\Phi_{\geq -1}$ of elements
all pairwise $c$-compatible in this sense. Note that it follows from repeated
application of (II) that $\alpha||_c\beta$ if and only if
$c(\alpha)||_c c(\beta)$ for all $\alpha,\beta\in\Phi_{\geq -1}$.
In type $A_n$, we fix the Coxeter element $c=s_ns_{n-1}\cdots s_1$, and
in the sequel we will drop the $c$ from $c$-cluster.

We set $\CC_c(n)$ to be the set of clusters of type $A_{n-1}$.

Fix $n\in\mathbb{N}$ and $l\in\mathbb{N}$ such that $n-l$ is even.
By a \emph{tagged} cluster we mean a cluster together with a
tag associated to each negative root (known as \emph{root tags})
together with an additional (\emph{global}) tag,
which can be $+$ or empty. If a negative root is tagged, we
write a bar over it; otherwise we leave it unadorned.

If $n$ and $l$ are odd, let $Y_c(n;l)$ denote the set of tagged
clusters of type $A_{(n-1)/2}$ with $(l-1)/2$ root tags and no global tag.
If $n$ and $l$ are even, let $Y_c(n;l)$ denote the set of tagged
clusters of type $A_{n/2-1}$ with $l/2$ tags (one of which can be a global
tag).


\emph{1. Bra-ket extraction}

Suppose first that $n$ is odd and let $X$ be a cluster of type $A_{n-1}$.
Let $r=(n-1)/2$.
For each $k\in \{r+1,\ldots ,n-1\}$ consider the sets
$$X_k:=\{\alpha\in X\,:\,\alpha>0,\ \alpha=\alpha_i+\cdots +\alpha_k,\ 1\leq i\leq r\},$$
$$Y_k:=\{\alpha\in X\,:\,\alpha>0,\ \alpha=\alpha_i+\cdots +\alpha_k\},$$
$$I_k:=\{i\in \{1,2,\ldots ,r\}\,:\,\alpha_i+\cdots + \alpha_k\in X_k\}.$$

Suppose first that $|X_k|=1$.
Then $X_k=\{\alpha_i+\cdots +\alpha_k\}$ for some $i\in [1,r]$.
If $|X_t|>1$ for some $t\in [r+1,k-1]$, choose a maximal such $t$;
otherwise let $t=r$.
Then the root $\alpha_i+\cdots +\alpha_k\in X$ is replaced by
$\alpha_{t+1}+\cdots +\alpha_k$.

Suppose secondly that $|X_k|>1$.
Then the roots in $Y_k\subseteq X$ are replaced by the roots
$-\alpha_k$ and the roots $-\alpha_{i-1}$ where $i\in I_k$ is not the
minimum element of $I_k$. The roots $-\alpha_{j-1}$ where $j$ is
the second smallest element of $I_k$ and $-\alpha_k$ are both tagged.
Note that as $|X_k|>1$ we must have $j\leq r$.

Let $X'$ be the new set of roots so obtained and let
\begin{eqnarray*}
C & = & \{\alpha\,:\,\alpha\in X',\alpha\subseteq \alpha_1+\cdots +\alpha_r\}
\\
D & = & \{\sigma(\alpha)\,:\,\alpha\in X',\alpha\subseteq \alpha_{r+1}+\cdots
+\alpha_{n-1}\},
\end{eqnarray*}
where $\sigma$ is the map taking $\alpha_i$ to $\alpha_{n-i}$, extended
linearly.

Then it can be shown that $C$ and $D$ are both clusters of type $A_r$ with the
same number of tags, i.e. elements of $Y_c(r;l)$ for some $l$.

Now suppose that $n$ is even and let $X$ be a cluster of type $A_{n-1}$.
Let $r=n/2-1$.
For each $k\in \{r+2,\ldots ,n-1\}$ consider the sets
$$X_k:=\{\alpha\in X\,:\,\alpha>0,\ \alpha=\alpha_i+\cdots +\alpha_k,\ 1\leq i\leq r+1\},$$
$$Y_k:=\{\alpha\in X\,:\,\alpha>0,\ \alpha=\alpha_i+\cdots +\alpha_k\},$$
$$I_k:=\{i\in \{1,2,\ldots ,r+1\}\,:\,\alpha_i+\cdots + \alpha_k\in X_k\}.$$
We also define:
$$X_{r+1}:=\{\alpha\in X\,:\,\alpha>0,\ \alpha=\alpha_i+\cdots +\alpha_k,\ 1\leq i\leq r\},$$
$$Y_{r+1}:=\{\alpha\in X\,:\,\alpha>0,\ \alpha=\alpha_i+\cdots +\alpha_k\},$$
$$I_{r+1}=\{i\in \{1,2,\ldots ,r\}\,:\,\alpha_i+\cdots + \alpha_k\in X_k\}.$$

If $k>r+1$ and $I_k\subseteq [r+2,n-1]$, or if $k=r+1$ and $I_k=\{r+1\}$,
then the elements of $Y_k$ are left alone.
Otherwise, if $|X_k|=1$ then $X_k=\{\alpha_i+\cdots +\alpha_k\}$ for some
$i\in [1,r+1]$ for $k>r+1$ or $i\in [1,r]$ for $k=r$.
If $k>r+1$ and $|X_t|>1$ for some
$t\in [r+2,k-1]$ then choose a maximal such $t$; otherwise let $t=r+1$,
except that if $k=r+1$ then let $t=r$.

Then the root $\alpha_i+\cdots +\alpha_k\in C$ is replaced by
$\alpha_{t+1}+\cdots +\alpha_k$.

If $|X_k|>1$ then the roots in $Y_k\subseteq X$ are replaced by the root
$-\alpha_k$ and the roots $-\alpha_{i-1}$ where $i\in I_k$ is not the
minimum element of $I_k$. The roots $-\alpha_{j-1}$ where $j$ is
the second smallest element of $I_k$ and $-\alpha_k$
are both tagged.

Let $X'$ be the new set of roots so obtained and let
\begin{eqnarray*}
C & = & \{\alpha\,:\,\alpha\in X',\alpha\subseteq \alpha_1+\cdots +\alpha_r\}
\\
D & = & \{\sigma(\alpha)\,:\,\alpha\in X',\alpha\subseteq \alpha_{r+2}+\cdots
+\alpha_{n-1}\},
\end{eqnarray*}
where $\sigma$ is the map taking $\alpha_i$ to $\alpha_{n-i}$, extended
linearly.

If $X'$ contains $\alpha_r$ then a global tag is
added to each of $C$ and $D$. If $X'$ contains $-\alpha_{r+1}$ then neither
$C$ nor $D$ gets a global tag. If $X'$ contains $-\overline{\alpha_{r+1}}$
then either $C$ or $D$ gets a global tag in such a way as to ensure that they
have the same number of tags.

Then $C$ and $D$ are both clusters of type $A_r$ with the same number of tags,
so they both lie in $Y_c(n;l)$ for some $l$.


\emph{2. Edge maps:}

Given a set $S$ of (possibly tagged) almost positive roots, we need a
function $L_S$ from the set $\{1,2,\ldots ,n\}$ to itself defined as follows.
Let $i\in\{1,2,\ldots ,n\}$. Then $L_S(i)$ is defined to be $j+1$ where
$j<i$ is maximal such that $-\alpha_j\in S$ (possibly tagged). If no such
$j$ exists then $L_S(i):=1$.

We now describe the maps $\phi_c^{\pm}$ on a tagged cluster $C$ in
$Y_c(n;l)$.
The map $\phi_c^+$ is given by adding a global tag if $n$ is odd.
If $n$ is even it is given by adding a new negative simple root
$-\alpha_{n/2}$ which is tagged if there was a global tag already
present (in which case the global tag is removed).

The map $\phi_c^-$ is given by the identity map if $n$ is odd
(and $l\geq 1$).
If $n$ is even it is given in the following way. Note that
$C$ is a tagged cluster of type $A_r$ where $r=n/2-1$, and we are
supposing that $l\geq 1$ so $C$ has at least one tag.

Suppose first that $C$ has no global tag. Let $k$ be maximal such that
$-\overline{\alpha_k}\in C$. (Note that such a $k$ must exist as $C$ must
have at least one tag and has no global tag).
Then every root of the form $-\alpha_l$ with
$l\geq k$ is replaced by $\alpha_{l+1}+\cdots +\alpha_{r+1}$.
The extra root $\alpha_{L_C(k)}+\cdots +\alpha_{r+1}$
is added to the set $C$.

Suppose secondly that $C$ does have a global tag.
Then the extra root $\alpha_{L_C(r)}+\cdots +\alpha_{r+1}$ is added to the
set $C$. In both cases, the global tag is removed.

The start of the array $Y_c$ is shown in Figure~\ref{pascalcluster}, and
an example of the bra-ket extraction is shown in Figure~\ref{clusterbraket}.


\begin{figure}
\begin{center}
$$
{\small
\xymatrix@R=10mm@C=-10mm{
&
\varphi \ar@{-}[dr]
\\
\phi \ar@{-}[ur] \ar@{-}[dr]
&&
\phi,+ \ar@{-}[dr]
\\
&
{\begin{array}{c} -\alpha_1 \\ \alpha_1 \end{array}}
\ar@{-}[ur] \ar@{-}[dr]
&&
-\overline{\alpha_1}
\ar@{-}[dr] 
\\
{\begin{array}{c} -\alpha_1 \\ \alpha_1 \end{array}}
\ar@{-}[ur] \ar@{-}[dr]
&&
{\begin{array}{c} -\alpha_1,+ \\ \alpha_1,+ \\ -\overline{\alpha_1}
\end{array}}
\ar@{-}[ur] \ar@{-}[dr] 
&&
-\overline{\alpha_1},+
\ar@{-}[dr] 
\\
&
{\begin{array}{c} -\alpha_1,-\alpha_2 \\ \alpha_1,-\alpha_2 \\
-\alpha_1,\alpha_2 \\ \alpha_1+\alpha_2, \alpha_1 \\
\alpha_1+\alpha_2,\alpha_2 \end{array}}
\ar@{-}[ur] \ar@{-}[dr] 
&&
{\begin{array}{c} -\overline{\alpha_1},-\alpha_2 \\
-\alpha_1,\overline{-\alpha_2} \\
\alpha_1,-\overline{\alpha_2} \\
-\overline{\alpha_1},\alpha_2 \end{array}}
\ar@{-}[ur] \ar@{-}[dr] 
&&
-\overline{\alpha_1},-\overline{\alpha_2} 
\ar@{-}[dr]
\\
{\begin{array}{c} -\alpha_1,-\alpha_2 \\ \alpha_1,-\alpha_2 \\
-\alpha_1,\alpha_2 \\ \alpha_1, \alpha_1+\alpha_2 \\
\alpha_1+\alpha_2,\alpha_2 \end{array}}
\ar@{-}[ur] \ar@{-}[dr] 
&&
{\begin{array}{cc}
-\alpha_1,-\alpha_2,+ & -\overline{\alpha_1},-\alpha_2 \\
\alpha_1,-\alpha_2,+ & -\alpha_1,-\overline{\alpha_2} \\
-\alpha_1,\alpha_2,+ & \alpha_1,-\overline{\alpha_2} \\
\alpha_1,\alpha_1+\alpha_2,+ & -\overline{\alpha_1},\alpha_2 \\
\alpha_1+\alpha_2,\alpha_2,+ &
\end{array}}
\ar@{-}[ur] \ar@{-}[dr] 
&&
{\begin{array}{c}
-\overline{\alpha_1},-\alpha_2,+ \\
-\alpha_1,\overline{\alpha_2},+ \\
\alpha_1,-\overline{\alpha_2},+ \\
-\overline{\alpha_1},\alpha_2,+ \\
-\overline{\alpha_1},-\overline{\alpha_2}
\end{array}}
\ar@{-}[ur] \ar@{-}[dr] 
&&
-\overline{\alpha_1},-\overline{\alpha_2},+
\ar@{-}[dr] 
\\
&
{\begin{array}{cc}
-\alpha_1,-\alpha_2,-\alpha_3 \\
\alpha_1,-\alpha_2,-\alpha_3 \\
-\alpha_1,\alpha_2,-\alpha_3 \\
\alpha_1,\alpha_1+\alpha_2,-\alpha_3 \\
\alpha_1+\alpha_2,\alpha_2,-\alpha_3 \\
\\
-\alpha_1,-\alpha_2,\alpha_3 \\
\alpha_1,-\alpha_2,\alpha_3 \\
-\alpha_1,\alpha_2,\alpha_2+\alpha_3 \\
\alpha_1,\alpha_1+\alpha_2,\alpha_1+\alpha_2+\alpha_3 \\
\alpha_1+\alpha_2,\alpha_1+\alpha_2+\alpha_3,\alpha_2 \\
\\
\alpha_1+\alpha_2+\alpha_3,\alpha_2+\alpha_3,\alpha_3 \\
-\alpha_1,\alpha_2+\alpha_3,\alpha_3 \\
\alpha_1,\alpha_1+\alpha_2+\alpha_3,\alpha_3 \\
\alpha_1+\alpha_2+\alpha_3,\alpha_2,\alpha_2+\alpha_3
\end{array}}
\ar@{-}[ur]
&&
{\begin{array}{c} 
-\alpha_1,-\alpha_2,-\overline{\alpha_3} \\
\alpha_1,-\alpha_2,-\overline{\alpha_3} \\
-\alpha_1,\alpha_2,-\overline{\alpha_3} \\
\alpha_1,\alpha_1+\alpha_2,-\overline{\alpha_3} \\
\alpha_1+\alpha_2,\alpha_2,-\overline{\alpha_3} \\
\\
-\overline{\alpha_1},-\alpha_2,-\alpha_3 \\
-\alpha_1,-\overline{\alpha_2},-\alpha_3 \\
\alpha_1,-\overline{\alpha_2},-\alpha_3 \\
-\overline{\alpha_1},\alpha_2,-\alpha_3 \\
\\
-\overline{\alpha_1},-\alpha_2,\alpha_3 \\
-\alpha_1,-\overline{\alpha_2},\alpha_3 \\
\alpha_1,-\overline{\alpha_2},\alpha_3 \\
-\overline{\alpha_1},\alpha_2,\alpha_2+\alpha_3 \\
-\overline{\alpha_1},\alpha_2+\alpha_3,\alpha_3
\end{array}}
\ar@{-}[ur]
&&
{\begin{array}{c} 
-\overline{\alpha_1},-\alpha_2,-\overline{\alpha_3} \\
-\alpha_1,-\overline{\alpha_2},-\overline{\alpha_3} \\
\alpha_1,-\overline{\alpha_2},-\overline{\alpha_3} \\
-\overline{\alpha_1},\alpha_2,-\overline{\alpha_3} \\
-\overline{\alpha_1},-\overline{\alpha_2},-\alpha_3 \\
-\overline{\alpha_1},-\overline{\alpha_2},\alpha_3
\end{array}}
\ar@{-}[ur]
&&
-\overline{\alpha_1},-\overline{\alpha_2},-\overline{\alpha_3}
}}
$$
\caption{The start of the array $Y_c$ of tagged clusters}
\label{pascalcluster}
\end{center}
\end{figure} 


\begin{figure}
\begin{center}
$$\begin{array}{c|c|c}
\text{$A_3$-cluster} & C & D \\
\hline
-\alpha_1,-\alpha_2,-\alpha_3 &-\alpha_1 & -\alpha_1 \\
-\alpha_1,-\alpha_2,\alpha_3 &-\alpha_1 & \alpha_1  \\
\alpha_1,-\alpha_2,-\alpha_3 & \alpha_1 & -\alpha_1 \\
\alpha_1,-\alpha_2,\alpha_3 & \alpha_1 & \alpha_1 \\
\hline
-\alpha_1,\alpha_2,-\alpha_3 & \alpha_1,+ & -\alpha_1,+ \\
-\alpha_1,\alpha_2,\alpha_2+\alpha_3 & -\alpha_1,+ & \alpha_1,+ \\
-\alpha_1,\alpha_2+\alpha_3,\alpha_3 & -\alpha_1,+ & -\overline{\alpha_1} \\
\alpha_1,\alpha_1+\alpha_2,-\alpha_3 & \alpha_1,+ & -\alpha_1,+ \\
\alpha_1,\alpha_1+\alpha_2,\alpha_1+\alpha_2+\alpha_3 & \alpha_1,+ & \alpha_1,+ \\
\alpha_1,\alpha_1+\alpha_2+\alpha_3,\alpha_3 & \alpha_1,+ & -\overline{\alpha_1} \\
\alpha_1+\alpha_2,\alpha_2,-\alpha_3 & -\overline{\alpha_1} & -\alpha_1,+ \\
\alpha_1+\alpha_2,\alpha_2,\alpha_1+\alpha_2+\alpha_3 & -\overline{\alpha_1} & \alpha_1,+ \\
\alpha_1+\alpha_2+\alpha_3,\alpha_2+\alpha_3,\alpha_3 & -\overline{\alpha_1} & -\overline{\alpha_1} \\
\hline
\alpha_1+\alpha_2+\alpha_3,\alpha_2,\alpha_2+\alpha_3 &
-\overline{\alpha_1},+ & -\overline{\alpha_1},+ 
\end{array}$$
\end{center}
\caption{Bra-ket extraction for clusters of type $A_3$}
\label{clusterbraket}
\end{figure}


The following Theorem holds. The proof is nontrivial, but will be
given in a separate publication~\cite{MarshMartin07} as it does not fit
our purposes here.

\begin{theorem}
The array $Y_c(n;l)$ of tagged clusters is a Pascal $(A_{\infty},0)$-array
underlying the sequence $\CC_c(n)$ of root clusters.
\end{theorem}

\Subsection{Remarks}
Note that we have {\em not} given a procedure for deriving a 
Pascal array from a Catalan sequence.
Rather we have developed one in each case by some strategy.

Our main strategy is to use ideas from representation theory.
Next we look explicitly at the representation theory side of things.

\section{Pascal arrays arising from algebras} \label{s:Pascalalgebras}
\Subsection{Towers of algebras
associated to rooted graphs
\label{new axioms}}
Our objective in this section is to describe a set of axioms on a tower of
algebras which is enough to associate a Pascal array to the tower
(in the way suggested in section~\ref{s:TLdiagrams}). 
The aim, of course, is to find an axiom set weak enough that concrete
examples exist.


Let $k$ be a field and $A$ be a  $k$-algebra. 
Let $\Delta$ be a set of $A$-modules. 
A filtration of an $A$-module $M$ is called a {\em $\Delta$-filtration}
if it has sections in  $\Delta$.
An $A$-module $M$ is said to be {\em $\Delta$-good} provided it has a
$\Delta$-filtration and the multiplicity of $L \in \Delta$ in any
$\Delta$-filtration of $M$ is independent of the choice of filtration.
We denote this multiplicity by $[M:L]$.

(This definition is motivated by the definition of $\Delta$-good modules over
a quasihereditary algebra).


Let $k$ be a field. Suppose we are given a tower $A_{\bullet}$ of $k$-algebras
$A_0,A_1,A_2,\ldots $. We consider the following axioms:

(N1) For all $n\in\mathbb{N}$, there is an algebra embedding
$A_n\hookrightarrow A_{n+1}$.

(N2) For all $n\in\mathbb{N}$, there is a set $\Lambda_n$, with an
$A_n$-module $\Delta_n(\lambda)$ for each $\lambda\in\Lambda_n$,
such that if $\lambda\not=\mu$ then $\Delta_n(\lambda)\not\cong
\Delta_n(\mu)$.

(N3) There is a fixed $N\in\mathbb{N}$ such that for all $n\in\mathbb{N}$
with $n\geq N$, $\Lambda_{n-N}\subseteq\Lambda_n$. Furthermore,
whenever $\Lambda_n\cap\Lambda_m\not=\phi$, $n\equiv m\mod N$.

By (N1), we have a restriction functor
${}_{A_{n-1}}res_{A_n}:A_n\modules\rightarrow A_{n-1}\modules$ for all 
$n\geq 1$.
For given $n$, set 
$\Delta  = \{ \Delta_n(\lambda) \; | \; \lambda \in \Lambda_n \}$. 

(N4) For all $n\geq 1$ and all $\lambda\in\Lambda_n$,
${}_{A_{n-1}}res_{A_n}(\Delta_n(\lambda))$ is a $\Delta$-good
$A_{n-1}$-module.

(N5) The cardinality of $\Lambda_0$ is $1$. We denote its unique element by
$0$. The module $\Delta_0(0)$ is a one-dimensional simple
$A_0$-module.

(N6) For all $n\in\mathbb{N}$, $n\geq N$, there is an exact functor
$F_n:A_n\modules\rightarrow A_{n-N}\modules$.

(N7) For all $\lambda\in \Lambda_n$, $F_n(\Delta_n(\lambda))$ is
isomorphic to $\Delta_{n-N}(\lambda)$ if $\lambda\in\Lambda_{n-N}$, and
otherwise is zero.

(N8) For all $n\in\mathbb{N}$, $n\geq N+1$, and all $\lambda\in\Lambda_n$,
we have that
$$(F_{n-1}\circ {}_{A_{n-1}}res_{A_n})(\Delta_n(\lambda))\cong
({}_{A_{n-N-1}}res_{A_{n-N}}\circ F_n)(\Delta_n(\lambda)),$$
provided $F_n(\Delta_n(\lambda))\not=0$.


The $\Delta$-{\em Bratteli diagram} 
of $A_{\bullet}$  is defined to be the array
$(\Delta_n(\lambda))_{n,\lambda}$ of modules.
We draw $m$ arrows from $\Delta_{n}(\lambda)$ to $\Delta_{n+1}(\mu)$
whenever $\Delta_{n}(\lambda)$ appears as a factor in
${}_{A_{n}}res_{A_{n+1}}(\Delta_{n+1}(\mu))$ with multiplicity $m>0$.

(N9) There is a graph $G$ with vertices (indexed by)
$\Lambda:=\cup_{n\in\mathbb{N}}\Lambda_n$ with the following property.
For every arrow $\alpha$ from $\lambda$ to $\mu$ in $G$ there
is a corresponding arrow $\alpha_n$ in the Bratelli diagram from
$\Delta_n(\lambda)$ to $\Delta_{n+1}(\mu)$ for any $n$ such that
$\lambda\in\Lambda_n$. These arrows are all distinct and,
as $\lambda,\mu$ vary over $\Lambda$, they exhaust the arrows in the
Bratteli diagram.

Note that $(G,0)$ is a rooted graph.

If $A_{\bullet}$ satisfies axioms (N1)-(N9) then we say that $A_{\bullet}$ is a
\emph{$(G,0)$-tower of algebras}.

We shall see later that, for any tower of algebras satisfying (N1)-(N8),
such a graph $G$ exists, so that (N9) is automatic for the right choice
of rooted graph $(G,0)$.

There is a notable class of examples of (N1-4) with $N=2$ given by
Jones' Basic Construction \cite[\S2.4]{GoodmanDelaharpeJones89}.
(Indeed $F_n$ can be related to the conditional expectation
\cite{Martin2000},
but we will not elaborate here.)


\begin{lemma} \label{Flift}
Suppose we have a tower $A_{\bullet}$ of algebras satisfying axioms (N1)-(N8)
as above. Let $\lambda\in\Lambda_{n-N}$, $\mu\in\Lambda_{n-N+1}$.
Suppose that there are $m>0$ arrows in the Bratelli diagram from
$\Delta_{n-N}(\lambda)$ to $\Delta_{n-N+1}(\mu)$.
Then there are $m>0$ arrows in the Bratteli diagram from
$\Delta_{n}(\lambda)$ to $\Delta_{n+1}(\mu)$. 
\end{lemma}

{\bf Proof:} By assumption, $F_{n}(\Delta_{n}(\lambda))$ appears exactly
$m>0$ times as a quotient in any $\Delta$-filtration of
${}_{A_{n-N}}res_{A_{n-N+1}}F_{n+1}(\Delta_{n+1}(\mu))$. Suppose that
${}_{A_{n}}res_{A_{n+1}}(\Delta_{n+1}(\mu))$
has a filtration with quotients $\Delta_{n}(\nu)$ for $\nu\in X$ for
some index set $X$, each occurring with multiplicity $m_{\nu}>0$.
Then, since $F_{n}$ is exact,
$F_{n}{}_{A_{n}}res_{A_{n+1}}(\Delta_{n+1}(\mu))$ has a filtration with
quotients
$F_{n}\Delta_{n}(\nu)$ for $\nu\in X'$, each with multiplicity $m_{\nu}$,
where $X'\subseteq X$ is the subset of those $\nu\in X$ such that
$F_{n}\Delta_{n}(\nu)\not=0$.
Hence ${}_{A_{n-N}}res_{A_{n-N+1}}F_{n+1}(\Delta_{n+1}(\mu))$ has such a filtration (by axiom (N8)).
It follows that $F_{n}\Delta_{n}(\lambda)\cong F_{n}\Delta_{n}(\nu)$
for some $\nu\in X'\subseteq X$, and $m=m_{\mu}$.
Hence $\Delta_{n}(\lambda)\cong \Delta_{n}(\nu)$
by axiom (N7). Hence $\mu=\nu$, by axiom (N2). So $\Delta_{n}(\lambda)$
appears in a filtration of ${}_{A_{n}}res_{A_{n+1}}\Delta_{n+1}(\mu)$, with
multiplicity $m_{\mu}$, as required.~$\Box$


The next lemma shows that there is a lot of uniformity in the Bratteli
diagram of a tower of algebras satisfying (N1)-(N8).

\begin{lemma} \label{Bratellilift}
Suppose that we have a tower $A_{\bullet}$ of algebras satisfying axioms
(N1)-(N8).
Suppose that there are $n_1,n_2\in\mathbb{N}$ such that
$\lambda\in\Lambda_{n_1}\cap \Lambda_{n_2}$.
Suppose further that there are $m$ arrows from $\Delta_{n_1}(\lambda)$ to
$\Delta_{n_1+1}(\mu)$ in the Bratteli diagram.
Then there are $m$ arrows from $\Delta_{n_2}(\lambda)$ to
$\Delta_{n_2+1}(\mu)$ in the Bratteli diagram.
\end{lemma}

{\bf Proof:}
For simplicity, we drop the subscripts from the functors $F_n$ and
${}_{A_{n-1}}res_{A_n}$. Since $\lambda\in\Lambda_{n_1}\cap \Lambda_{n_2}$,
we have that $n_1$ is congruent to $n_2$ modulo $N$ by axiom (N3).

Case (I): We first suppose that $n_1=n_2+rN$ for some integer $r>0$,
so that $F^r(\Delta_{n_1}(\lambda))\cong \Delta_{n_2}(\lambda)$.
We note that $F^r(\Delta_{n_1+1}(\mu))$ is either zero or
isomorphic to $\Delta_{n_2+1}(\mu)$, by axiom (N7). But
$\Delta_{n_1}(\lambda)$ appears in a $\Delta$-filtration of
$res\Delta_{n_1+1}(\lambda)$
and $F^r(\Delta_{n_1}(\lambda))=\Delta_{n_2}(\lambda)\not=0$, so, since $F$ is
exact, $F^r(\Delta_{n_1+1})(\mu)\cong \Delta_{n_2+1}(\mu)$ is not zero.

Hence $M:=F^r res(\Delta_{n_1+1}(\mu)) \cong res F^r(\Delta_{n_1+1}(\mu))=
\Delta_{n_2+1}(\mu)$, the isomorphism by axiom (N8).
By the exactness of $F$, this module has
$F^r(\Delta_{n_1}(\lambda)) \cong \Delta_{n_2}(\lambda)$ as a quotient
in a $\Delta$-filtration; it appears with multiplicity $m$,
so there are $m$ arrows from $\Delta_{n_2}(\lambda)$ to $\Delta_{n_2+1}(\mu)$
in the Bratteli diagram as required.

Case (II): We now suppose that $n_1=n_2-rN$ for some $r>0$, so that
$F^r\Delta_{n_2}(\lambda)\cong \Delta_{n_1}(\lambda)$. Since
$\mu\in\Lambda_{n_1+1}$, we know that
$F^r\Delta_{n_2+1}(\mu)\cong \Delta_{n_1+1}(\mu)$. Furthermore,
there are $m$ arrows in the Bratteli diagram from $\Delta_{n_1}(\lambda)$ to
$\Delta_{n_1+1}(\mu)$, i.e. from $F^r(\Delta_{n_2}(\lambda))$ to
$F^r(\Delta_{n_2+1}(\mu))$, and that both of these are non-zero. 
Repeated application of Lemma~\ref{Flift} gives $m$ arrows in the Bratteli
diagram from $\Delta_{n_2}(\lambda)$ to $\Delta_{n_2+1}(\mu)$ as
required.~$\Box$

\del(graphG)
Let $G$ be the graph with vertices
$\Lambda:=\cup_{n\in\mathbb{N}}\Lambda_n$. Suppose that $\lambda\in\Lambda_n$
and $\mu\in\Lambda_{n+1}$. Then let $G$ have $m$ arrows from $\lambda$ to
$\mu$, where $m$ is the number of arrows in the Bratteli diagram from
$\Delta_n(\lambda)$ to $\Delta_{n+1}(\mu)$ (note that this well-defined
by Lemma~\ref{Bratellilift}). We choose a correspondence
$\alpha\rightarrow \alpha_n$ between the arrows from $\lambda$ to $\mu$ in
$G$ and the arrows in the Bratteli diagram from $\Delta_n(\lambda)$ to
$\Delta_{n+1}(\mu)$ in the Bratteli diagram.
\end{de}

We have the following:

\begin{cor} \label{corgraphG}
Suppose that $A_{\bullet}$ is a tower of algebras satisfying axioms (N1)-(N8).
Let $G$ be the graph defined as above. Then $A_{\bullet}$ is a $(G,0)$-tower
of algebras.
\end{cor}

{\bf Proof:} This is immediate from the definition of $G$.
\Qed

{\de{\label{d:YA}
Define the array $Y_A$ of sets by setting $Y_A(n;\lambda)$ to be the set
of paths in the Bratteli diagram from $\Delta_0(0)$ to $\Delta_n(\lambda)$,
for $n\in\mathbb{N}$ and $\lambda$ a vertex of $G$ such that
$\lambda\in\Lambda_n$.
}}


\begin{theorem}
Let $(G,0)$ be a rooted graph, $A_0,A_1,\ldots $ a
$(G,0)$-tower of algebras and $Y_A$ the
corresponding array of sets defined above.
Then the array $Y_A$ is equivalent to $Y_{G,0}$. In other words,
$Y_A$ is a Pascal $G$-array (see Definition \ref{p:array}).
\end{theorem}

{\bf Proof:}
It is clear that $Y_A$ has the same index sets ($\Lambda$ and $\mathbb{N}$)
as $Y=Y_{G,0}$. By definition, for $\lambda=0\in\Lambda_0$,
there are $m>0$ arrows from $0$ to $\mu$ in $G$ if and only if there
are $m>0$ arrows from $\Delta_0(0)$ to $\Delta_1(\mu)$ in the Bratteli
diagram. It follows that there are pointwise bijections between
the sets in layer $n=1$ of $Y_A$ and $Y_{G,0}$.

Suppose that there is an arrow $\alpha$ from $\lambda$ to $\mu$ in $G$.
Define $\edge{\alpha}:A(n;\lambda)\rightarrow A(n+1,\mu)$ composing
a path from $\Delta_0(0)$ to $\Delta_{n}(\lambda)$ in the Bratelli diagram
with the arrow $\alpha_n$ in the Bratteli diagram
from $\Delta_{n}(\lambda)$ to $\Delta_{n+1}(\mu)$. Since the resulting path
together with the arrow $\alpha_n$ determine the original path to
$\Delta_{n-1}(\lambda)$, it is clear that $\edge{\alpha}$ is injective.

Any path $\pi$ in the Bratteli diagram from $\Delta_0(0)$ to $\Delta_n(\mu)$
for some $\mu\in\Lambda_n$ must have visited some module
$\Delta_{n-1}(\lambda)$ immediately before finishing at $\Delta_n(\mu)$.
Let $\pi'$ be the path $\pi$ with the final arrow truncated,
finishing at $\Delta_{n-1}(\lambda)$.
The last arrow in $\pi$, which goes from $\Delta_{n-1}(\lambda)$ to
$\Delta_n(\mu)$,
must be of the form $\alpha_n$ where $\alpha$ is an arrow from $\lambda$
to $\mu$ in $G$, and it is clear that $\edge{\alpha}(\pi')=\pi$.

Finally, suppose that $\edge{\alpha'}(\pi'')=\pi$ for some path $\pi''$ in the
Bratelli diagram from $\Delta_0(0)$ to $\Delta_{n-1}(\nu)$ and some arrow
$\alpha'$ from $\nu$ to $\mu$ in $G$. Then the last arrow in $\pi$ must be
$\alpha'_n$, so $\alpha'_n=\alpha_n$, so $\alpha=\alpha'$. Then $\pi''=\pi'$
since $\edge{\alpha}$ is injective. The result is proved.
\Qed

\begin{cor} \label{Gtowercorrespondance}
The array of sets $Y_A$ is in pointwise explicit bijection with
the array $Y_{G,0}$, i.e. there is a bijection between the set of oriented
paths in the Bratteli diagram from $\Delta_0(0)$ to $\Delta_n(\lambda)$
and the set of paths in $G$ from $0$ to $\lambda$. \Qed
\end{cor}


We have the following:
\begin{prop} \label{towerbases}
Let $(G,0)$ be a rooted graph, $A_{\bullet}$ a $(G,0)$-tower of algebras and
$Y_A$ the corresponding array of sets defined above. Then there is a
one-to-one correspondence between $Y_{G,0}(n;\lambda)$ and a basis for
$\Delta_n(\lambda)$.
\end{prop}

{\bf Proof:}
The elements of $Y_{G,0}(n;\lambda)$ are in one-to-one correspondance with
$Y_{A}(n;\lambda)$ by Corollary~\ref{Gtowercorrespondance}, i.e. with paths
in the Bratteli diagram from $\Delta_0(0)$ to $\Delta_n(\lambda)$, so we have
to show that such paths correspond to a basis for $\Delta_n(\lambda)$.
We do this by induction on $n$; since $\Delta_0(0)$ is one-dimensional
(axiom N5), it holds for $n=0$, so suppose that it holds for $n-1$.
The result for $\Delta_n(\lambda)$ then follows from the fact that it
holds for $n-1$, and thus for each section in a $\Delta$-composition
series for ${}_{A_{n-1}}res_{A_{n}}\Delta_n(\lambda)$, noting that each
such section appears the same number of times as the number of arrows
in the Bratteli diagram from that section to $\Delta_n(\lambda)$.
\Qed


As we shall see, the $G$-tower axiom set is satisfied by many interesting
examples. 
However it is not claimed that it is unique in meeting the stated
objective. In melding the combinatorial array and representation
theoretic aspects of the problem there are various possibilities as to
how one might impose a suitable structure on the representation theory
in a single layer. One could choose to require semisimplicity
(which would bring us close to the `classical' towers of 
\cite{GoodmanDelaharpeJones89}); or cellularity \cite{GrahamLehrer96};
or tabularity \cite{GreenRM01}; or quasiheredity \cite{DlabRingel92}.
Here we have a setup which works (at least) with quasiheredity,
as we shall see in Section~\ref{towersofrecollement}.
This is simply because the tower structure we need is already partly
integrated with quasiheredity, in \cite{\cmpx}. 



\Subsection{Undirected Case}
Recall~\eqref{e:sumofsquares} that if $G$ is undirected with vertex $v_0$ then
$$
N_{G,v_0}(2n;v_0)=\sum_{x\in V} N_{G,v_0}(n;x)^2.
$$
A path of length $2n$ from $v_0$ to $v_0$ can be regarded as two paths
of length $n$ from $v_0$ to some vertex $l$.
Let $(C(n))_{n\in \mathbb{N}}$ be a Catalan $(G,v_0)$-sequence. Then
by definition there is a bijection between $C(n)$ and the set of such pairs of
paths, and therefore with $P_G(2n;v_0,v_0)$, and thus with the set
$Y(2n;v_0)$ for any Pascal $(G,v_0)$-array $Y$: thus the Catalan
$(G,v_0)$-sequence appears in the Pascal array itself.

Suppose that $A_{\bullet}$, together with modules $\Delta_n(\mu)$, is
a $(G,v_0)$-tower of algebras. Then for each $\mu\in \Lambda_n$,
in the restriction of $\Delta_{2n}(v_0)$ to $A_n$, the multiplicity of
any module $\Delta_n(\mu)$ in a $\Delta$-filtration is equal to the
number of paths from $v_0$ to $\mu$, i.e. $\dim(\Delta_n(\mu))$. 
That is,
the multiplicity of any module $\Delta_n(\mu)$ is equal to its own dimension.

This situation is particularly nice 
if each $A_n$ is semisimple and the
modules $\Delta_n(\mu)$ are the simple $A_n$-modules. We see that the
module $\Delta_{2n}(v_0)$ is isomorphic to the regular module ${}_{A_n}A_n$.
Better still is the case in which the modules can be defined freely over
some small ground ring, with extensions over fields which include both
semisimple and non-semisimple cases (cf. \cite{JamesKerber81}). 
The semisimple cases then ensure the above combinatoric
for the non-semisimple cases.
This is exemplified in the case of the tower $TL_{\bullet}$
(see Sections~\ref{s:TLdiagrams},\ref{TLA}).

Let us consider this in a little more detail.
Let $R$ be a ring. 
Suppose that we have a sequence of finite rank 
$R$-algebras $\{ A^R_n \}_n$, 
such that each $ A^R_n$ is a free $R$-module, with basis $B_n$. 
For each homomorphism $\phi$ from $R$ to 
a field $k$ which makes $k$ an $R$-algebra
we get a $k$-algebra by $A^{\phi}_n = k \otimes_R A^R_n$. 
Suppose further that 
there is a filtration of the left regular module $RB_n$ of 
$A^R_n$ by submodules, each of which is a free $R$-module, and that
there is a set of left modules
$\Delta^R(\lambda)$ with basis, 
such that every section in the filtration is isomorphic to some
$\Delta^R(\lambda)$; 
and that these pass to a complete set of simple modules in
some semisimple algebra $A^{\phi}_n$. 
(This might underly a cellularity or a quasiheredity 
\cite{GrahamLehrer96,DlabRingel92} say.)
Then the rank of  $A^R_n$ is the sum of squares of ranks of the
$\Delta^R(\lambda)$s. 
Passing to any  $A^{\phi}_n$ that is
not semisimple we have a sum of squares formula like (\ref{sumLL}) for the 
$\Delta$s, and a distinct formula (\ref{sumPL}) for the simples and
projectives
giving the same sum. 
We shall give examples in Section~\ref{s:truncation}.





\Subsection{Towers of recollement\label{towersofrecollement}}
We next explain the link between the systems of algebras described above
and the towers of recollement considered in~\cite{\cmpx}.
We firstly recall the set-up from \cite{\cmpx}. Let $k$ be an algebraically
closed field, $A$ a finite dimensional algebra over $k$ and $e$ an idempotent
in $A$. 
We then have the {\em localisation} functor $F:A\modules\rightarrow eAe\modules$
taking an $A$-module $M$ to $eM$, 
and the {\em globalisation} functor $G:eAe\modules
\rightarrow A\modules$ taking an $eAe$-module $N$ to $Ae\otimes_{eAe} N$. We have:

\begin{theorem} (GREEN~\cite{Green80}) \label{greenb} \\
Let $\{L(\lambda), \lambda\in \Lambda\}$ be a full set of simple $A$-modules,
and set $\Lambda^e=\{\lambda\in \Lambda\,:\,eL(\lambda)\not=0\}$. Then
$\{eL(\lambda)\,:\,\lambda\in\Lambda^e\}$ 
is a full set of simple $eAe$-modules.
The remaining simple modules $L(\lambda)$ (with $\lambda\in\Lambda\setminus
\Lambda^e$) are a full set of simple $A/AeA$-modules.~$\Box$
\end{theorem}

Suppose we have a family of finite dimensional algebras over $k$,
$\mathcal{A}=A_0,A_1,A_2,\ldots$, with idempotents $e_n\in A_n$. Let
$\Lambda_n$ denote the full set of simple $A_n$ modules, for each $n$,
and let $\Lambda^n$ denote the full set of simple $A_n/A_ne_nA_n$-modules.
In \cite{\cmpx} the following axioms are considered.

(A1) For each $n\geq 2$, there is an isomorphism $\Phi_n:A_{n-2}\rightarrow
e_nA_ne_n$.

We remark that by Theorem~\ref{greenb}, we have
$$\Lambda_n=\Lambda^n\sqcup \Lambda_{n-2}.$$
For $m,n\in\mathbb{N}$ with $m-n$ even, we set $\Lambda^n_m=\Lambda^n$ regarded
as a subset of $\Lambda_m$ if $m\geq n$ and $\Lambda^n_m=\phi$ otherwise.
Setting $e_{n,0}=1\in A_n$, for $1\leq i\leq n/2$ we define new idempotents
in $A_n$ by setting $e_{n,i}=\Phi_n(e_{n-2,i-1})$. We also set $A_{n,i}=
A_n/(A_ne_{n,i+1}A_n)$.

(A2) (i) The algebra $A_n/A_ne_nA_n$ is semisimple. \\
(ii) For each $n\geq 0$ and $0\leq i\leq n/2$, the surjective multiplication
map $A_ne_{n,i}\otimes_{e_nA_ne_n} e_{n,i}A_n\rightarrow A_ne_{n,i}A_n$ is
bijective.

It is remarked in~\cite{\cmpx} 
that axioms (A1) and (A2) are enough to obtain that
each $A_i$ is a quasihereditary algebra. Let $\Delta_n(\lambda)$,
$\lambda\in \Lambda_n$, denote the standard modules for $A_n$.

(A3) For each $n\geq 0$, the algebra $A_n$ can be identified with a subalgebra
of $A_{n+1}$.

This enables us to define the restriction functor
$\mbox{res}_n:A_n\modules\rightarrow A_{n-1}\modules$
and the induction functor $\mbox{ind}_n:A_n\modules\rightarrow A_{n+1}\modules$
given by $ind_n(M)=A_{n+1}\otimes_{A_n} M$.

(A4) For all $n\geq 1$, we have that $A_ne_n\cong A_{n-1}$ as a left $A_{n-1}$-,
right $A_{n-2}$- bimodule.

Here we will also consider the following opposite version of (A4):

(A4') For all $N\geq 1$, we have that $e_nA_n\cong A_{n-1}$ as a right
$A_{n-1}$, left $A_{n-2}$-bimodule.

We divide the next axiom into two parts, in order to better relate it
to (N1)--(N8):

(A5) (a) For each $\lambda\in\Lambda_n^m$, we have that
$res(\Delta_n(\lambda))$ has a $\Delta$-filtration. \\
(b) Furthermore, 
$$supp(res(\Delta_n(\lambda)))\subseteq \Lambda_{n-1}^{m-1}\sqcup
\Lambda_{n-1}^{m+1}.$$
This implies that
$$supp(ind(\Delta_n(\lambda)))\subseteq \Lambda_{n+1}^{m-1}\sqcup
\Lambda_{n+1}^{m+1}.$$

(A6) For each $\lambda\in \Lambda_n^m$ there exists $\mu\in\Lambda_{n+1}^{m-1}$
such that $\lambda\in supp(res(\Delta_{n+1}(\mu)))$.


A tower of algebras of satisfying (A1),(A2),(A3),(A4),(A5) and (A6) is
called a \emph{tower of recollement} in \cite{\cmpx}.
We can now relate $(G,0)$-towers of algebras with towers of recollement.

\begin{prop}\label{ToR-GT}
Suppose that $k$ is an algebraically closed field, and $A_0,A_1,A_2,\ldots$
is a sequence of $k$-algebras satisfying axioms (A1),(A2),(A3),(A4),(A4'),
(A5)(a) and (N5). Then $A_0,A_1,A_2,\ldots $ is a $(G,0)$-tower of algebras
for some graph $G$.
\end{prop}

{\bf Proof:}
By (A3), (N1) holds.
By (A1) and (A2) we know that each $A_i$ is quasihereditary, and we can
take the standard modules $\Delta_n(\lambda)$ for $\lambda\in\Lambda_n$
for (N2). By the comment after (A1), we see that $\Lambda_{n-2}\subseteq
\Lambda_n$ for all $n\geq 2$. It is clear from (A1) and Green's theorem,
Theorem \ref{greenb}, that, by appropriate choice of the $\Lambda^n$, we
can ensure that the $\Lambda_n$ are chosen so that
$\Lambda_n\cap \Lambda_m\not=\phi$ implies that $n-m$ is even, so (N3) holds.
The axiom (N4) follows from axiom (A5)(a). For (N6), we take $F_n$ to be
the localisation functor taking an $A_n$-module $M$ to $e_nM$; note
that this an exact functor. The axiom (N7) follows from axioms (A1) and (A2)
(as stated in~\cite[p3]{\cmpx}, this follows, for example,
from \cite[A1]{donkin98}). For (N8), we have:
\begin{eqnarray*}
F_{n+1}ind_n M &=& e_{n+1}(A_{n+1}\otimes_{A_n} M) \\
&=& e_{n+1}A_{n+1}\otimes_{A_n} M \\
&\cong & A_n \otimes_{A_n} M ={}_{A_{n-1}}res_{A_n} M
\end{eqnarray*}
as $A_{n-1}$-modules.
If $F_n\Delta_n(\lambda)\not=0$, then
\begin{eqnarray*}
F_{n-1}{}_{A_{n-1}}res_{A_n}\Delta_n(\lambda) &\cong
    & F_{n-1}{}_{A_{n-1}}res_{A_n} G_{n-2}\Delta_{n-2}(\lambda) \\
&\cong & F_{n-1} ind_{n-2}\Delta_{n-2}(\lambda) \\
&\cong & {}_{A_{n-3}}res_{A_{n-2}}\Delta_{n-2}(\lambda) \\
&\cong & {}_{A_{n-3}}res_{A_{n-2}}F_n\Delta_n(\lambda),
\end{eqnarray*}
using the fact that $G_{n-2}\Delta_{n-2}(\lambda)\cong \Delta_n(\lambda)$
which relies on (A1) and (A2), and the fact that
$${}_{A_{n-1}}res_{A_n} G_{n-2}(\Delta_{n-2}(\lambda))\cong ind_{n-2}(\Delta_{n-2}(\lambda)),$$
which relies on (A1), (A2), (A3) and (A4) (see~\cite{\cmpx}).
If we define the rooted graph $(G,0)$ as in Definition~\ref{graphG}, we are
done (see Corollary~\ref{corgraphG}). $\Box$

\Subsection{The Temperley-Lieb algebra\label{TLA}}
\newcommand{\llambda}{l}
For each $\delta\in k$,
the Temperley-Lieb algebras, $TL_n=TL_n(\delta)$,
$n=0,1,2,\ldots $ form a tower of recollement $TL_{\bullet}$
(provided $\delta\in k\setminus \{0\}$;
see for example \cite[1.2]{\cmpx}).

\begin{figure}[htbp]
\[
\includegraphics{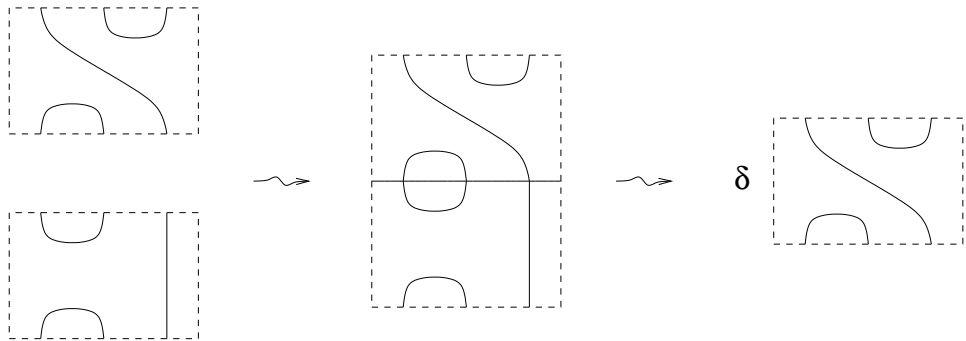}
\]
\caption{Example of multiplication in $TL_3$.}
\label{f:TLcomposition}
\end{figure}

These algebras may be defined as follows
(see \cite{Martin91}, and cf. \cite{TemperleyLieb71}
for the original definition by presentation). 
Suppose that $\delta\in k$. 
The Temperley-Lieb algebra $TL_n$ 
over $k$ has the set $\TLD(n,n)$ of diagrams as basis. 
To multiply two basis elements, the diagrams are combined, with the southern
edge of the first meeting the northern edge of the second. 
If this process creates
any loops (not connected to a vertex), they are removed and the resulting basis
element is multiplied by $\delta$ for each one. (We remark that the order
in which such loops are removed does not affect the result, which is
required for this to be well-defined; see~\cite{Bergman}).
See Figure~\ref{f:TLcomposition} for an example.


Note that the number of propagating lines cannot increase in
composition. 
This means that 
the set $\TLD_l(n,n)$ is a basis for a $TL_n$-module, where the action
is algebra multiplication, modulo diagrams with fewer propagating
lines.
Now recall that each diagram $D\in \TLD(n,n)$ 
can be expressed in the bra-ket decomposition
(\ref{ppm01}). Write 
\[
D = | D \rangle \; \langle D |
\] 
for this.

{\lemma{\label{l:TL}\cite{Martin91}
 If $DD' = \delta^r D''$ and $D'' \in D_l(n,n)$ then 
$D'' =  | D \rangle \; \langle D' |$ 
(i.e. it is independent of  $ | D' \rangle, \; \langle D |$). \Qed
}}

It follows that the module  $k\TLD_l(n,n)$ above decomposes as
a sum of isomorphic left submodules:
\[
k\TLD_l(n,n) \cong \bigoplus_{D \in \TLD_l(l,n)} k | \TLD_l(n,l) \rangle \; \langle D |
\]
In particular, in 
each of these submodules the right-hand part of the component
($|D\rangle$) is irrelevant for the action of the Temperley-Lieb
algebra.
Thus $\TLD_l(n,l)$ can be considered as a basis for a left
$TL_n$-module. 
Varying $l$ over the possible numbers of propagating lines, 
these will turn out below to give the set of standard modules,
so we will denote these modules $\Delta_n(l)$. 

By Lemma~\ref{l:TL} there is a well-defined inner product 
$ \langle D | | D' \rangle $
on $\Delta_{n}(l) = k D_{l}(n,l)$ 
given by  
\eql(gram1)
 \langle D | | D' \rangle = \left\{ \begin{array}{ll}
\delta^r
& \mbox{in the case in Lemma~\ref{l:TL} with $D'' \in D_l(n,n)$ } 
\\
0 
& \mbox{
if $D'' \not\in  D_l(n,n)$}
\end{array} \right.
\eq
(see \cite{Martin91} for details).


Next we consider the tower structure. 
For $TL_n$ we have the indexing set
$$
\Lambda_n=\{n,n-2,n-4,\ldots ,0\mbox{\ or\ } 1\}.
$$


By~\cite[1.2]{\cmpx} we know that axioms (A1) to (A5)(a) hold.
An argument similar to that for axiom (A4) shows that (A4') holds in
this case also, and it is clear that (N5) holds. 
By Proposition~\ref{ToR-GT}
it follows that the \TL\ algebras (with the above set of
modules) form a $(G,0)$-tower of algebras for some rooted graph $(G,0)$.
Since we have the short exact sequences
$$
0\rightarrow \Delta_{n-1}(l-1) \rightarrow res \Delta_n(l) \rightarrow
\Delta_{n-1}(l+1)\rightarrow 0
$$
for $0\leq l<n$ and $res\Delta_n(n)\cong \Delta_{n-1}(n-1)$
(see~\cite{Martin91,CoxGrahamMartin03} for details), it follows from the
definition
of $G$ (see Definition~\ref{graphG}) that $(G,0)=(A_{\infty},0)$,
and hence that the array associated to the tower of Temperley-Lieb
algebras (see Definition~\ref{d:YA}),
is a Pascal $(A_{\infty},0)$-array of sets and thus equivalent to
the array $Y_{TL}$ considered in Section~\ref{s:TLdiagrams}. We see
that the basis of $\Delta_n(\lambda)$ provided by Proposition~\ref{towerbases} 
can be explicitly parametrized by the half diagrams $D_l(n,l)$,
recovering the above description of the modules $\Delta_n(l)$.

Remark:
 if all the inner products (\ref{gram1}) are non-degenerate (which
will depend on $\delta$), the algebra has a semisimple multimatrix
structure as in (\ref{sumLL}).
This inner product structure 
is the original reason for the  bra-ket nomenclature.
{{\rem
Given a rooted directed graph satisfying suitable
axioms, \cite[p52]{GoodmanDelaharpeJones89}
describes a tower of algebras with that graph as its
Bratteli diagram (due to
V.~S.~Sunder \cite{Sunder87} and A.~Ocneanu \cite{Ocneanu86}). 
It follows from this construction that, given any rooted
directed graph $(G,0)$, there is a tower of algebras whose Bratteli
diagram is the directed cover of $(G,0)$ (see
Definition~\ref{d:directedcover}).
}}

\section{A first generalisation: 
$(A^{\infty}_{\infty},0)$ and the blob algebra}
\label{s:blob}

A {\em blob diagram} is a \TL\ diagram 
with the additional possibility of decorating any arc 
exposed to the western end of the diagram with a single blob
\cite{MartinSaleur94a,MartinWoodcock2000}. 
We denote the set of blob diagrams with $n$ vertices on the north edge
and $m$ on the south edge by $\TLD^b(n,m)'$. 
For example, 
let $e$ denote the diagram in $\TLD^b(n,n)'$ with $n$ propagating lines and the
westernmost line decorated with a blob.
For $l\geq 0$ 
we write $\TLD^b_l(n,m)'$  for the subset of blob diagrams
with $l$ propagating lines, one of which (if $l>0$) is decorated;
and  write $\TLD^b_{-l}(n,m)'$  for the subset
with $l$ propagating lines, all undecorated.


For $\delta,\delta' \in k$, 
the sequence of blob algebras~\cite{MartinSaleur94a}
$b_{\bullet}= \{ b_n(\delta,\delta') \}_{n=0,1,2,\ldots} $ 
is a sequence of $k$-algebras defined as follows. 
A diagram basis of $b_n$ is  $\TLD^b(n,n)'$.
If the process of composition as in the \TL\ case
creates loops then
these are again removed: 
each undecorated loop is replaced with the scalar $\delta$; 
each decorated loop is replaced with the scalar $\delta'$. 
Any arc with two blobs on it is replaced
with the same arc with a single blob. 
It will be evident that this composition is again associative and unital.


We may retain the bra-ket construction as in 
(1) from Section~\ref{s:TLdiagrams}, with the
proviso that when starting with a propagating 
arc with a single blob on it, this is first
replaced with two blobs; one each is then distributed to the 
bra and ket parts of the diagram. 
The half diagrams obtained from $\TLD^b(n,n)'$ with $l$
propagating lines thus lie in $\TLD^b_{\pm l}(n,l)'$.
The constructions in (2) in Section~\ref{s:TLdiagrams} go through
unchanged.


As in~\cite{MartinSaleur94a}, for symmetry we introduce a generator $f=1-e$. 
As a shorthand in the diagram calculus
we may represent $f$ as a propagating line carrying a square box. 
Two boxes on an arc can thus again be replaced by one. 
Now let $\TLD^b(n,m)$ be a set consisting of the same tangles
as  $\TLD^b(n,m)'$, except that
every arc exposed to the western end of the rectangle carries 
either a blob or a square. 
Then 
another basis of $b_n$ is 
 $\TLD^b(n,n)$. 
For $l\in\{n,n-2,n-4,\ldots ,2-n,-n\}$, we denote by
$\TLD^b_l(n,l)$ the set of half diagrams in $\TLD^b(n,l)$ with $|l|$ propagating
lines including one propagating line with a blob if $l>0$, respectively one
propagating line with a square if $l<0$.

The sequence $b_{\bullet}$ is  known to form a tower of 
recollement~\cite{CoxMartinParkerXi06} for suitable $\delta,\delta'$.

In this case the index set for simples  
of $b_n$ 
is given by $\Lambda_n=\{n,n-2,n-4,\ldots, 2-n,-n\}$.
Let $l\in \Lambda_n$. A basis for the standard module
$\Delta^b_n(l)$ is given by the set 
$\TLD^b(n,l)$.
See~\cite[Proposition 2]{MartinSaleur94a}.

\begin{figure}[htbp]
{
\newcommand{\li}
{\begin{picture}(10,10)(0,0)
\put(0,-3){\line(0,1){10}}
\end{picture}}
\newcommand{\sq}{\begin{picture}(10,10)(0,0)
\put(0,-3){\line(0,1){10}}
\put(-2,2){${}_{\Box}$}
\end{picture}}
\newcommand{\bl}
{\begin{picture}(10,10)(0,0)
\put(0,-3){\line(0,1){10}}
\put(-2,0){$\bullet$}
\end{picture}}
\newcommand{\dbl}
{\begin{picture}(10,15)(0,0)
\put(0,-3){\line(0,1){15}}
\put(-2,0){$\bullet$}
\put(-2,4){$\bullet$}
\end{picture}}
\newcommand{\dsq}
{\begin{picture}(10,15)(0,0)
\put(0,-3){\line(0,1){15}}
\put(-2,1){${}_\Box$}
\put(-2,7){${}_\Box$}
\end{picture}}
\newcommand{\cus}{\bigcup_{\!\!\!\! \Box}}
\newcommand{\cub}{\bigcup_{\!\!\!\! \bullet}}
\newcommand{\sms}{
\smile \!\! \circ}
\newcommand{\smb}{\smile \!\!
\bullet}
\newcommand{\rec}{\begin{picture}(0,0)(0,0)
\put(-15,60){\dashbox{1.0}(140,40){}}
\end{picture}}
\renewcommand{\Cup}{
\begin{picture}(20,15)(-7,-7)
\put(0,0){\oval(10,15)[b]}
\end{picture}}
\newcommand{\Cus}{
\begin{picture}(20,20)(-7,-7)
\put(0,0){\oval(10,15)[b]}
\put(-3,-8){${}_{\Box}$}
\end{picture}}
\newcommand{\Cub}{
\begin{picture}(20,20)(-7,-7)
\put(0,0){\oval(10,15)[b]}
\put(-3,-10){$\bullet$}
\end{picture}}
\newcommand{\scus}{
\begin{picture}(20,20)(-7,-7)
\put(0,0){\oval(10,10)[b]}
\put(0,0){\oval(30,20)[b]}
\put(0,-10){${}_{\Box}$}
\end{picture}}
\newcommand{\scup}{
\begin{picture}(40,15)(-15,-7)
\put(0,0){\oval(10,10)[b]}
\put(0,0){\oval(30,20)[b]}
\end{picture}}
\newcommand{\scub}{
\begin{picture}(20,20)(-7,-7)
\put(0,0){\oval(10,10)[b]}
\put(0,0){\oval(30,20)[b]}
\put(0,-12){$\bullet$}
\end{picture}}
\[
\begin{array}{cccccccccccc}
&&&&\sq&&\bl
\\ \\
&&&\sq |&&\Cus&&\bl |
\\
&&&     &&\Cub
\\ \\
&&\sq \li | && \sq \Cup && \Cus \bl && \bl \li |
\\
&&        && \Cus \sq    && \Cub \bl
\\
&&        && \Cub \sq    && \bl \Cup
\\ \\
&\sq \li \li |&&\sq \li \Cup && \scus    
                       && \Cus \bl | && \bl \li \li |
\\
&       &&\sq \Cup | && \Cus \Cus && \Cub \bl |
\\
&       &&\Cus \sq | && \Cub \Cus && \bl \Cup |
\\
&       &&\Cub \sq | && \Cus \Cub && \bl \li \Cup
\\
&       &&           && \Cub \Cub
\\
&        &&          && \scub
\\
\end{array}
\]
}
\caption{The start of the array of blob diagrams}
\label{f:blobdiagrams}
\end{figure}
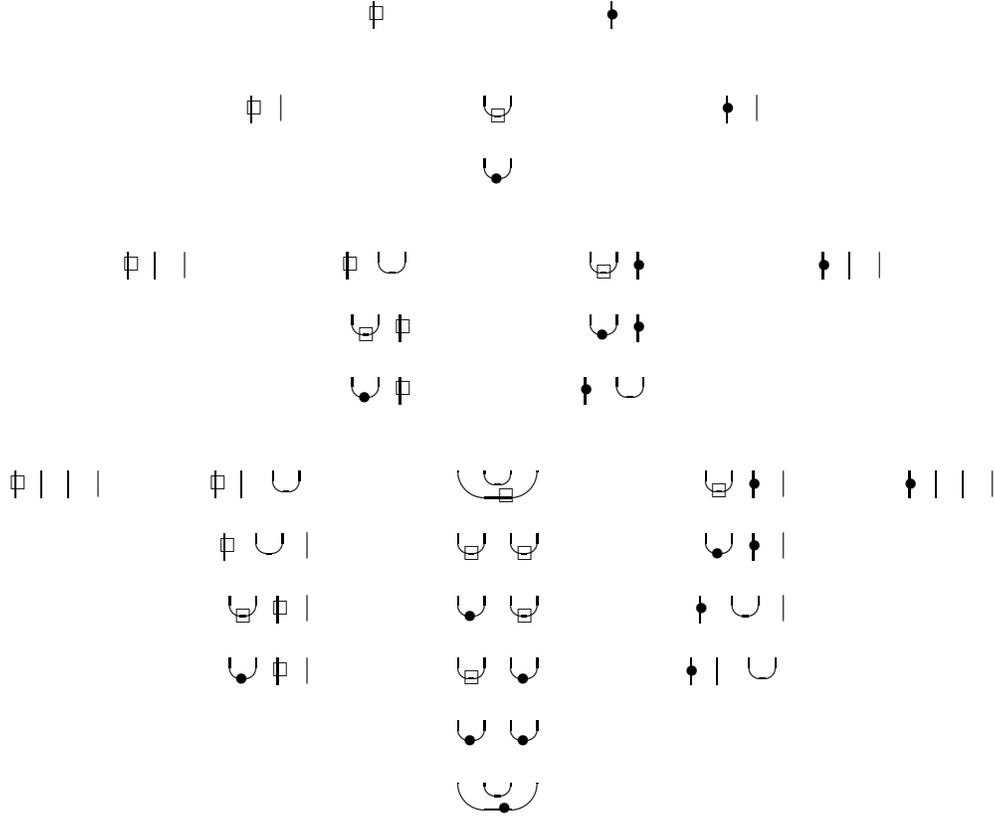


\begin{prop}
The array of sets 
$((\TLD^b_l(n,l))_{l\in A_{\infty}^{\infty}})_n$ 
is a Pascal $(A^{\infty}_{\infty},0)$--array.
(a) If $i>0$, the edge maps $\edge{i,i+1}$ are given
by $\phi^1$ and the edge maps $\edge{i,i-1}$ are given by $\phi^u$
(see Section~\ref{s:TLdiagrams} for the definitions). \\
(b) If $i=0$, the edge maps $\edge{i,i+1}$ are given by
$\phi^1$ together with a blob on the new propagating line, and the
edge maps $\edge{i,i-1}$ are given by $\phi^1$ together with a square on
the new propagating line. \\
(c) If $i<0$, the edge maps $\edge{i,i+1}$ are given
by $\phi^u$ and the edge maps $\edge{i,i-1}$ are given by $\phi^1$.
\\
This array underlies the basis sequence $(\TLD^b(n,n))_n$. 
\end{prop}

{\bf Proof:} It suffices to note that \\
(a) If $l>0$ then
$$
\TLD^b_l(n,l) 
\; = \; \phi^u( \TLD^b_{l+1}(n-1,l+1) ) \; \bigsqcup \; \phi^1(
\TLD^b_{l-1}(n-1,l-1) ).$$
(b) If $l=0$ then
$$\TLD^b_l(n,0) 
\; = \; \phi^u( \TLD^b_{+1}(n-1,1) ) \; 
   \bigsqcup \; \phi^u(\TLD^b_{-1}(n-1,1) ).$$

(c) If $l<0$ then
$$
\TLD^b_l(n,l) 
\; = \; \phi^u( \TLD^b_{l-1}(n-1,l-1) ) \; \bigsqcup \; \phi^1(
\TLD^b_{l+1}(n-1,l+1) ).$$
\Qed

The same algebraic interpretation of equation~\eqref{N=sumN2} exists in this
case (cf equation~\eqref{dimTL}).
The Pascal array
$((\TLD^b_l(n,l))_{l\in A_{\infty}^{\infty}})_n$ 
appears in~\cite{MartinSaleur94a}. We reproduce it in
Figure~\ref{f:blobdiagrams} for the convenience of the reader.

Note that the array of cardinalities for this Pascal array is not the
Catalan triangle but the ordinary Pascal triangle. 
(This example is the reason for our use of the term in general.)

In Section~\ref{s:furcated_graphs} we extend this generalisation. 

\section{Pascal array examples: Rooted trees} \label{s:furcated_graphs}
\newcommand{\groot}{{\mathcal A }}%
\newcommand{\dirc}{{\mathcal G }}%

Our aim in this section is to provide a large class of examples of
Pascal arrays associated to more general rooted graphs
(still trees, at this stage), 
while at the same time
indicating towers of algebras which have equivalent Pascal arrays.
We start with some notation for specifying rooted trees
which allows us to extend the generalisation in Section~\ref{s:blob}
(i.e. that from $A_{\infty}$ to $A^{\infty}_{\infty}$) 
in an appropriate way. 


\del(d:directedcover)
The \emph{directed cover}
of a rooted directed graph $(\Gamma,0)$ is the
directed graph $\dirc(\Gamma)$ with vertex set given by paths on $\Gamma$
starting at $0$. There is an edge from path $p$ to $p'$ provided that
$p=(e_1,e_2,\ldots e_k)$ and $p'=(e_1,e_2,\ldots ,e_k,e)$ for some edge
$e$ of $\Gamma$. 
The root of $\dirc(\Gamma)$ is the trivial path at $0$,
also denoted $0$. 
\end{de}

\begin{figure}[htbp]
\[
\includegraphics{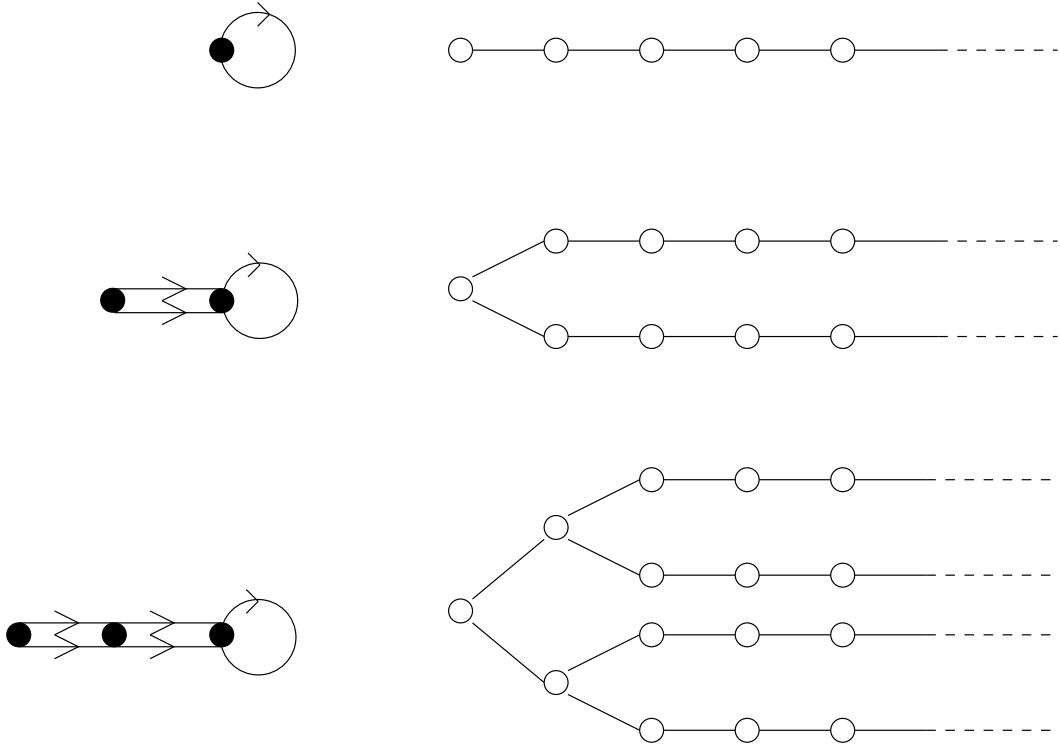}
\]
\caption{Some graphs and the underlying unoriented graphs of 
their directed covers.
The graphs on the right are the rooted graphs 
$\groot(1)=A_{\infty}$, $\groot(2,1)=A_{\infty}^{\infty}$ 
and $\groot(2,2,1)$. 
These are the underlying unoriented graphs of the
directed covers of the corresponding graphs on the left.}
\label{f:furcated_graphs}
\end{figure}

Figure~\ref{f:furcated_graphs} shows some
examples of graphs $\Gamma$ and the underlying unoriented graph
of their corresponding directed covers.


\begin{de}
For $\lambda=(\lambda_1,\lambda_2,\ldots )$ a sequence of natural
numbers,  the graph $\Gamma(\lambda)$ 
is defined to have vertices $0,1,2,\ldots $
and $\lambda_i$ directed edges from $i-1$ to $i$ for all $i=1,2,\ldots$.

If $\lambda=(\lambda_1,\lambda_2,\ldots ,\lambda_k)$ is finite 
and the last entry $\lambda_k$ is positive or zero 
then $\Gamma(\lambda)$ is the graph with vertices $0,1,\ldots ,k-1$, 
and $\lambda_i$ directed edges from $i-1$ to $i$ for $i=1,\ldots ,k-1$, 
and $\lambda_k$ directed edges from $k-1$ to $k-1$.
(See the left hand side of Figure~\ref{f:furcated_graphs} for some examples.)

The {\em rooted $\lambda$-ary tree} $\groot(\lambda)$ 
is the underlying unoriented graph of the
directed cover $\dirc(\Gamma(\lambda))$ of $\Gamma(\lambda)$. 
\end{de}


Note that 
$\groot(\lambda)$ is a tree with a root and one (undirected)
edge from the root to each of $\lambda_1$ children; each child has
this single edge to the parent, then one (undirected)
edge from it to each of $\lambda_2$ of its own children; and so on. 
If $\lambda=(\lambda_1,\lambda_2,\ldots ,\lambda_k)$ is finite 
and the last entry $\lambda_k$ in $\lambda$ is positive 
then the sequence of multiplicities of children in 
this construction of $\groot(\lambda)$
is obtained by treating $\lambda_k$ as infinitely repeating 
(i.e. $\lambda=(1)$ is equivalent to $\lambda=(1,1,1,...)$). 
On the other hand, 
if the last entry in $\lambda$ 
is $\lambda_k=0$ then the graph $\groot(\lambda)$ terminates 
at the $k$-th generation of children.   


Examples: 
$\groot((1))=A_{\infty}$,  $\groot((2,1))=A_{\infty}^{\infty}$, 
and  $\groot((2,2,1))$ is the last graph illustrated in 
Figure~\ref{f:furcated_graphs}. 

By ordering the children of each vertex, we can label the
vertices of $\groot(\lambda)$ 
in the $j$-th generation (any $j$) 
by sequences of the form
$(d_1,d_2,\ldots,d_j)$, where $1\leq d_i\leq \lambda_i$ for all $i$.
{
\example{{ \label{ex:vlabel}
The vertex labels for vertices in the $j$-th generation
in  $\groot(2,1)$ are of the form $(1,1,1,\ldots,1)$
or $(2,1,1,\ldots,1)$, corresponding to the two branches. 
}}}

Remark:
Our definition here was 
partly motivated by the discussion in~\cite{NagnibedaWoess01}
of graphs with \emph{finitely many cone types} (introduced in~\cite{Cannon84}).
The graphs $\groot(\lambda_1,\lambda_2,\ldots ,\lambda_k)$ all have finitely many cone types.


We now consider $Y_{(\groot(\lambda),0)}$ 
and various equivalent arrays,
in case $\lambda=(d^k,1)$. 
The
Catalan $(A_{\infty},0)$-sequences 
in Sections~\ref{s:TLdiagrams},~\ref{ss:Cb} and~\ref{ss:Ct} have one (or more)
generalisation corresponding to each of the graphs
generalising $(A_{\infty},0) = (\groot((1)),0)$ here. 
In particular they all have sequences of algebras associated to them.
We give a description of the combinatorial objects below, as well as
a brief 
indication of the corresponding algebras.

\Subsection{Parentheses\label{s:parenthesesgeneral}}
The set $\CC_b(n)$ of properly nested bracket sequences 
(Section~\ref{ss:Cb}) can be generalised,
for example, by
allowing $d$ different kinds of bracket pairs. If the choice for each pair
of brackets is free then the set of matched bracket sequences
containing $2n$ brackets simply has degree $d^n C(n)$
(choose a basic matched bracket sequence, then choose the type of
bracket for each pair). 


On the other hand suppose that only the brackets in the outermost
layer of nesting can be chosen freely, so that the `internal' brackets
must be of the same type as their external brackets.
We have, for example, with $d=2$ and $n=2$:
\[
\{ (()), ()(), [](), ()[], [][], [[]] \}.
\]
Write  $\{\CC_{b,(2,1)}(n)\}_n$ for the  sequence of such sets.
As with the blob diagrams from Section~\ref{s:blob}, 
the sequence $\{\CC_{b,(2,1)}(n)\}_n$ 
has an underlying Pascal  $Y_{(\groot(2,1),0)}$ array
(i.e. 
Pascal's triangle itself). 
This can be seen by constructing the Pascal array,
which we denote $Y_{b,(2,1)}$. 
Noting Example~\ref{ex:vlabel}, we may more simply
label the vertices of $\groot(2,1)$ ($ = A^{\infty}_{\infty}$) from
$\mathbb{Z}$, so that it has a positive branch and a negative branch
(cf. the picture in Figure~\ref{f:furcated_graphs}). 
For $l\in\mathbb{Z}$ then, 
the set $Y_{b,(2,1)}(n;l)$ is given by
the set of sequences of $n$ brackets with $|l|$ more open-brackets than
close-brackets with the additional property that once a bracket is
opened (from a choice of two types), 
any subsequent brackets must be of the same type, until it is closed,
or the end of the sequence is reached. 
The sign of $l$ is (say) $+$ if the open brackets are round, and $-$
if they are square. 
\newline
(1) Decomposition:
Evidently separating a sequence into two $n$-component subsequences
gives 
\[
\CC_{b,(2,1)}(n) \rightarrow 
    \cup_l Y_{b,(2,1)}(n;l) \times  Y_{b,(2,1)}^o(n;l)
\]
(where $ Y_{b,(2,1)}^o(n;l)$ is the set of opposite sequences).
\newline
(2) Edge maps:
For $i\geq 0$, the map
$\edge{i,i+1}$ is given by adding a round open-bracket and the map
$\edge{i+1,i}$ is given by adding a round close-bracket.
For $i\leq 0$, the map
$\edge{i,i-1}$ is given by adding a square open-bracket and the map
$\edge{i-1,i}$ is given by adding a square close-bracket.

\newcommand{\sst}[1]{ \left\{ \begin{array}{c} #1 
                                       \end{array} \right\} }

\begin{figure}
\xymatrix@R=5pt{ \emptyset \\ 
 & { \begin{array}{cc} \sst{(} \\ & \sst{[} \end{array} } \\
{ \mas () \\ {} [] \sam } 
 && { \begin{array}{cc} {\sst{ (( }} \\ & {\sst{ [[ }} \end{array} }
 \\
 &  { \begin{array}{cc} \mas ()( \\ {} []( \\ (() \sam \\ 
                            &  \mas ()[ \\ {}[][ \\ {} [[] \sam \end{array} } 
   && { \begin{array}{cc} {\sst{ ((( }} \\ & {\sst{ [[[ }} \end{array} }  
 \\ 
{ \begin{array}{c} \mas ()() \\ (()) \\ {} []() \\
                           {}  ()[] \\ {}[][] \\ {} [[]] \sam \end{array} } 
 && {\ldots}
}
\caption{\label{b21} The array $Y_{b,(2,1)}$ 
with layer vertices from $\groot(2,1)$ arranged as in 
Figure~\ref{f:furcated_graphs}.}
\end{figure}
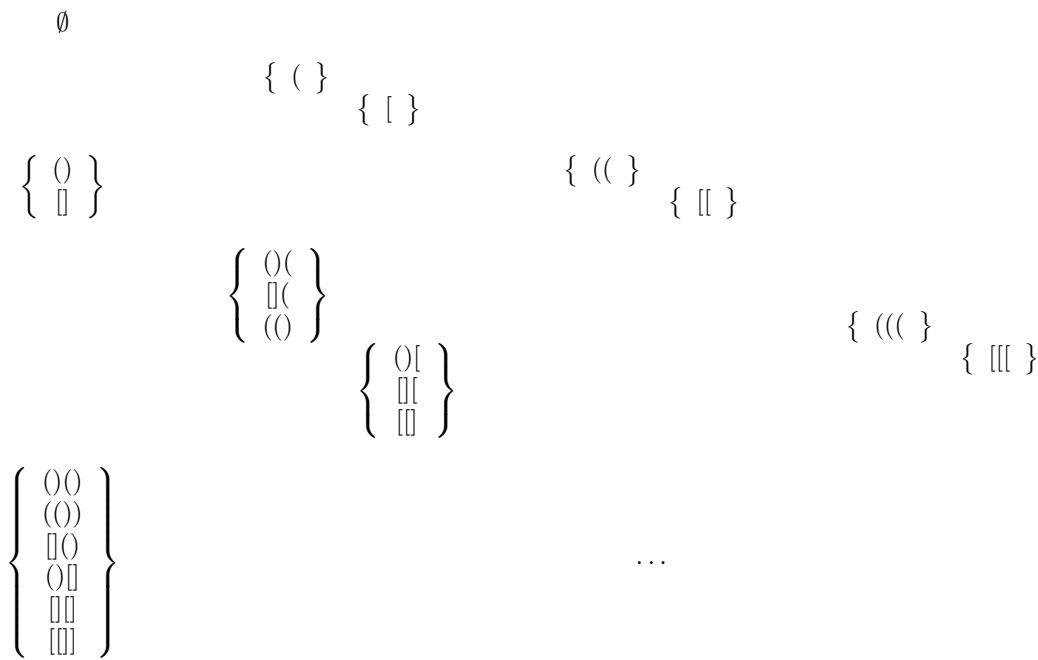

Altogether the array begins as in Figure~\ref{b21}. 
That is, it is a form of the ordinary Pascal triangle,
which is the Pascal array for $(A_{\infty}^{\infty},0)$ folded over 
(purely to match the way $A_{\infty}^{\infty}$ is drawn in
Figure~\ref{f:furcated_graphs})
onto its right-hand side. 
The corresponding Catalan sequence is of course 1, 2, 6, 20, \ldots .

More generally, given $\lambda=(\lambda_1,\lambda_2,\ldots )$ as at the
start of this section, 
consider the case where there are $\lambda_1$
possible types of brackets for the first layer of nesting, $\lambda_2$
for the second, and so on. For $\lambda=(d^k,1)$ this corresponds to
the case in which only the brackets in the outermost $k$ layers of nesting
can be chosen freely; beyond that point all the brackets must be of
the same type as those in the $k$th layer of nesting.
Write $\CC_{b,\lambda}$ for the Catalan sequence, and 
 $Y_{b,\lambda}$ for the Pascal array in this case. 

\begin{prop}
The Pascal array $Y_{b,\lambda}$ 
is equivalent to $Y_{(\groot(\lambda),0)}$.
\end{prop}
{\em Proof:}
(1)
For a vertex $(d_1,d_2,\ldots ,d_j)$ of $\groot(\lambda)$,
the set $Y_{b,\lambda}(n;d_1,d_2,\ldots ,d_j)$ is given by
the set of sequences of $n$ brackets (with the above restriction on
types of brackets) with $j$ unmatched open brackets
of types $d_1,d_2,\ldots ,d_j$ in order from left to right.
\newline
(2)
The edge maps corresponding to moving along an edge away from the root
are given by adding an open-bracket corresponding to the branch taken,
and the edge map corresponding to moving along an edge back towards the
root is given by adding a close-bracket of the same type as the open
bracket added on going out.
\Qed

\begin{figure}
\[
\includegraphics{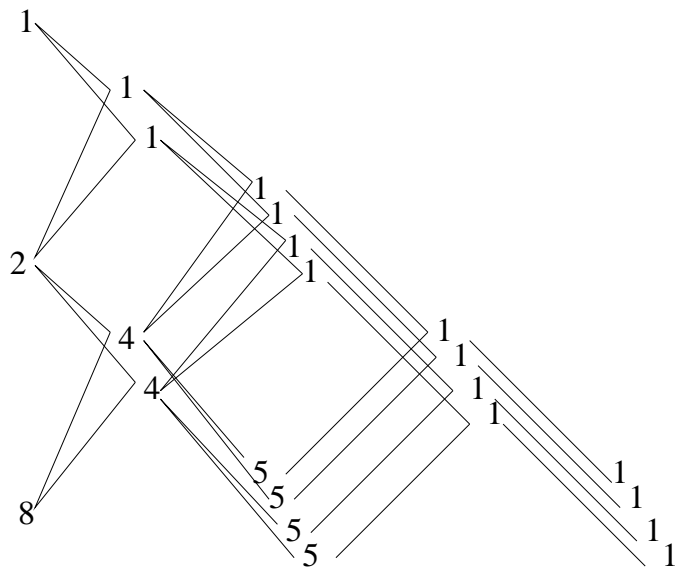}
\]
\caption{\label{furc22} Cardinalities in the array $Y_{b,(2,2,1)}$. }
\end{figure}

For example, in the case $d=2$, $\lambda=(2,2,1)$, the Pascal array
corresponds to the last graph shown explicitly in the list in
Figure~\ref{f:furcated_graphs}. The array of cardinalities begins
as in Figure~\ref{furc22}, 
so that the Catalan  numbers are 
\[
1,2,8,36,168,796,...
\] 
In Appendix~\ref{s:generatingfunctions}
we use the corresponding generalisation of the rooted tree
combinatorics to determine a generating function for this sequence.
\footnote{This sequence has one unrelated incarnation in 
Sloane \cite{sloane1}: A084868.}


{{\rem
We remark that if we have the rule that after certain sequences of
open-brackets which have not been closed, we are not allowed to make the
bracket nest any deeper, we obtain a Pascal array corresponding to a
truncated version of the $\lambda$-ary tree, in which the corresponding
branches of the tree have been removed.
We will return to this point in Section~\ref{s:truncation}. 
}}


{{\rem
We note that there is bijection between the matched bracket sequences
of length $2n$ in $C_{b,(d,d,d,...)}$ discussed at the beginning of this
section and words of length $2n$ in the usual generators of the Cuntz
algebra which evaluate to 1 
\cite[\S2.2]{Matsumoto06}. The above provides an
underlying Pascal array for the sequence of sets of such words. In this
case the bra-ket extraction amounts to writing such a word as a
concatenation of two words, expressions for each of a pair of
elements whose product is the identity.
}}



\Subsection{Contour algebra diagrams}
A \TL\ diagram can be thought of as a partition of the rectangular
interval of the plane bounded by its frame. 
A line segment in a \TL\ diagram is {\em 0-covered} if it forms
part of the boundary of the part of the rectangular interval containing the
western end of the frame of the diagram. 
A line segment in a \TL\ diagram is {\em $k$-covered} if it forms part of the
boundary of a part of the rectangular interval also bounded by 
an $(k-1)$-covered line and it is not $k'$-covered for any $k'<k$.

{\de{ 
\label{d:cont alg diagram}
A level-$k$ $d$-colour contour algebra diagram 
is a \TL\ diagram in which a line
may be decorated (with up to $d-1$  `blobs') 
if it is $k'$-covered for some $k'<k$. 
These diagrams may be constructed analogously to the \TL\ 
diagrams (which are the level-$0$ $1$-colour case).
}}


\begin{prop}
The Pascal array for the sequence of sets $\CC_{k,d}(n)$ of level-$k$
$d$-colour contour diagrams with total $2n$ vertices is 
equivalent to that of the rooted graph $(\groot(d^k,1),0)$. 
\end{prop}
{\em Proof}: A vertex $l$ of $\groot(d^k,1)$ is determined by
a sequence $(d_1,d_2,\ldots ,d_j)$ of integers
where $1\leq d_i\leq d$ for $i\leq k$ and $d_i=1$ for $i>k$.
The $d_i$ determine which branches are taken in a
path to $l$. So for $n\in\mathbb{N}$, we define
$Y_{TL,(d^k,1)}(n;d_1,d_2,\ldots ,d_j)$ to be the set of $d$-colour
level-$k$ half-diagrams with $2n$ vertices, $j$ propagating lines and
$d_i-1$ blobs on the $i$th propagating line from the left of the diagram
for $i=1,2,\ldots ,j$.
The edge map for stepping away from the root at
a branch point adds a propagating line to the right of the diagram,
with the apporpriate number of blobs;
at a non-branch point it adds an undecorated propagating line.
The edge map for stepping towards the root bends over
the rightmost propagating line in the diagram.

If a diagram in $Y_{TL,(d^k,1)}(n;d_1,d_2,\ldots ,d_j)$ has a propagating
line on the right, it is clearly in the image of a edge map for the
edge from $(n;d_1,d_2,\ldots ,d_{j-1})$. The only alternative
is that its rightmost arc $\gamma$ is non-propagating, and then it is
in the image of the edge map for the edge from
$(d_1,d_2,\ldots ,d_j,d_{j+1})$, where $d_{j+1}-1$ is the number of
blobs decorating $\gamma$,
noting the fact that $\gamma$ must be $j$-covered, so that $d_{j+1}=1$
if $j\geq k$. Since the edge maps above are clearly injective, we see that
$Y_{TL,(d^k,1)}(n;d_1,d_2,\ldots ,d_j)$ is the disjoint union of the images
of the appropriate edge maps of the edges incident with
$(n;d_1,d_2,\ldots ,d_j)$. Thus the sets of diagrams
$Y_{TL,(d^k,1)}(n;d_1,d_2,\ldots ,d_j)$ form a Pascal array of type
$(\groot(d^k,1))$ as required.

Finally, we note for the bra-ket construction that cutting a level $k$
$d$-colour contour algebra diagram into two parts horizontally gives
two diagrams in $Y_{TL,(d_k,1)}(n;d_1,d_2,\ldots ,d_j)$, where 
$d_1-1,d_2-1,\ldots,d_j-1$ 
are the numbers of blobs decorating the propagating lines counting from the
western end of
the diagram (so that a cut propagating line keeps the same number of
blobs on each of the parts). Any two such
diagrams can be stitched back together in an obvious way, and it is clear
that these two operations are inverse to each other. Thus the set of
level $k$, $d$-colour contour algebra diagrams on $2n$ vertices is equal
to the union of the sets of stitchings of pairs of diagrams from the sets
$Y_{TL,(d^k,1)}(n;d_1,d_2,\ldots ,d_j)$ over all possible vertices
$(d_1,d_2,\ldots ,d_j)$ and we are done. $\Box$


The example corresponding to the rooted graph $(\groot(2,2,2,1),0)$ has
Pascal $G$-array as illustrated in Figure~\ref{1}
(cf. \cite{MartinSaleur94a} and Figure~\ref{f:blobdiagrams}
for the corresponding figure for
$(\groot(2,1),0)$).  We remark that the combinatorics here can also be
generalised to arbitrary $\lambda$-ary graphs in the same way as the
parentheses considered in the last section, by allowing $\lambda_i$ colours
on arcs which are $(i-1)$-covered.

{\oddsidemargin -1cm
\setlength{\unitlength}{.25mm}
\newcommand{\aar}{\ar@{-}}         
\newcommand{\headroom}{80}       
\newcommand{\raised}{-\headroom} 
\newcommand{\pp}[2]{\begin{picture}(#1,\headroom)(0,\raised)
    \put(-20,0){#2} \end{picture}}
\newcommand{\ppp}[4]{\begin{picture}(#1,#2)(0,-#3) \put(-20,0){#4}
    \end{picture}}       

\newcommand{\putch}{\put(-.3,-26)}
\newcommand{\poutch}{\put(-.3,-20)}

\begin{figure}
\xymatrix%
@R=35pt@C=15pt@M=0pt%
{
{\begin{picture}(0,10)(10,30) 
\put(0,0){ {\Large $\emptyset$ } 
} 
\end{picture}}
\\
\aar[dr]+U \aar[drr]+U 
\\
 & 
\ppp{.2}{10}{-5}{
{\put(20,0){
\put(0, 0){\dashbox{1.0}(20,0){}}
\thicklines
\put(10, 0){\line(0,-1){20}}
}}
} 
& 
\ppp{.30}{10}{-5}{
{\put(20,0){
\put(0, 0){\dashbox{1.0}(20,0){}}
\thicklines
\put(10, 0){\line(0,-1){20}}
\put(10,-10.0){\circle*{4}}
}}
} 
\\
&\aar[dl]+U \aar[drr]+U \aar[drrr]+U
&\aar[dll]+U \aar[drrr]+U \aar[drrrr]+U
\\
\ppp{15}{15}{5}{
{\put(20,0){
\put(0, 0){\dashbox{1.0}(30,0){}}
\thicklines
\put(15.0, 0){\oval(10,10)[b]}
}}
\poutch{
{\put(20,0){
\put(0, 0){\dashbox{1.0}(30,0){}}
\thicklines
\put(15.0, 0){\oval(10,10)[b]}
\put(15.0,-5.0){\circle*{4}}
}}
}}
&&&
\ppp{20}{10}{-5}{
{\put(20,0){
\put(0, 0){\dashbox{1.0}(30,0){}}
\thicklines
\put(10, 0){\line(0,-1){20}}
\put(20, 0){\line(0,-1){20}}
}}
}
&
\ppp{20}{10}{-5}{
{\put(20,0){
\put(0, 0){\dashbox{1.0}(30,0){}}
\thicklines
\put(10, 0){\line(0,-1){20}}
\put(20, 0){\line(0,-1){20}}
\put(20,-10.0){\circle*{4}}
}}
}
&
\ppp{20}{10}{-5}{
{\put(20,0){
\put(0, 0){\dashbox{1.0}(30,0){}}
\thicklines
\put(10, 0){\line(0,-1){20}}
\put(10,-10.0){\circle*{4}}
\put(20, 0){\line(0,-1){20}}
}}
}
&
\ppp{20}{10}{-5}{
{\put(20,0){
\put(0, 0){\dashbox{1.0}(30,0){}}
\thicklines
\put(10, 0){\line(0,-1){20}}
\put(10,-10.0){\circle*{4}}
\put(20, 0){\line(0,-1){20}}
\put(20,-10.0){\circle*{4}}
}}
}
&
\\ 
\aar[dr]+U \aar[drr]+U 
&&&\aar[dll]+U\aar[drr]+U\aar[drrr]+U
&\aar[dlll]+U\aar[drrr]+U\aar[drrrr]+U  
&\aar[dlll]+U\aar[drrrr]+U\aar[drrrrr]+U 
&\aar[dllll]+U \aar[drrrrr]+U\aar[drrrrrr]+U 
\\
\newcommand{\pppshift}[1]{{\ppp{20}{.3}{-15}{##1}}}
&
\ppp{30}{80}{67}{
{\put(20,0){
\put(0, 0){\dashbox{1.0}(40,0){}}
\thicklines
\put(15.0, 0){\oval(10,10)[b]}
\put(30, 0){\line(0,-1){20}}
}}
\putch{
{\put(20,0){
\put(0, 0){\dashbox{1.0}(40,0){}}
\thicklines
\put(15.0, 0){\oval(10,10)[b]}
\put(15.0,-5.0){\circle*{4}}
\put(30, 0){\line(0,-1){20}}
}}
\putch{
{\put(20,0){
\put(0, 0){\dashbox{1.0}(40,0){}}
\thicklines
\put(10, 0){\line(0,-1){20}}
\put(25.0, 0){\oval(10,10)[b]}
}}
\putch{
{\put(20,0){
\put(0, 0){\dashbox{1.0}(40,0){}}
\thicklines
\put(10, 0){\line(0,-1){20}}
\put(25.0, 0){\oval(10,10)[b]}
\put(25.0,-5.0){\circle*{4}}
}}
}}}}
&
\ppp{30}{70}{60}{
{\put(20,0){
\put(0, 0){\dashbox{1.0}(40,0){}}
\thicklines
\put(15.0, 0){\oval(10,10)[b]}
\put(30, 0){\line(0,-1){20}}
\put(30,-10.0){\circle*{4}}
}}
\putch{
{\put(20,0){
\put(0, 0){\dashbox{1.0}(40,0){}}
\thicklines
\put(15.0, 0){\oval(10,10)[b]}
\put(15.0,-5.0){\circle*{4}}
\put(30, 0){\line(0,-1){20}}
\put(30,-10.0){\circle*{4}}
}}
\putch{
{\put(20,0){
\put(0, 0){\dashbox{1.0}(40,0){}}
\thicklines
\put(10, 0){\line(0,-1){20}}
\put(10,-10.0){\circle*{4}}
\put(25.0, 0){\oval(10,10)[b]}
}}
\putch{
{\put(20,0){
\put(0, 0){\dashbox{1.0}(40,0){}}
\thicklines
\put(10, 0){\line(0,-1){20}}
\put(10,-10.0){\circle*{4}}
\put(25.0, 0){\oval(10,10)[b]}
\put(25.0,-5.0){\circle*{4}}
}}
}}}}
&&&
\ppp{20}{-23}{-39}{
{\put(20,0){
\put(0, 0){\dashbox{1.0}(40,0){}}
\thicklines
\put(10, 0){\line(0,-1){20}}
\put(20, 0){\line(0,-1){20}}
\put(30, 0){\line(0,-1){20}}
}}
}
&
\ppp{20}{.3}{-34}{
{\put(20,0){
\put(0, 0){\dashbox{1.0}(40,0){}}
\thicklines
\put(10, 0){\line(0,-1){20}}
\put(20, 0){\line(0,-1){20}}
\put(30, 0){\line(0,-1){20}}
\put(30,-10.0){\circle*{4}}
}}
}
&
\ppp{20}{.3}{-29}{
{\put(20,0){
\put(0, 0){\dashbox{1.0}(40,0){}}
\thicklines
\put(10, 0){\line(0,-1){20}}
\put(20, 0){\line(0,-1){20}}
\put(20,-10.0){\circle*{4}}
\put(30, 0){\line(0,-1){20}}
}}
}
&
\ppp{20}{.3}{-25}{
{\put(20,0){
\put(0, 0){\dashbox{1.0}(40,0){}}
\thicklines
\put(10, 0){\line(0,-1){20}}
\put(20, 0){\line(0,-1){20}}
\put(20,-10.0){\circle*{4}}
\put(30,-10.0){\circle*{4}}
\put(30, 0){\line(0,-1){20}}
}}
}
&
\ppp{20}{.3}{-20}{
{\put(20,0){
\put(0, 0){\dashbox{1.0}(40,0){}}
\thicklines
\put(10, 0){\line(0,-1){20}}
\put(20, 0){\line(0,-1){20}}
\put(10,-10.0){\circle*{4}}
\put(30, 0){\line(0,-1){20}}
}}
}
&
\ppp{20}{.3}{-15}{
{\put(20,0){
\put(0, 0){\dashbox{1.0}(40,0){}}
\thicklines
\put(10, 0){\line(0,-1){20}}
\put(20, 0){\line(0,-1){20}}
\put(10,-10.0){\circle*{4}}
\put(30,-10.0){\circle*{4}}
\put(30, 0){\line(0,-1){20}}
}}
}
&
\ppp{20}{.3}{-10}{
{\put(20,0){
\put(0, 0){\dashbox{1.0}(40,0){}}
\thicklines
\put(10, 0){\line(0,-1){20}}
\put(20, 0){\line(0,-1){20}}
\put(10,-10.0){\circle*{4}}
\put(20,-10.0){\circle*{4}}
\put(30, 0){\line(0,-1){20}}
}}
}
&
\ppp{20}{.3}{-5}{
{\put(20,0){
\put(0, 0){\dashbox{1.0}(40,0){}}
\thicklines
\put(10, 0){\line(0,-1){20}}
\put(20, 0){\line(0,-1){20}}
\put(10,-10.0){\circle*{4}}
\put(20,-10.0){\circle*{4}}
\put(30,-10.0){\circle*{4}}
\put(30, 0){\line(0,-1){20}}
}}
}
&
\\
&\aar[dl]+U\aar[drr]+U\aar[drrr]+U 
&\aar[dll]+U\aar[drrr]+U\aar[drrrr]+U 
&&&\aar[dll]+U &\aar[dlll]+U &\aar[dlll]+U &\aar[dllll]+U 
&\aar[dllll]+U 
&\aar[dlllll]+U 
&\aar[dlllll]+U 
&\aar[dllllll]+U 
\\
\ppp{50}{100}{90}{
{\put(20,0){
\put(0, 0){\dashbox{1.0}(50,0){}}
\thicklines
\put(15.0, 0){\oval(10,10)[b]}
\put(35.0, 0){\oval(10,10)[b]}
}}
\poutch{
{\put(20,0){
\put(0, 0){\dashbox{1.0}(50,0){}}
\thicklines
\put(15.0, 0){\oval(10,10)[b]}
\put(15.0,-5.0){\circle*{4}}
\put(35.0, 0){\oval(10,10)[b]}
}}
\poutch{
{\put(20,0){
\put(0, 0){\dashbox{1.0}(50,0){}}
\thicklines
\put(25.0, 0){\oval(30,18.0)[b]}
\put(25.0, 0){\oval(10,10)[b]}
}}
\poutch{
{\put(20,0){
\put(0, 0){\dashbox{1.0}(50,0){}}
\thicklines
\put(25.0, 0){\oval(30,18.0)[b]}
\put(25.0, 0){\oval(10,10)[b]}
\put(25.0,-5.0){\circle*{4}}
}}
\poutch{
{\put(20,0){
\put(0, 0){\dashbox{1.0}(50,0){}}
\thicklines
\put(15.0, 0){\oval(10,10)[b]}
\put(35.0, 0){\oval(10,10)[b]}
\put(35.0,-5.0){\circle*{4}}
}}
\poutch{
{\put(20,0){
\put(0, 0){\dashbox{1.0}(50,0){}}
\thicklines
\put(15.0, 0){\oval(10,10)[b]}
\put(15.0,-5.0){\circle*{4}}
\put(35.0, 0){\oval(10,10)[b]}
\put(35.0,-5.0){\circle*{4}}
}}
\poutch{
{\put(20,0){
\put(0, 0){\dashbox{1.0}(50,0){}}
\thicklines
\put(25.0, 0){\oval(30,18.0)[b]}
\put(25.0,-9.0){\circle*{4}}
\put(25.0, 0){\oval(10,10)[b]}
}}
\poutch{
{\put(20,0){
\put(0, 0){\dashbox{1.0}(50,0){}}
\thicklines
\put(25.0, 0){\oval(30,18.0)[b]}
\put(25.0,-9.0){\circle*{4}}
\put(25.0, 0){\oval(10,10)[b]}
\put(25.0,-5.0){\circle*{4}}
}}
}}}}}}}}
&&&
\ppp{40}{100}{90}{
{\put(20,0){
\put(0, 0){\dashbox{1.0}(50,0){}}
\thicklines
\put(15.0, 0){\oval(10,10)[b]}
\put(30, 0){\line(0,-1){20}}
\put(40, 0){\line(0,-1){20}}
}}
\putch{
{\put(20,0){
\put(0, 0){\dashbox{1.0}(50,0){}}
\thicklines
\put(15.0, 0){\oval(10,10)[b]}
\put(15.0,-5.0){\circle*{4}}
\put(30, 0){\line(0,-1){20}}
\put(40, 0){\line(0,-1){20}}
}}
\putch{
{\put(20,0){
\put(0, 0){\dashbox{1.0}(50,0){}}
\thicklines
\put(10, 0){\line(0,-1){20}}
\put(25.0, 0){\oval(10,10)[b]}
\put(40, 0){\line(0,-1){20}}
}}
\putch{
{\put(20,0){
\put(0, 0){\dashbox{1.0}(50,0){}}
\thicklines
\put(10, 0){\line(0,-1){20}}
\put(25.0, 0){\oval(10,10)[b]}
\put(25.0,-5.0){\circle*{4}}
\put(40, 0){\line(0,-1){20}}
}}
\putch{
{\put(20,0){
\put(0, 0){\dashbox{1.0}(50,0){}}
\thicklines
\put(10, 0){\line(0,-1){20}}
\put(20, 0){\line(0,-1){20}}
\put(35.0, 0){\oval(10,10)[b]}
}}
\putch{
{\put(20,0){
\put(0, 0){\dashbox{1.0}(50,0){}}
\thicklines
\put(10, 0){\line(0,-1){20}}
\put(20, 0){\line(0,-1){20}}
\put(35.0, 0){\oval(10,10)[b]}
\put(35.0,-5.0){\circle*{4}}
}}
}}}}}}
&
\ppp{40}{100}{90}{
{\put(20,0){
\put(0, 0){\dashbox{1.0}(50,0){}}
\thicklines
\put(15.0, 0){\oval(10,10)[b]}
\put(30, 0){\line(0,-1){20}}
\put(40, 0){\line(0,-1){20}}
\put(40,-10.0){\circle*{4}}
}}
\putch{
{\put(20,0){
\put(0, 0){\dashbox{1.0}(50,0){}}
\thicklines
\put(15.0, 0){\oval(10,10)[b]}
\put(15.0,-5.0){\circle*{4}}
\put(30, 0){\line(0,-1){20}}
\put(40, 0){\line(0,-1){20}}
\put(40,-10.0){\circle*{4}}
}}
\putch{
{\put(20,0){
\put(0, 0){\dashbox{1.0}(50,0){}}
\thicklines
\put(10, 0){\line(0,-1){20}}
\put(25.0, 0){\oval(10,10)[b]}
\put(40, 0){\line(0,-1){20}}
\put(40,-10.0){\circle*{4}}
}}
\putch{
{\put(20,0){
\put(0, 0){\dashbox{1.0}(50,0){}}
\thicklines
\put(10, 0){\line(0,-1){20}}
\put(25.0, 0){\oval(10,10)[b]}
\put(25.0,-5.0){\circle*{4}}
\put(40, 0){\line(0,-1){20}}
\put(40,-10.0){\circle*{4}}
}}
\putch{
{\put(20,0){
\put(0, 0){\dashbox{1.0}(50,0){}}
\thicklines
\put(10, 0){\line(0,-1){20}}
\put(20, 0){\line(0,-1){20}}
\put(20,-10.0){\circle*{4}}
\put(35.0, 0){\oval(10,10)[b]}
}}
\putch{
{\put(20,0){
\put(0, 0){\dashbox{1.0}(50,0){}}
\thicklines
\put(10, 0){\line(0,-1){20}}
\put(20, 0){\line(0,-1){20}}
\put(20,-10.0){\circle*{4}}
\put(35.0, 0){\oval(10,10)[b]}
\put(35.0,-5.0){\circle*{4}}
}}
}}}}}}
&
\ppp{40}{100}{90}{
{\put(20,0){
\put(0, 0){\dashbox{1.0}(50,0){}}
\thicklines
\put(15.0, 0){\oval(10,10)[b]}
\put(30, 0){\line(0,-1){20}}
\put(40, 0){\line(0,-1){20}}
\put(30,-10.0){\circle*{4}}
}}
\putch{
{\put(20,0){
\put(0, 0){\dashbox{1.0}(50,0){}}
\thicklines
\put(15.0, 0){\oval(10,10)[b]}
\put(15.0,-5.0){\circle*{4}}
\put(30, 0){\line(0,-1){20}}
\put(40, 0){\line(0,-1){20}}
\put(30,-10.0){\circle*{4}}
}}
\putch{
{\put(20,0){
\put(0, 0){\dashbox{1.0}(50,0){}}
\thicklines
\put(10, 0){\line(0,-1){20}}
\put(25.0, 0){\oval(10,10)[b]}
\put(40, 0){\line(0,-1){20}}
\put(10,-10.0){\circle*{4}}
}}
\putch{
{\put(20,0){
\put(0, 0){\dashbox{1.0}(50,0){}}
\thicklines
\put(10, 0){\line(0,-1){20}}
\put(25.0, 0){\oval(10,10)[b]}
\put(25.0,-5.0){\circle*{4}}
\put(40, 0){\line(0,-1){20}}
\put(10,-10.0){\circle*{4}}
}}
\putch{
{\put(20,0){
\put(0, 0){\dashbox{1.0}(50,0){}}
\thicklines
\put(10, 0){\line(0,-1){20}}
\put(20, 0){\line(0,-1){20}}
\put(10,-10.0){\circle*{4}}
\put(35.0, 0){\oval(10,10)[b]}
}}
\putch{
{\put(20,0){
\put(0, 0){\dashbox{1.0}(50,0){}}
\thicklines
\put(10, 0){\line(0,-1){20}}
\put(20, 0){\line(0,-1){20}}
\put(10,-10.0){\circle*{4}}
\put(35.0, 0){\oval(10,10)[b]}
\put(35.0,-5.0){\circle*{4}}
}}
}}}}}}
&
\ppp{40}{100}{90}{
{\put(20,0){
\put(0, 0){\dashbox{1.0}(50,0){}}
\thicklines
\put(15.0, 0){\oval(10,10)[b]}
\put(30, 0){\line(0,-1){20}}
\put(40, 0){\line(0,-1){20}}
\put(30,-10.0){\circle*{4}}
\put(40,-10.0){\circle*{4}}
}}
\putch{
{\put(20,0){
\put(0, 0){\dashbox{1.0}(50,0){}}
\thicklines
\put(15.0, 0){\oval(10,10)[b]}
\put(15.0,-5.0){\circle*{4}}
\put(30, 0){\line(0,-1){20}}
\put(40, 0){\line(0,-1){20}}
\put(30,-10.0){\circle*{4}}
\put(40,-10.0){\circle*{4}}
}}
\putch{
{\put(20,0){
\put(0, 0){\dashbox{1.0}(50,0){}}
\thicklines
\put(10, 0){\line(0,-1){20}}
\put(25.0, 0){\oval(10,10)[b]}
\put(40, 0){\line(0,-1){20}}
\put(10,-10.0){\circle*{4}}
\put(40,-10.0){\circle*{4}}
}}
\putch{
{\put(20,0){
\put(0, 0){\dashbox{1.0}(50,0){}}
\thicklines
\put(10, 0){\line(0,-1){20}}
\put(25.0, 0){\oval(10,10)[b]}
\put(25.0,-5.0){\circle*{4}}
\put(40, 0){\line(0,-1){20}}
\put(10,-10.0){\circle*{4}}
\put(40,-10.0){\circle*{4}}
}}
\putch{
{\put(20,0){
\put(0, 0){\dashbox{1.0}(50,0){}}
\thicklines
\put(10, 0){\line(0,-1){20}}
\put(20, 0){\line(0,-1){20}}
\put(10,-10.0){\circle*{4}}
\put(20,-10.0){\circle*{4}}
\put(35.0, 0){\oval(10,10)[b]}
}}
\putch{
{\put(20,0){
\put(0, 0){\dashbox{1.0}(50,0){}}
\thicklines
\put(10, 0){\line(0,-1){20}}
\put(20, 0){\line(0,-1){20}}
\put(10,-10.0){\circle*{4}}
\put(20,-10.0){\circle*{4}}
\put(35.0, 0){\oval(10,10)[b]}
\put(35.0,-5.0){\circle*{4}}
}}
}}}}}}
\\
\hspace{50pt}&\hspace{30pt}
}
\vspace{25pt}
\caption{\label{1} Pascal array for level-3 2-colour contour diagrams. The
figure has been truncated on the right (in the last row only) and below.}
\end{figure}
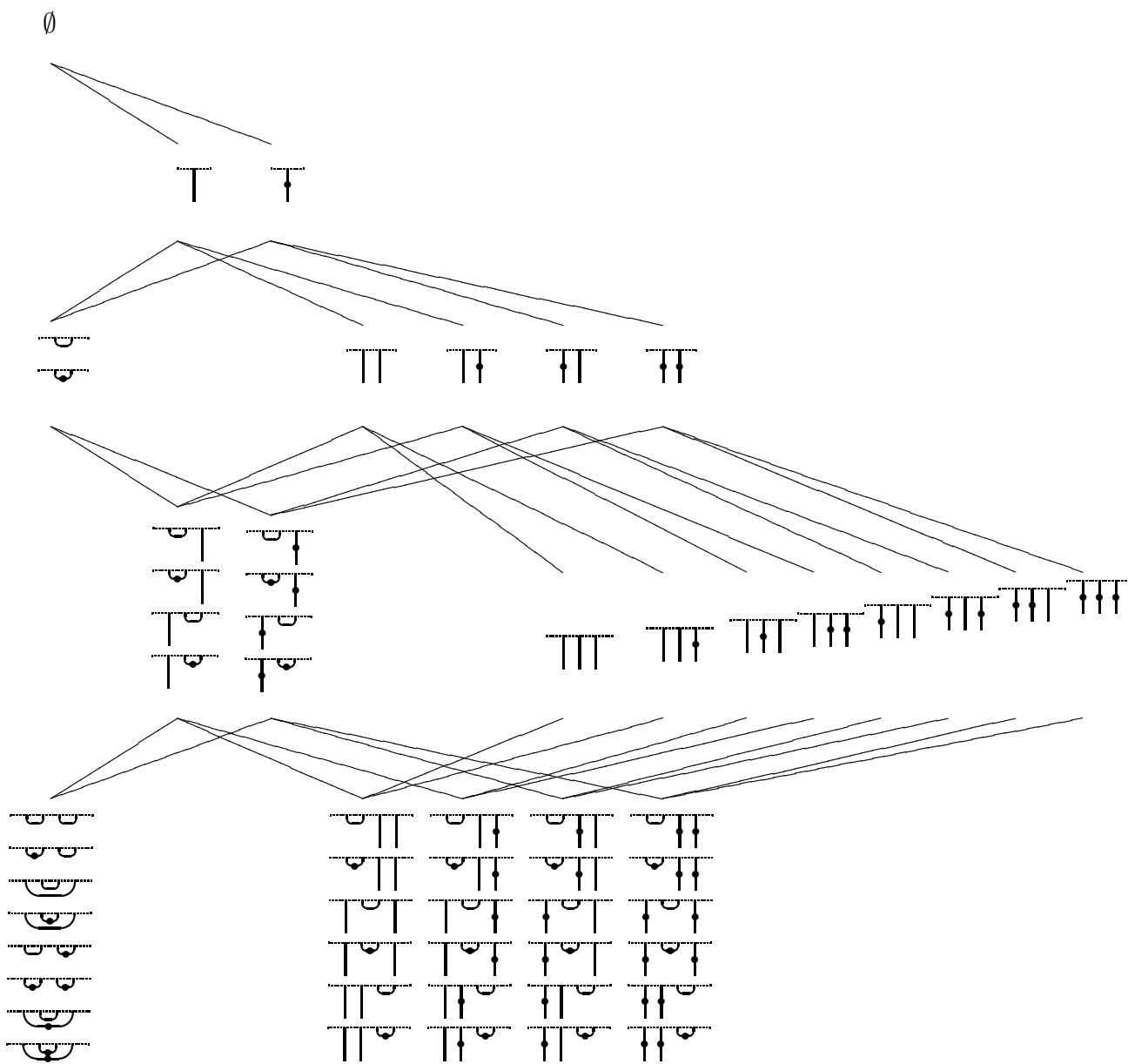
}


Just as these decorated diagrams generalise the \TL\ diagrams,
so there is a generalisation of diagram composition, giving us a
generalised \TL\ algebra with these diagrams as basis. 
It will be evident that juxtaposition of these decorated diagrams
in a manner analogous to \TL\ diagram juxtaposition 
gives another decorated diagram, 
with the underlying diagram (ignoring blobs) being
given as in Section~\ref{TLA}.
The new features are that decorated lines of the composite may have
too many blobs, and that loops may also be decorated, as in
Section~\ref{s:blob}.
Thus we need rules to reduce such pseudodiagrams into linear
combinations of legal diagrams. 
In fact contour diagrams can be seen as a basis for contour algebras
in a number of different ways, and correspondingly different sets of
rules arise. In our case it will be sufficient to proceed as follows.
Let $(\delta^0, \delta^1, \ldots, \delta^{d-1})$ be a $d$-tuple of
elements of the field $k$. Then the contour algebra 
$c_n(\delta^0, \delta^1, \ldots, \delta^{d-1})$
arises from one of the following reductions. Either
\\
(a) Cyclotomic reduction:
If a line in a pseudodiagram has $d' > d-1$ blobs on it then 
 reduce modulo $d$; Or
\\ 
(b) Blob reduction: When a line segment with a run of 
$d'$ blobs on it gets concatenated with a
 line with $d'' \geq d'$ blobs on it then 
discard the run of $d'$ blobs.
\\
Then, if a loop has $l$ blobs on it, replace it with the scalar
factor $\delta^l$. 

The case $d=2$, $k=1$ with the blob reduction is the blob algebra 
(briefly reviewed in Section~\ref{s:blob}).
For other cases see \cite{CoxMartinParkerXi06}. 

\Subsection{Forests with special spinneys}
Next we consider rooted planar trees 
in which 
some subset of the vertices at distance
one from the root (i.e. in layer $1$) are designated as a marked subset. 
Alternatively we can think of this as a generalisation in which the
vertices in layer $1$ are {\em coloured}, from a choice of two colours. 
The consequence is that, for example, there are now six distinct
trees with $n=2$ edges:
\[
\xymatrix@R=22pt{{\circ_{}} \\ {\circ_{red}} \ar[u] \\  {\circ} \ar[u]} 
\hspace{.3in} 
\xymatrix@R=22pt{{\circ_{}} \\ {\circ_{blue}} \ar[u] \\  {\circ} \ar[u]} 
\hspace{.351in} 
\xymatrix@C=6pt{\\ {\circ_{red}} && {\circ_{red}} \\ & {\circ} \ar[ul]\ar[ur]} 
\hspace{.351in} 
\xymatrix@C=5pt{\\ {\circ_{red}} && {\circ_{blue}} \\ & {\circ} \ar[ul]\ar[ur]} 
\hspace{.351in} 
\xymatrix@C=5pt{\\ {\circ_{blue}} && {\circ_{red}} \\ & {\circ} \ar[ul]\ar[ur]} 
\hspace{.351in} 
\xymatrix@C=6pt{\\ {\circ_{blue}} && {\circ_{blue}} \\ & {\circ} \ar[ul]\ar[ur]} 
\] 

\begin{figure}
\includegraphics{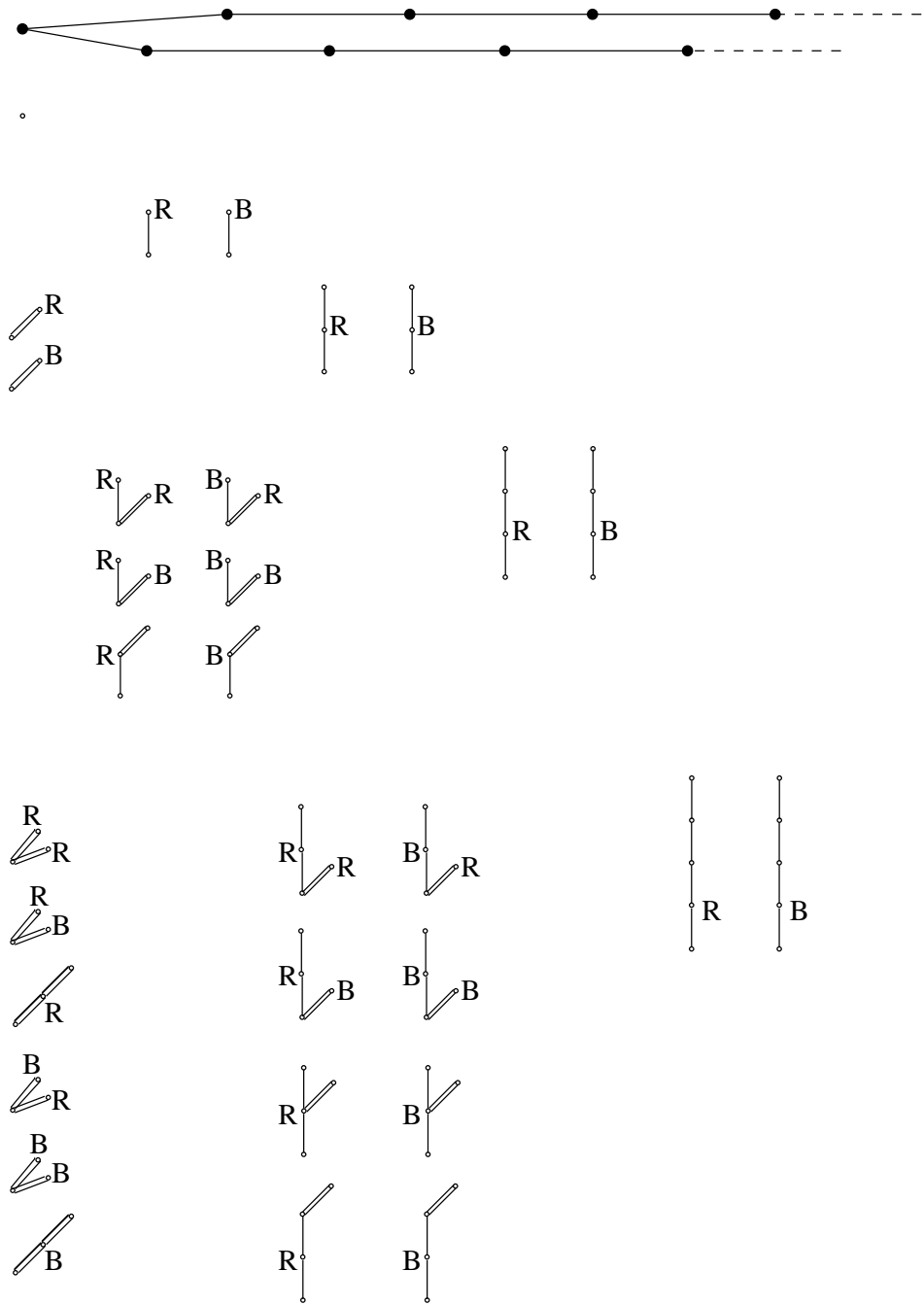}
\caption{\label{f:colouredtrees} The Pascal $(A_{\infty}^{\infty},0)$--array 
$Y_{t,(2,1)}$ of half-trees with coloured vertices in layer $1$.}
\end{figure}

For $l$ a vertex of $A_{\infty}^{\infty}$, we define $Y_{t,(2,1)}(n;l)$
to be like the set of half-trees in $Y_t(n;|l|)$,
generalised to have coloured vertices in layer $1$.
For $l\not=0$, the vertex of the trunk in layer $1$ is
coloured by the sign of $l$. For $l=0$, the trunk has only one vertex,
so no vertex in layer $1$, and there is no restriction on the colouring.
The edge maps for moving away from the root are given by adding an
edge to the trunk (labelled by the appropriate colour if the added vertex
is in layer $1$); for moving back towards the root they are given by
adding an extra edge to the top edge of the trunk (these two edges then
become part of the uppermost branch of the tree). In this way we obtain
a Pascal array of type $A_{\infty}^{\infty}$ --- see
Figure~\ref{f:colouredtrees}. The rules for splicing the
trees together are the same as for the uncoloured trees; a vertex in
layer $1$ of the trunk inherits the same colouring as the colours of the
corresponding vertices in the trunks of the half-trees.

It is clear that this set-up can be generalised to more colours and also
to allowing colours in other layers. For a vertex $(d_1,d_2,\ldots ,d_j)$ of
$\groot(\lambda)$, we define $Y_{t,\lambda}(d_1,d_2,\ldots ,d_j)$ to be
the set of half-trees from $Y_t(n;j)$ except that the vertices in the first
$j$ layers of the trunk have colours $d_1,d_2,\ldots ,d_j$ moving from
bottom to top.
The edge maps are given by adding a vertex (and corresponding edge)
of the appropriate colour to the trunk (moving away from the root) or
doubling the top edge of the trunk (moving towards the root). In this way
we obtain a Pascal array of type $\groot(\lambda)$, and we see that the
set of rooted planar trees with the vertices in layer $i$ coloured by
$\lambda_i$ colours is a Catalan sequence of sets of type $\groot(\lambda)$.

Remark: There is a bijection between
$Y_{t,\lambda}(d_1,d_2,\ldots ,d_j)$
and the set $Y_{b,\lambda}(d_1,d_2,\ldots ,d_j)$
considered in Section~\ref{s:parenthesesgeneral} given by following the
boundary tree of a half-tree anticlockwise from the root; the colours
correspond to the types of bracket. This induces a bijection on the
level of the trees themselves. 

An advantage in considering
trees is that it helps to determine the generating functions for the
corresponding sequences of generalised Catalan numbers
(see Appendix~\ref{s:generatingfunctions}).
 
\Subsection{A generalisation to locally-finite rooted trees}
We remark that the combinatorial constructions considered above can be
generalised
by replacing the trees $\groot(\lambda)$ with arbitrary locally-finite
rooted trees, i.e. trees in which each vertex has finite valency.
\section{Pascal array examples: Other rooted graphs} \label{s:reentrant}

In our examples to this point each graph $G$ has been a tree.
The purpose of this section is to show that the notion of Pascal
$G$-arrays underlying generalised Catalan sequences 
extends usefully beyond those cases in which $G$ is a tree graph.

The ultimate 
guide and source for all our examples here is representation theory.
But this is simply so that we can select graphs 
whose arrays we know in advance to have wider combinatorial interest. 

It is appropriate to begin by introducing the graphs that we use.

\begin{figure}
\includegraphics{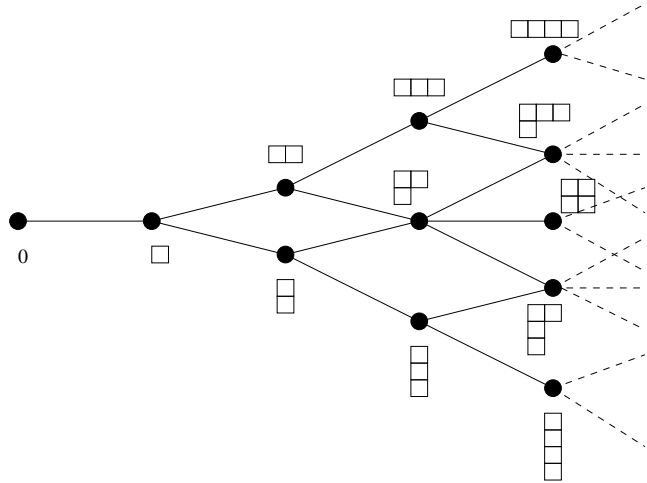}
\caption{\label{Younggraph} The beginning of the Young graph in Young
  diagram labelling.}
\end{figure}

\begin{de}
\label{Young graph}
The {\em Young graph} $\YounG$ is defined as follows.
The vertex set is the set $\Gamma$ of integer partitions,
that is, the set of integer sequences $\lambda=(\lambda_1,\lambda_2,\ldots)$
with $\lambda_i \geq \lambda_{i+1}$ for all $i$ and all but finitely many
entries zero. There is an edge from $\lambda$ to $\lambda'$ if 
$\lambda-\lambda'$ is of the form $(0,0,\ldots,\pm 1,0,0,\ldots)$. 
\end{de}
See Figure~\ref{Younggraph} for an illustration, 
with integer partitions displayed as Young diagrams.

Remarks: 
\newline (1) The Young graph is infinite, undirected, simple and loop-free. 
\newline
(2) There is an inclusion of symmetric groups 
$I:S_n \hookrightarrow S_{n+1}$ 
given by $I(p)(i)=p(i)$ for $i \leq n$ and $I(p)(n+1)=n+1$. 
The Young graph  
is the Bratteli diagram for the sequence of symmetric
group algebras $\C \subseteq \C S_1 \subset \C S_2 \subset \ldots$
under $I$,
with each $\nu\in\Gamma$ labelling an irreducible $S_{|\nu|}$-module 
$V_{\nu}$ in the usual way \cite{JamesKerber81,VershikOkunkov96}.  

{\de{ 
Let $\Gamma^+$ denote another distinct copy of $\Gamma$. 
The {\em double Young graph} $\YounGG$ is the simple undirected graph with
vertex set $\Gamma\cup\Gamma^+$. For $\lambda,\mu\in \Gamma$,
there is an edge between $\lambda$ and $\nu+$ if $\lambda=\nu$ or  
$\lambda-\nu$ is of the form $(0,0,\ldots, 1,0,0,\ldots)$. 
}}
\newline
It is clear
that the graph $\YounGG$ is bipartite, with decomposition $\Gamma\cup\Gamma^+$.
It is illustrated in Figure~\ref{Young1}.


\def\shadow{ shadow }
\def\SS(#1){{\mathcal S}_{#1}}                
\newcommand{\svtnrm}{}
\def\Gbarn{{\tilde{\Gamma}_{n}}{}}          

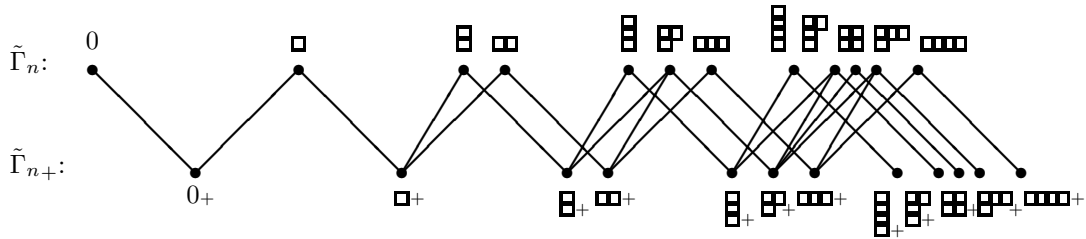
\begin{figure}

\setlength{\unitlength}{0.0054in}%
\begin{picture}(1040,220)(40,120)
\thicklines
\put(120,280){\circle*{10}}
\put(320,280){\circle*{10}}
\put(480,280){\circle*{10}}
\put(520,280){\circle*{10}}
\put(640,280){\circle*{10}}
\put(680,280){\circle*{10}}
\put(720,280){\circle*{10}}
\put(800,280){\circle*{10}}
\put(840,280){\circle*{10}}
\put(860,280){\circle*{10}}
\put(880,280){\circle*{10}}
\put(920,280){\circle*{10}}
\put(220,180){\circle*{10}}
\put(420,180){\circle*{10}}
\put(580,180){\circle*{10}}
\put(620,180){\circle*{10}}
\put(740,180){\circle*{10}}
\put(780,180){\circle*{10}}
\put(820,180){\circle*{10}}
\put(900,180){\circle*{10}}
\put(940,180){\circle*{10}}
\put(960,180){\circle*{10}}
\put(980,180){\circle*{10}}
\put(1020,180){\circle*{10}}
\put(635,300){\framebox(10,10){}}
\put(635,310){\framebox(10,10){}}
\put(635,320){\framebox(10,10){}}
\put(780,310){\framebox(10,10){}}
\put(780,320){\framebox(10,10){}}
\put(780,330){\framebox(10,10){}}
\put(780,300){\framebox(10,10){}}
\put(475,300){\framebox(10,10){}}
\put(475,310){\framebox(10,10){}}
\put(510,300){\framebox(10,10){}}
\put(520,300){\framebox(10,10){}}
\put(670,300){\framebox(10,10){}}
\put(670,310){\framebox(10,10){}}
\put(680,310){\framebox(10,10){}}
\put(705,300){\framebox(10,10){}}
\put(715,300){\framebox(10,10){}}
\put(725,300){\framebox(10,10){}}
\put(845,300){\framebox(10,10){}}
\put(845,310){\framebox(10,10){}}
\put(855,310){\framebox(10,10){}}
\put(855,300){\framebox(10,10){}}
\put(880,300){\framebox(10,10){}}
\put(880,310){\framebox(10,10){}}
\put(890,310){\framebox(10,10){}}
\put(900,310){\framebox(10,10){}}
\put(925,300){\framebox(10,10){}}
\put(935,300){\framebox(10,10){}}
\put(945,300){\framebox(10,10){}}
\put(955,300){\framebox(10,10){}}
\put(810,300){\framebox(10,10){}}
\put(810,310){\framebox(10,10){}}
\put(810,320){\framebox(10,10){}}
\put(820,320){\framebox(10,10){}}
\put(315,300){\framebox(10,10){}}
\put(120,300){\makebox(0,0)[b]{\raisebox{0pt}[0pt][0pt]{\svtnrm 0}}}
\put(575,140){\framebox(10,10){}}
\put(575,150){\framebox(10,10){}}
\put(610,150){\framebox(10,10){}}
\put(620,150){\framebox(10,10){}}
\put(735,130){\framebox(10,10){}}
\put(735,140){\framebox(10,10){}}
\put(735,150){\framebox(10,10){}}
\put(770,140){\framebox(10,10){}}
\put(770,150){\framebox(10,10){}}
\put(780,150){\framebox(10,10){}}
\put(805,150){\framebox(10,10){}}
\put(815,150){\framebox(10,10){}}
\put(825,150){\framebox(10,10){}}
\put(880,130){\framebox(10,10){}}
\put(880,140){\framebox(10,10){}}
\put(880,150){\framebox(10,10){}}
\put(880,120){\framebox(10,10){}}
\put(910,130){\framebox(10,10){}}
\put(910,140){\framebox(10,10){}}
\put(910,150){\framebox(10,10){}}
\put(920,150){\framebox(10,10){}}
\put(980,140){\framebox(10,10){}}
\put(980,150){\framebox(10,10){}}
\put(990,150){\framebox(10,10){}}
\put(1000,150){\framebox(10,10){}}
\put(945,140){\framebox(10,10){}}
\put(945,150){\framebox(10,10){}}
\put(955,150){\framebox(10,10){}}
\put(955,140){\framebox(10,10){}}
\put(1025,150){\framebox(10,10){}}
\put(1035,150){\framebox(10,10){}}
\put(1045,150){\framebox(10,10){}}
\put(1055,150){\framebox(10,10){}}
\put(120,280){\line( 1,-1){100}}
\put(220,180){\line( 1, 1){100}}
\put(320,280){\line( 1,-1){100}}
\put(420,180){\line( 1, 1){100}}
\put(520,280){\line( 1,-1){100}}
\put(620,180){\line( 1, 1){100}}
\put(720,280){\line( 1,-1){100}}
\put(820,180){\line( 1, 1){100}}
\put(920,280){\line( 1,-1){100}}
\put(420,180){\line( 3, 5){ 60}}
\put(480,280){\line( 1,-1){100}}
\put(580,180){\line( 1, 1){100}}
\put(680,280){\line( 1,-1){100}}
\put(780,180){\line( 1, 1){100}}
\put(880,280){\line( 1,-1){100}}
\put(580,180){\line( 3, 5){ 60}}
\put(640,280){\line( 1,-1){100}}
\put(740,180){\line( 1, 1){100}}
\put(840,280){\line( 1,-1){100}}
\put(620,180){\line( 3, 5){ 60}}
\put(740,180){\line( 3, 5){ 60}}
\put(800,280){\line( 1,-1){100}}
\put(780,180){\line( 3, 5){ 60}}
\put(780,180){\line( 4, 5){ 80}}
\put(860,280){\line( 1,-1){100}}
\put(820,180){\line( 3, 5){ 60}}
\put(415,150){\framebox(10,10){}}
\put( 40,280){\makebox(0,0)[lb]{\raisebox{0pt}[0pt][0pt]{\svtnrm $\Gbarn$:}}}
\put( 40,180){\makebox(0,0)[lb]{\raisebox{0pt}[0pt][0pt]{\svtnrm ${\Gbarn}_+$:}}}
\put(225,150){\makebox(0,0)[b]{\raisebox{0pt}[0pt][0pt]{\svtnrm 0\tiny +}}}
\put(435,150){\makebox(0,0)[b]{\raisebox{0pt}[0pt][0pt]{\tiny +}}}
\put(595,140){\makebox(0,0)[b]{\raisebox{0pt}[0pt][0pt]{\tiny +}}}
\put(640,150){\makebox(0,0)[b]{\raisebox{0pt}[0pt][0pt]{\tiny +}}}
\put(755,130){\makebox(0,0)[b]{\raisebox{0pt}[0pt][0pt]{\tiny +}}}
\put(795,140){\makebox(0,0)[b]{\raisebox{0pt}[0pt][0pt]{\tiny +}}}
\put(845,150){\makebox(0,0)[b]{\raisebox{0pt}[0pt][0pt]{\tiny +}}}
\put(900,120){\makebox(0,0)[b]{\raisebox{0pt}[0pt][0pt]{\tiny +}}}
\put(930,130){\makebox(0,0)[b]{\raisebox{0pt}[0pt][0pt]{\tiny +}}}
\put(974,140){\makebox(0,0)[b]{\raisebox{0pt}[0pt][0pt]{\tiny +}}}
\put(1015,140){\makebox(0,0)[b]{\raisebox{0pt}[0pt][0pt]{\tiny +}}}
\put(1075,150){\makebox(0,0)[b]{\raisebox{0pt}[0pt][0pt]{\tiny +}}}
\end{picture}

\caption{
\label{Young1}
The beginning of the double Young graph in Young diagram labelling.
}

\end{figure}



\Subsection{Bell numbers and the partition algebras\label{bellnumbers}}
The Bell number $b(n)$ is the number of ways to partition a set of $n$
elements. The sequence begins:
$$1,1,2,5,15,\ldots .$$


For a set $S$ denote by $E(S)$ the set of partitions of $S$.
Define $\underline{n} = \{1,2,...,n \}$ and
$\underline{n'} = \{1',2',...,n' \}$.
Consider the sequence of sets $\CC_p(0),\CC_p(0+),\CC_p(1),
\CC_p(1+),\ldots$ defined as follows. Each
$\CC_p(n)$ is the set 
$E(\underline{n}\cup \underline{n'})$ 
of partitions of $\underline{n}\cup \underline{n'}$;
and $\CC_p(n+)$ is the set of partitions of
$\underline{n+1}\cup \underline{(n+1)'}$ in which $n+1$ and $(n+1)'$ always
lie in the same part
(thus $\CC_p(n+) \subset \CC_p(n+1)$;
we will also consider as canonical the inclusion of 
$\CC_p(n)$ in $\CC_p(n+)$ by $p \mapsto p \cup \{\{n+1,n+1'\}\}$).
This sequence can be considered to be an example of a {\em Bell sequence}.
The question, as before, is whether there is an underlying Pascal array.
 
Following~\cite{MartinRollet98}, 
we shall define an array $\{ Y_p(n;\nu) \}_{n;\nu}$
of sets equivalent to $Y_{\YounGG,0}$. 
Here
$\nu$ is of the form $\lambda$ or $\lambda+$ for an integer
partition $\lambda$ of $l$ where $l\leq n$.
Remarkably, it will turn out that
the corresponding Catalan sequence is the Bell sequence above. 


{\de{
We define a {\em half-partition} of $\underline{n}$ to be an ordered
partition of a partition of this set into two parts. The elements of
the second part are known as propagating components.
}}

For $n\in\mathbb{N}$ and an integer partition $\lambda$ of $l\in\mathbb{N}$,
we define
$Y_p(2n;\lambda)$ to be the set of pairs $(\sigma,T)$ where $\sigma$ is a
half-partition of $\underline{n}$ with $l$ propagating components and
$T$ is a standard tableau filling of the Young diagram corresponding to
$\lambda$. We define $Y_p(2n+1;\lambda+)$ to be the set of pairs
$(\sigma,T)$ where $\sigma$ is a half partition of $\underline{n+1}$ with
$l+1\leq n+1$ propagating components including a part which has $n+1$ as a
member and $T$ is a filling of the Young diagram of $\lambda$ as in the
even case.


(1) Bra-ket decomposition. \\
Given a partition of
$\underline{n}\cup \underline{n'}$ in $\CC_p(n)$, we take the set of parts
contained entirely within $\underline{n}$ as the first component of a
half-partition. For the second part we consider the intersections of the
other parts of our original partition with $\underline{n}$ (excluding the
empty set). In this way we obtain (say) $l\leq n$ propagating components.
We totally order these by the lowest number they contain.
Repeating the same procedure with $\underline{n'}$ we obtain a
pair of half-partitions. The parts of the original partition match the
$l$ (ordered) propagating components of the first half-partition with those of
the second half-partition and thus define an element of the symmetric group of
degree $l$. Applying the Robinson-Schensted correspondence we obtain a pair of
tableaux of the same shape $\lambda$ (an integer partition of $l$) and
associate the first to the first half-partition and the second to the second
half-partition. In this way we obtain a pair of elements of $Y_p(2n;\lambda)$.
(See also \cite{Halverson05}.)

Given a partition of $\underline{n+1}\cup \underline{(n+1)'}$ in $\CC_p(n+)$,
we follow the same procedure as above, to obtain a pair of half-partitions
with $l+1$ propagating components in $Y_p(2n+1;\lambda)$, where $\lambda$ is
an integer partition of $l$.

\begin{figure}
\[
\includegraphics{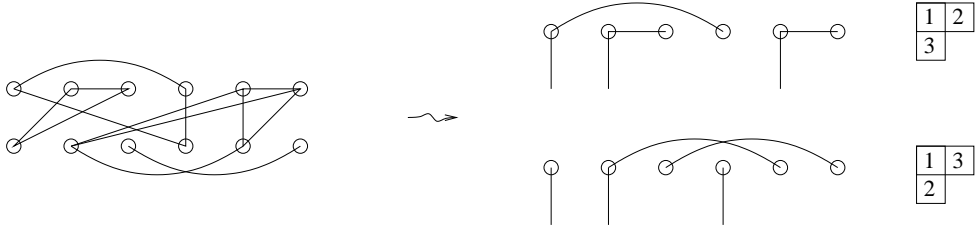}
\]
\caption{\label{f:partitionbraket} Example of bra-ket decomposition for partitions.}
\end{figure}

For an example, take the partition
$$\{\{1,4,4'\},\{2,3,1'\},\{5,6,2',5'\},\{3',6'\}\}$$
of $\underline{6}\cup \underline{6'}$. We obtain the two
half partitions $(\phi,\{\{1,4\},\{2,3\},\{5,6\}\})$ and
$(\{\{3,6\}\},\{\{1\},\{2,5\},\{4\}\})$, each with $3$ propagating components.
The permutation in $S_3$ is given by $1\to 3$, $2\to 1$ and $3\to 2$.
The half partitions, together with their corresponding tableaux
(from the Robinson-Schensted correspondence applied to this permutation)
are shown in Figure~\ref{f:partitionbraket}.


(2) Edge maps. \\
If $\lambda$ is an integer partition of $l\leq n$, we define
$\edge{(2n;\lambda),(2n+1,\lambda+)}$ to be the map that takes a pair
$(\sigma,T)$ to the half partition $\sigma$ with an extra propagating part
$\{n+1\}$ and the same tableau $T$.

Let $\nu$ be an integer partition such that $\lambda-\nu$ is of the form
$(0,0,\ldots ,1,0,0,\ldots )$.
Before defining $\edge{(2n;\lambda),(2n+1,\nu+)}$ we need some notation.
Let $V_{\nu}$ denote the corresponding
irreducible $S_{l-1}$-module. Then ${}_{S_{l}}Ind_{S_{l-1}}V_{\nu}=
\mathbb{C}S_{l}\otimes_{\mathbb{C}S_{l-1}} V_{\nu}$ has a composition
series of irreducible $S_{l-1}$-modules $V_{\lambda'}$ where $\lambda'$
varies over the integer partitions such that $\lambda'-\nu$ is of the
form $(0,0,\ldots ,1,0,0,\ldots )$ (thus including $\lambda$). It has a
basis consisting of the elements $\pi_k\otimes U$ where $U$ is
(a basis element corresponding to) a standard tableau filling of the Young
diagram corresponding to $\nu$ and
$$\pi_k:=\left(\begin{array}{cccccccc} 1 & 2 & \cdots & l-1-k & l-k & l+1-k &
\cdots & l \\ 1 & 2 & \cdots & l-1-k & l & l-k & \cdots & l-1 \end{array}
\right),$$
for $k=1,2,\ldots ,l$; note that $\pi_1,\pi_2,\ldots ,\pi_l$ is a
system of coset representatives for $S_{l-1}$ as a subgroup of $S_l$.

We choose a bijection between this basis for ${}_{S_{l}}Ind_{S_{l-1}}V_{\nu}$
and a basis compatible with a composition series for this module as
$S_l$-module. Such a bijection, on restriction, gives rise to
a bijection between a subset of the set of elements $\pi_k\otimes U$ as
above and standard tableau fillings of a Young diagram of
shape $\lambda$. This defines a map $\edge{(2n;\lambda),(2n+1,\nu+)}$
in the following way. Given a pair $(\sigma,T)$ in $\CC_p(n)$, let
$\pi_k\otimes U$ be the corresponding element under this bijection.
Then $\edge{(2n;\lambda),(2n+1,\nu+)}$ takes $(\sigma,T)$ to
to the half partition $\sigma$ transformed into
a half partition $\sigma'$ of $\underline{n+1}$ by adding $n+1$ to the
$k$th propagating part of $\sigma$, together with the tableau $U$.

For $\lambda$ an integer partition of $l$, we define
$\edge{(2n+1;\lambda+),(2n+2,\lambda)}$ to be the map which transforms the
propagating component containing $n+1$ into a
non-propagating component and leaves $T$ unchanged.
Let $\mu$ be an integer partition such that
$\mu-\lambda$ is of the form $(0,0,\ldots ,1,0,0,\ldots )$. We define
$\edge{(2n+1;\lambda+),(2n+2,\mu)}$ to be the map which leaves $\sigma$
unchanged and labels the new box in $\mu$ with $l+1$.


We have:

\prl(p:partitionalgebraarray)
The sets $Y_p(n;\lambda)$ described above form a
$(\YounGG,0)$-array of sets with corresponding Catalan sequence of
of sets given by $\CC_p(0)\subseteq \CC_p(0+)\subseteq \CC_p(1)
\subseteq \CC_p(1+)\subseteq \cdots$. \Qed
\end{pr}

\begin{figure}
\[
\includegraphics{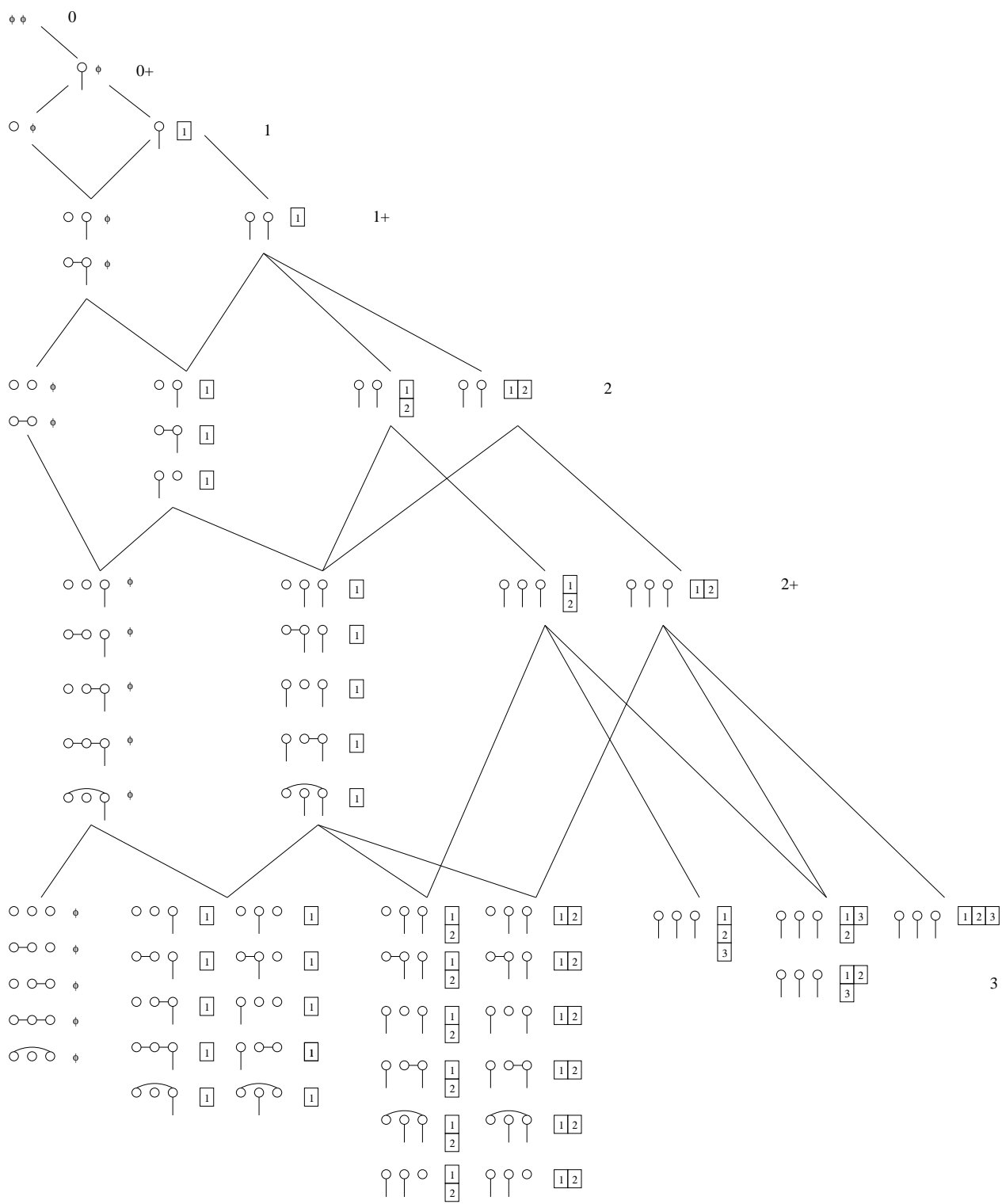}
\]
\caption{\label{f:partitionarray} The start of the array of half-partitions.}
\end{figure}

\begin{figure}

\setlength{\unitlength}{0.009in}%
\begin{picture}(567,741)(-30,88)
\thicklines
\put( 20,820){\line( 3,-2){ 60}}
\put( 80,780){\line(-3,-2){ 60}}
\put( 20,740){\line( 3,-2){ 60}}
\put( 80,700){\line(-3,-2){ 60}}
\put( 20,660){\line( 3,-2){ 60}}
\put( 80,620){\line(-3,-2){ 60}}
\put( 20,580){\line( 3,-2){ 60}}
\put( 80,540){\line(-3,-2){ 60}}
\put( 20,500){\line( 3,-2){ 60}}
\put( 80,460){\line(-3,-2){ 60}}
\put( 20,420){\line( 3,-2){ 60}}
\put( 80,380){\line(-3,-2){ 60}}
\put( 20,340){\line( 3,-2){ 60}}
\put( 80,300){\line(-3,-2){ 60}}
\put( 20,260){\line( 3,-2){ 60}}
\put( 80,220){\line(-3,-2){ 60}}
\put( 20,180){\line( 3,-2){ 60}}
\put( 80,140){\line(-3,-2){ 60}}
\put( 80,780){\line( 3,-2){ 60}}
\put(140,740){\line(-3,-2){ 60}}
\put( 80,700){\line( 3,-2){ 60}}
\put(140,660){\line(-3,-2){ 60}}
\put( 80,700){\line( 3,-2){ 60}}
\put(140,660){\line(-3,-2){ 60}}
\put( 80,620){\line( 3,-2){ 60}}
\put(140,580){\line(-3,-2){ 60}}
\put( 80,540){\line( 3,-2){ 60}}
\put(140,500){\line(-3,-2){ 60}}
\put( 80,460){\line( 3,-2){ 60}}
\put(140,420){\line(-3,-2){ 60}}
\put( 80,380){\line( 3,-2){ 60}}
\put(140,340){\line(-3,-2){ 60}}
\put( 80,300){\line( 3,-2){ 60}}
\put(140,260){\line(-3,-2){ 60}}
\put( 80,220){\line( 3,-2){ 60}}
\put(140,180){\line(-3,-2){ 60}}
\put(140,740){\line( 3,-2){ 60}}
\put(200,700){\line(-3,-2){ 60}}
\put(140,660){\line( 3,-2){ 60}}
\put(200,620){\line(-3,-2){ 60}}
\put(140,580){\line( 3,-2){ 60}}
\put(200,540){\line(-3,-2){ 60}}
\put(140,500){\line( 3,-2){ 60}}
\put(200,460){\line(-3,-2){ 60}}
\put(140,420){\line( 3,-2){ 60}}
\put(200,380){\line(-3,-2){ 60}}
\put(140,340){\line( 3,-2){ 60}}
\put(200,300){\line(-3,-2){ 60}}
\put(140,260){\line( 3,-2){ 60}}
\put(200,220){\line(-3,-2){ 60}}
\put(200,700){\line( 3,-2){ 60}}
\put(260,660){\line(-3,-2){ 60}}
\put(200,620){\line( 3,-2){ 60}}
\put(260,580){\line(-3,-2){ 60}}
\put(200,540){\line( 3,-2){ 60}}
\put(260,500){\line(-3,-2){ 60}}
\put(200,460){\line( 3,-2){ 60}}
\put(260,420){\line(-3,-2){ 60}}
\put(200,380){\line( 3,-2){ 60}}
\put(260,340){\line(-3,-2){ 60}}
\put(200,300){\line( 3,-2){ 60}}
\put(260,260){\line(-3,-2){ 60}}
\put(260,660){\line( 3,-2){ 60}}
\put(320,620){\line(-3,-2){ 60}}
\put(260,580){\line( 3,-2){ 60}}
\put(320,540){\line(-3,-2){ 60}}
\put(260,500){\line( 3,-2){ 60}}
\put(320,460){\line(-3,-2){ 60}}
\put(260,420){\line( 3,-2){ 60}}
\put(320,380){\line(-3,-2){ 60}}
\put(260,340){\line( 3,-2){ 60}}
\put(320,300){\line(-3,-2){ 60}}
\put(320,540){\line( 3,-2){ 60}}
\put(380,500){\line(-3,-2){ 60}}
\put(320,460){\line( 3,-2){ 60}}
\put(380,420){\line(-3,-2){ 60}}
\put(320,380){\line( 3,-2){ 60}}
\put(380,340){\line(-3,-2){ 60}}
\put(380,580){\line( 3,-2){ 60}}
\put(440,540){\line(-3,-2){ 60}}
\put(380,500){\line( 3,-2){ 60}}
\put(440,460){\line(-3,-2){ 60}}
\put(380,420){\line( 3,-2){ 60}}
\put(440,380){\line(-3,-2){ 60}}
\put(440,540){\line( 3,-2){ 60}}
\put(500,500){\line(-3,-2){ 60}}
\put(440,460){\line( 3,-2){ 60}}
\put(500,420){\line(-3,-2){ 60}}
\put(500,500){\line( 3,-2){ 60}}
\put(560,460){\line(-3,-2){ 60}}
\put(200,700){\line( 1,-1){ 40}}
\put(240,660){\line(-1,-1){ 40}}
\put(200,620){\line( 1,-1){ 40}}
\put(240,580){\line(-1,-1){ 40}}
\put(240,660){\line( 3,-2){ 60}}
\put(300,620){\line(-3,-2){ 60}}
\put(300,620){\line( 3,-2){ 60}}
\put(360,580){\line(-3,-2){ 60}}
\put(320,620){\line( 1,-1){ 40}}
\put(360,580){\line(-1,-1){ 40}}
\put(300,620){\line( 1,-1){ 40}}
\put(340,580){\line(-1,-1){ 40}}
\put(240,580){\line( 3,-2){ 60}}
\put(300,540){\line(-3,-2){ 60}}
\put(200,540){\line( 1,-1){ 40}}
\put(240,500){\line(-1,-1){ 40}}
\put(200,460){\line( 1,-1){ 40}}
\put(240,420){\line(-1,-1){ 40}}
\put(240,500){\line( 3,-2){ 60}}
\put(300,460){\line(-3,-2){ 60}}
\put(300,540){\line( 3,-2){ 60}}
\put(360,500){\line(-3,-2){ 60}}
\put(340,580){\line( 3,-2){ 60}}
\put(400,540){\line(-3,-2){ 60}}
\put(300,540){\line( 1,-1){ 40}}
\put(340,500){\line(-1,-1){ 40}}
\put(320,540){\line( 1,-1){ 40}}
\put(360,500){\line(-1,-1){ 40}}
\put(200,380){\line( 1,-1){ 40}}
\put(240,340){\line(-1,-1){ 40}}
\put(200,300){\line( 1,-1){ 40}}
\put(240,260){\line(-1,-1){ 40}}
\put(240,420){\line( 3,-2){ 60}}
\put(300,380){\line(-3,-2){ 60}}
\put(240,340){\line( 3,-2){ 60}}
\put(300,300){\line(-3,-2){ 60}}
\put(300,460){\line( 3,-2){ 60}}
\put(360,420){\line(-3,-2){ 60}}
\put(300,380){\line( 3,-2){ 60}}
\put(360,340){\line(-3,-2){ 60}}
\put(360,500){\line( 3,-2){ 60}}
\put(420,460){\line(-3,-2){ 60}}
\put(360,420){\line( 3,-2){ 60}}
\put(420,380){\line(-3,-2){ 60}}
\put(420,540){\line( 3,-2){ 60}}
\put(480,500){\line(-3,-2){ 60}}
\put(420,460){\line( 3,-2){ 60}}
\put(480,420){\line(-3,-2){ 60}}
\put(480,500){\line( 3,-2){ 60}}
\put(540,460){\line(-3,-2){ 60}}
\put(300,460){\line( 1,-1){ 40}}
\put(340,420){\line(-1,-1){ 40}}
\put(300,380){\line( 1,-1){ 40}}
\put(340,340){\line(-1,-1){ 40}}
\put(320,460){\line( 1,-1){ 40}}
\put(360,420){\line(-1,-1){ 40}}
\put(320,380){\line( 1,-1){ 40}}
\put(360,340){\line(-1,-1){ 40}}
\put(440,540){\line( 1,-1){ 40}}
\put(480,500){\line(-1,-1){ 40}}
\put(440,460){\line( 1,-1){ 40}}
\put(480,420){\line(-1,-1){ 40}}
\put(340,500){\line( 3,-2){ 60}}
\put(400,460){\line(-3,-2){ 60}}
\put(340,420){\line( 3,-2){ 60}}
\put(400,380){\line(-3,-2){ 60}}
\put(400,540){\line( 3,-2){ 60}}
\put(460,500){\line(-3,-2){ 60}}
\put(420,540){\line( 1,-1){ 40}}
\put(460,500){\line(-1,-1){ 40}}
\put(400,540){\line( 1,-1){ 40}}
\put(440,500){\line(-1,-1){ 40}}
\put(420,460){\line( 1,-1){ 40}}
\put(460,420){\line(-1,-1){ 40}}
\put(400,460){\line( 1,-1){ 40}}
\put(440,420){\line(-1,-1){ 40}}
\put(420,540){\line( 5,-4){ 50}}
\put(470,500){\line(-5,-4){ 50}}
\put(470,500){\line( 3,-2){ 60}}
\put(530,460){\line(-3,-2){ 60}}
\put(360,580){\line( 3,-2){ 60}}
\put(420,540){\line(-3,-2){ 60}}
\put(321,619){\line( 3,-2){ 60}}
\put(381,579){\line(-3,-2){ 60}}
\put(440,500){\line( 3,-2){ 60}}
\put(500,460){\line(-3,-2){ 60}}
\put(400,460){\line( 3,-2){ 60}}
\put(460,420){\line(-3,-2){ 60}}
\put(460,500){\line( 3,-2){ 60}}
\put(520,460){\line(-3,-2){ 60}}
\put(420,460){\line( 5,-4){ 50}}
\put(470,420){\line(-5,-4){ 50}}
\put( 10,820){\makebox(0,0)[lb]{\raisebox{0pt}[0pt][0pt]{\twlrm 1}}}
\put( 10,745){\makebox(0,0)[lb]{\raisebox{0pt}[0pt][0pt]{\tiny 1}}}
\put(135,750){\makebox(0,0)[lb]{\raisebox{0pt}[0pt][0pt]{\tiny 1}}}
\put(195,710){\makebox(0,0)[lb]{\raisebox{0pt}[0pt][0pt]{\tiny 1}}}
\put(320,630){\makebox(0,0)[lb]{\raisebox{0pt}[0pt][0pt]{\tiny 1}}}
\put(280,615){\makebox(0,0)[lb]{\raisebox{0pt}[0pt][0pt]{\tiny 1}}}
\put(380,585){\makebox(0,0)[lb]{\raisebox{0pt}[0pt][0pt]{\tiny 1}}}
\put(325,575){\makebox(0,0)[lb]{\raisebox{0pt}[0pt][0pt]{\tiny 1}}}
\put(380,535){\makebox(0,0)[lb]{\raisebox{0pt}[0pt][0pt]{\tiny 1}}}
\put(570,460){\makebox(0,0)[lb]{\raisebox{0pt}[0pt][0pt]{\tiny 1}}}
\put(485,455){\makebox(0,0)[lb]{\raisebox{0pt}[0pt][0pt]{\tiny 1}}}
\put(435,550){\makebox(0,0)[lb]{\raisebox{0pt}[0pt][0pt]{\tiny 1}}}
\put( 75,790){\makebox(0,0)[lb]{\raisebox{0pt}[0pt][0pt]{\tiny 1}}}
\put(255,670){\makebox(0,0)[lb]{\raisebox{0pt}[0pt][0pt]{\tiny 1}}}
\put( 75,710){\makebox(0,0)[lb]{\raisebox{0pt}[0pt][0pt]{\tiny 2}}}
\put( 10,665){\makebox(0,0)[lb]{\raisebox{0pt}[0pt][0pt]{\tiny 2}}}
\put(135,670){\makebox(0,0)[lb]{\raisebox{0pt}[0pt][0pt]{\tiny 3}}}
\put(220,660){\makebox(0,0)[lb]{\raisebox{0pt}[0pt][0pt]{\tiny 1}}}
\put( 75,630){\makebox(0,0)[lb]{\raisebox{0pt}[0pt][0pt]{\tiny 5}}}
\put(135,590){\makebox(0,0)[lb]{\raisebox{0pt}[0pt][0pt]{\tiny 10}}}
\put( 10,585){\makebox(0,0)[lb]{\raisebox{0pt}[0pt][0pt]{\tiny 5}}}
\put(195,630){\makebox(0,0)[lb]{\raisebox{0pt}[0pt][0pt]{\tiny 5}}}
\put( 75,550){\makebox(0,0)[lb]{\raisebox{0pt}[0pt][0pt]{\tiny 15}}}
\put(220,580){\makebox(0,0)[lb]{\raisebox{0pt}[0pt][0pt]{\tiny 6}}}
\put(275,580){\makebox(0,0)[lb]{\raisebox{0pt}[0pt][0pt]{\tiny 6}}}
\put(195,550){\makebox(0,0)[lb]{\raisebox{0pt}[0pt][0pt]{\tiny 22}}}
\put(135,510){\makebox(0,0)[lb]{\raisebox{0pt}[0pt][0pt]{\tiny 37}}}
\put(280,540){\makebox(0,0)[lb]{\raisebox{0pt}[0pt][0pt]{\tiny 9}}}
\put(220,500){\makebox(0,0)[lb]{\raisebox{0pt}[0pt][0pt]{\tiny 31}}}
\put(275,500){\makebox(0,0)[lb]{\raisebox{0pt}[0pt][0pt]{\tiny 31}}}
\put(195,470){\makebox(0,0)[lb]{\raisebox{0pt}[0pt][0pt]{\tiny 99}}}
\put( 75,470){\makebox(0,0)[lb]{\raisebox{0pt}[0pt][0pt]{\tiny 52}}}
\put(135,430){\makebox(0,0)[lb]{\raisebox{0pt}[0pt][0pt]{\tiny 151}}}
\put( 75,390){\makebox(0,0)[lb]{\raisebox{0pt}[0pt][0pt]{\tiny 203}}}
\put(  5,505){\makebox(0,0)[lb]{\raisebox{0pt}[0pt][0pt]{\tiny 15}}}
\put(  5,425){\makebox(0,0)[lb]{\raisebox{0pt}[0pt][0pt]{\tiny 52}}}
\put(  5,345){\makebox(0,0)[lb]{\raisebox{0pt}[0pt][0pt]{\tiny 203}}}
\put(365,580){\makebox(0,0)[lb]{\raisebox{0pt}[0pt][0pt]{\tiny 2}}}
\put(335,540){\makebox(0,0)[lb]{\raisebox{0pt}[0pt][0pt]{\tiny 9}}}
\put(495,510){\makebox(0,0)[lb]{\raisebox{0pt}[0pt][0pt]{\tiny 1}}}
\put(425,500){\makebox(0,0)[lb]{\raisebox{0pt}[0pt][0pt]{\tiny 1}}}
\put(425,540){\makebox(0,0)[lb]{\raisebox{0pt}[0pt][0pt]{\tiny 2}}}
\put(325,500){\makebox(0,0)[lb]{\raisebox{0pt}[0pt][0pt]{\tiny 10}}}
\put(355,505){\makebox(0,0)[lb]{\raisebox{0pt}[0pt][0pt]{\tiny 20}}}
\put(390,500){\makebox(0,0)[lb]{\raisebox{0pt}[0pt][0pt]{\tiny 10}}}
\put(280,460){\makebox(0,0)[lb]{\raisebox{0pt}[0pt][0pt]{\tiny 61}}}
\put(335,460){\makebox(0,0)[lb]{\raisebox{0pt}[0pt][0pt]{\tiny 61}}}
\put(445,500){\makebox(0,0)[lb]{\raisebox{0pt}[0pt][0pt]{\tiny 3}}}
\put(485,500){\makebox(0,0)[lb]{\raisebox{0pt}[0pt][0pt]{\tiny 3}}}
\put(460,500){\makebox(0,0)[lb]{\raisebox{0pt}[0pt][0pt]{\tiny 2}}}
\put(385,460){\makebox(0,0)[lb]{\raisebox{0pt}[0pt][0pt]{\tiny 14}}}
\put(420,465){\makebox(0,0)[lb]{\raisebox{0pt}[0pt][0pt]{\tiny 28}}}
\put(455,460){\makebox(0,0)[lb]{\raisebox{0pt}[0pt][0pt]{\tiny 14}}}
\put(505,460){\makebox(0,0)[lb]{\raisebox{0pt}[0pt][0pt]{\tiny 3}}}
\put(545,460){\makebox(0,0)[lb]{\raisebox{0pt}[0pt][0pt]{\tiny 3}}}
\put(525,460){\makebox(0,0)[lb]{\raisebox{0pt}[0pt][0pt]{\tiny 2}}}
\put(215,420){\makebox(0,0)[lb]{\raisebox{0pt}[0pt][0pt]{\tiny 160}}}
\put(280,420){\makebox(0,0)[lb]{\raisebox{0pt}[0pt][0pt]{\tiny 160}}}
\put(320,420){\makebox(0,0)[lb]{\raisebox{0pt}[0pt][0pt]{\tiny 75}}}
\put(355,425){\makebox(0,0)[lb]{\raisebox{0pt}[0pt][0pt]{\tiny 150}}}
\put(275,380){\makebox(0,0)[lb]{\raisebox{0pt}[0pt][0pt]{\tiny 385}}}
\put(335,380){\makebox(0,0)[lb]{\raisebox{0pt}[0pt][0pt]{\tiny 385}}}
\put(195,390){\makebox(0,0)[lb]{\raisebox{0pt}[0pt][0pt]{\tiny 471}}}
\put(215,340){\makebox(0,0)[lb]{\raisebox{0pt}[0pt][0pt]{\tiny 856}}}
\put(275,340){\makebox(0,0)[lb]{\raisebox{0pt}[0pt][0pt]{\tiny 856}}}
\put(135,350){\makebox(0,0)[lb]{\raisebox{0pt}[0pt][0pt]{\tiny 674}}}
\put(190,315){\makebox(0,0)[lb]{\raisebox{0pt}[0pt][0pt]{\tiny 2386}}}
\put( 75,310){\makebox(0,0)[lb]{\raisebox{0pt}[0pt][0pt]{\tiny 877}}}
\put(130,275){\makebox(0,0)[lb]{\raisebox{0pt}[0pt][0pt]{\tiny 3263}}}
\put(  5,270){\makebox(0,0)[lb]{\raisebox{0pt}[0pt][0pt]{\tiny 877}}}
\put( 70,235){\makebox(0,0)[lb]{\raisebox{0pt}[0pt][0pt]{\tiny 4140}}}
\put(  5,190){\makebox(0,0)[lb]{\raisebox{0pt}[0pt][0pt]{\tiny 4140}}}
\put(485,415){\makebox(0,0)[lb]{\raisebox{0pt}[0pt][0pt]{\tiny 45}}}
\put(460,415){\makebox(0,0)[lb]{\raisebox{0pt}[0pt][0pt]{\tiny 30}}}
\put(505,410){\makebox(0,0)[lb]{\raisebox{0pt}[0pt][0pt]{\tiny 15}}}
\put(445,370){\makebox(0,0)[lb]{\raisebox{0pt}[0pt][0pt]{\tiny 135}}}
\put(390,420){\makebox(0,0)[lb]{\raisebox{0pt}[0pt][0pt]{\tiny 75}}}
\put(425,420){\makebox(0,0)[lb]{\raisebox{0pt}[0pt][0pt]{\tiny 15}}}
\put(446,417){\makebox(0,0)[lb]{\raisebox{0pt}[0pt][0pt]{\tiny 45}}}
\put(383,375){\makebox(0,0)[lb]{\raisebox{0pt}[0pt][0pt]{\tiny 135}}}
\put(414,385){\makebox(0,0)[lb]{\raisebox{0pt}[0pt][0pt]{\tiny 270}}}
\put(321,340){\makebox(0,0)[lb]{\raisebox{0pt}[0pt][0pt]{\tiny 520}}}
\put(393,337){\makebox(0,0)[lb]{\raisebox{0pt}[0pt][0pt]{\tiny 520}}}
\put(356,346){\makebox(0,0)[lb]{\raisebox{0pt}[0pt][0pt]{\tiny 1040}}}
\put(332,294){\makebox(0,0)[lb]{\raisebox{0pt}[0pt][0pt]{\tiny 2416}}}
\put(272,300){\makebox(0,0)[lb]{\raisebox{0pt}[0pt][0pt]{\tiny 2416}}}
\put(273,256){\makebox(0,0)[lb]{\raisebox{0pt}[0pt][0pt]{\tiny 4802}}}
\put(211,256){\makebox(0,0)[lb]{\raisebox{0pt}[0pt][0pt]{\tiny 4802}}}
\put(209,212){\makebox(0,0)[lb]{\raisebox{0pt}[0pt][0pt]{\tiny 12867}}}
\put(145,170){\makebox(0,0)[lb]{\raisebox{0pt}[0pt][0pt]{\tiny 17007}}}
\put( 89,131){\makebox(0,0)[lb]{\raisebox{0pt}[0pt][0pt]{\tiny 21147}}}
\put(  3, 88){\makebox(0,0)[lb]{\raisebox{0pt}[0pt][0pt]{\tiny 21147}}}
\end{picture}

\caption{
\label{Young2}
Counting walks on the double Young graph.
That is, 
part of the `Bratteli' diagram
for $P_0 \subset P_{0+} \subset P_{1} \subset P_{1+} \subset P_{2} \cdots$
(to be precise, all restrictions
of standard module $\SS({(0)}) (9)$),
in which each  $\SS(\lam)$  is represented by its dimension. 
}

\end{figure}
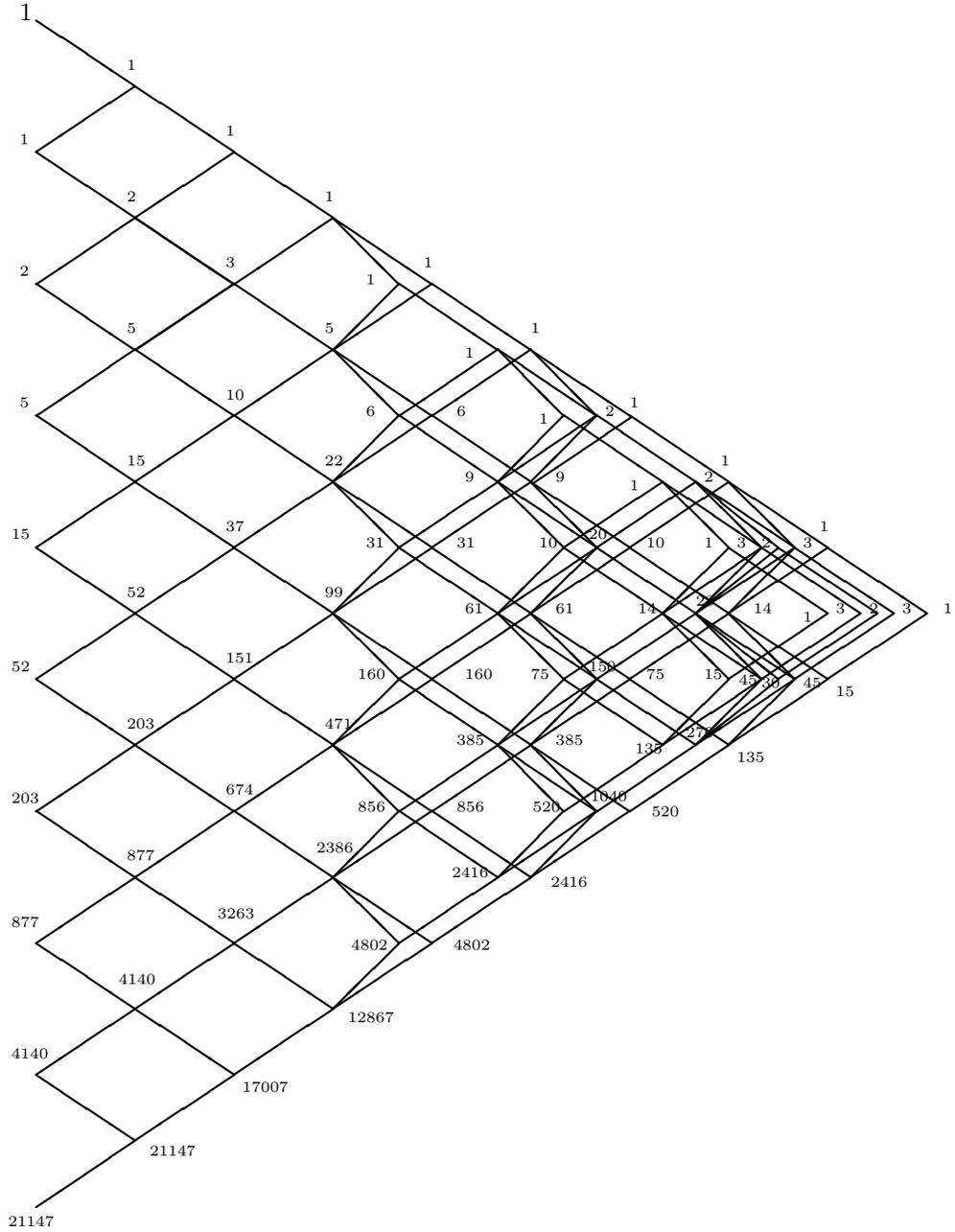

The start of the array is displayed in Figure~\ref{f:partitionarray}.
Note that propagating lines (indicating propagating components) have
always been drawn at the right hand end of the part. The
sizes of some of the sets in the array are displayed in Figure~\ref{Young2}.


As a consequence of the above and Proposition~\ref{p:equivalence} we have:


\prl(dYg){\rm \cite{MartinRollet98}}
Let $v_0$ denote the empty integer partition, a vertex of the double Young
graph, $\YounGG$. Then the sequence $\{ N_{\YounGG}(2n;v_0) \}_n$ is the
Bell numbers. 
\end{pr}


Just as for $(A_{\infty},0)$, 
the above array also arises in an algebraic context,
which we now describe.

The \emph{partition algebra} $P_n$ has
$\mathbb{C}$-basis $\CC_p(n)$. Let $\delta\in \mathbb{C}$.
The multiplication is defined as follows (see \cite{Martin94}).
For
$T\subseteq S$ sets, 
and $p\in E(S)$, let $p|_{T}$ denote the \emph{restriction} of
$p$ to $T$, obtained by removing all the elements of $S\setminus T$ from each
part of $p$.
For $x$ a relation, let $\overline{x}$ denote its transitive closure.

Let $\underline{n''}:=\{1'',2'',\ldots ,n''\}$.
For a partition $p$ of $\underline{n}\cup \underline{n'}$, let $p'$
denote the partition of $\underline{n'}\cup\underline{n''}$ obtained by
replacing each $i\in\underline{n}$ with $i'$ and each $i'\in\underline{n'}$
with $i''$ in every part of $p$.
For a partition $p$ of $\underline{n}\cup \underline{n''}$, let $p'$
denote the partition of $\underline{n}\cup\underline{n'}$ obtained by
replacing each $i''\in\underline{n''}$ with $i'$ 
and leaving each $i\in\underline{n}$ alone. 

\begin{figure}
\[
\includegraphics{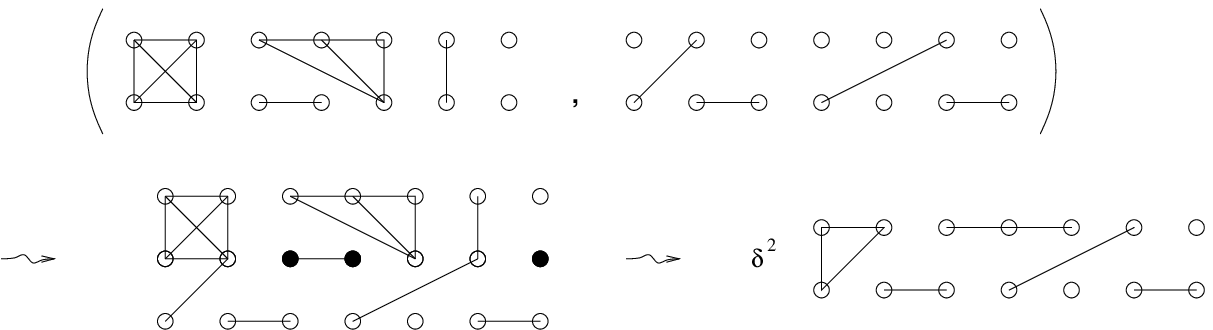}
\]
\caption{\label{partitionexample} Multiplication in the partition algebra.}
\end{figure}

Let $p,q$ be partitions in $\CC_p(n)$, regarded as equivalence
relations. The product of $p$ and $q$ in the
partition algebra is then defined to be $\delta^k$ times
$(\overline{p\cup q'}|_{\underline{n}\cup \underline{n''}})'$, 
where $k$ is the number of
parts of the transitive closure 
$\overline{p\cup q'}$
contained entirely in $\underline{n'}$.
For an example, see Figure~\ref{partitionexample}. The parts contributing
to the exponent of $\delta$ in the product are indicated by filled-in circles.

There is also a related algebra $P_n^+$, defined as the subalgebra
of $P_{n+1}$ with basis $\CC_p(n+)\subseteq \CC_p(n+1)$;
see \cite{Martin2000}. Following~\cite{Martin2000}, we have the following.

\prl(partitiontower)
The tower of algebras $P_0\subseteq P_0^+\subseteq P_1
\subseteq P_1^+\subseteq \cdots $ is a $(\YounGG ,0)$-tower of algebras
(recall that $\YounGG$ is the double Young graph defined above).
\end{pr}
{\bf Proof:}
We explain why the axioms (N1) to (N9) hold. Axiom (N1) is clear from the
definition. For $l=0,1,2,\ldots ,n$, consider the ideal $I_l$ of
$P_n$ generated (and indeed spanned) by partitions with at most $l$ propagating
components (i.e. parts containing both unprimed and primed numbers).

For each $l$, the quotient $I_l/I_{l-1}$ is a left $P_n$-module with basis
given by partitions with exactly $l$ propagating components.
By the definition of multiplication in $P_n$, it is clear that the
left $P_n$ action on such basis partitions does not depend on the parts
of the partition entirely contained in $\underline{n'}$, and cannot change
such parts. It follows that $I_l/I_{l-1}$ is the direct sum of isomorphic
copies of a module $M_n(l)$ indexed by the half-partitions of
$\underline{n'}$ with $l$ propagating components. By fixing such a half
partition $\sigma$ we obtain a copy of $M_n(l)$ which has basis given by
partitions of $\underline{n}\cup \underline{n'}$ whose ket component
contains the half-partition $\sigma$. Via the description of the bra-ket
decomposition we can see that these correspond to pairs consisting of a
half-partition of $\underline{n}$ together with an element of the symmetric
group $S_l$.

We embed $S_l$ into $S_n$ as permutations of the smallest elements of each
of the $l$ propagating components of $\sigma$. Since $S_n$ is naturally
embedded in $P_n$ (with each permutation $\rho$ regarded as a pair partition
with parts of the form $\{i,\rho(i)'\}$ for $i\in\underline{n}$),
we obtain an embedding of $S_l$ into $P_n$. With this embedding, the
right action of $P_n$ on $I_l/I_{l-1}$ restricts to a right action of
$S_l$ on $M_n(l)$. On the level of the basis pairs consisting of a half
partition and an element of $S_l$, described above, this is just the regular
action of $S_l$ on the second element of the pair, so $M_n(l)$ becomes a
free right $S_l$-module.

Let $F$ be the exact functor from left $S_l$-modules to left $P_n$-modules
taking an $S_l$-module $V$ to $M_n(l)\otimes_{S_l} V$. 
For $\lambda$ an integer partition of $l$, we define
the module $\Delta_n(\lambda)$ to be
$F(V_{\lambda})=M_n(l)\otimes_{S_l} V_{\lambda}$,
where, as usual, $V_{\lambda}$ denotes the irreducible $S_l$-module indexed by
$\lambda$. Then it is clear that a basis for $\Delta_n(\lambda)$ is given by
the elements of $Y_p(2n,\lambda)$.

For $P_n^+$, we consider, for each $l$, the left ideal $I_l$ generated by
partitions with at most $l+1$ propagating components.
The quotient $I_l/I_{l-1}$ is a
$P_n^+$-module which can be written as a direct sum of copies of a module
$M_{n+}(l+)$ indexed by half-partitions of $\underline{n'}$ with $l+1$
propagating components including one containing $(n+1)'$.
Each module $M_{n+}(l+)$ has basis given by pairs
consisting of a half partition of $\underline{n+1}$ with $l+1$ propagating
components including one which contains $n+1$, and an element of $S_l$.
It is naturally a right $S_l$-module (where $S_l$ acts regularly on the
second element of the pair). For $\lambda$ an integer partition of $l$,
we define the module $\Delta_n(\lambda)$ to be the tensor product
$M_{n+}(l+)\otimes_{S_l} V_{\lambda}$. A basis for $\Delta_n(\lambda)$ is
given by the elements of $Y_p(2n+1,\lambda+)$.

Let $e_n$ be the idempotent given by $1/\delta$ times the partition
$\{\{1,1'\},\{2,2'\},
\ldots ,\{n-1,(n-1)'\},\{n\},\{n'\}\}$, and let $F_n$ be the functor taking
a module $M$ to $e_nM$ as in~\cite{Martin2000}. Acting $e_n$ on the basis
elements for $\Delta_n(\lambda)$ gives all basis elements containing
a part $\{n\}$ (noting that $e_n$ kills any basis element which has a
propagating component whose intersection with $\underline{n}$ is $\{n\}$).
These span a subspace isomorphic to $\Delta_{n-1}(\lambda)$ over $P_{n-1}$.
Hence $F_n\Delta_n(\lambda)=\Delta_{n-1}(\lambda)$. If we define $e_{n+}$
to be the idempotent given by $1/\delta$ times the partition
$\{\{1,1'\},\{2,2'\},\ldots ,\{n-1,(n-1)'\},\{n\},\{n'\},\{n+1,(n+1)'\}\}$
and $F_{n+}M=e_{n+}M$ we obtain similarly that
$F_{n+}\Delta_{n+}(\lambda+)=\Delta_{(n-1)+}(\lambda+)$ if
$l\leq n-1$ and is zero otherwise. It is clear that $F_n$ and $F_{n+}$ are
both exact. In this way, (N3), (N6) and (N7) follow with $N=2$;
the labelling is designed for this to work.

It follows easily by considering the action of elements of $P_n$ with
various numbers of propagating components, that modules $\Delta_n(\lambda)$
with $\lambda$ having different numbers of propagating components cannot
be isomorphic.
Under the natural embedding of the symmetric group $S_n$ into
$P_n$ the module $\Delta_n(\lambda)$ becomes the $S_n$-module $V_{\lambda}$
if $\lambda$ is a partition of $n$. It follows that no two modules
$\Delta_n(\lambda),\Delta_n(\mu)$ such that $\lambda\not=\mu$ and
$\lambda$ and $\mu$ each have $n$ propagating components can be
isomorphic. Suppose that $\lambda$ and $\mu$ are distinct and have $k<n$
propagating components and $\Delta_n(\lambda)\cong \Delta_n(\mu)$. Then
by applying $F_kF_{k+1}\cdots F_n$, we obtain
$\Delta_k(\lambda)\cong \Delta_k(\mu)$, a contradiction to the above.
A similar argument applies to the $P_n^+$-modules $\Delta_{n+}(\lambda+)$.
Thus axiom (N2) holds.

To see axiom (N4), we first consider restricting from $P_n$ to $P_{(n-1)+}$.
It is clear that the basis elements of $\Delta_n(\lambda)$ with a propagating
component containing $\{n\}$ form a submodule isomorphic over $P_{(n-1)+}$ to
$M_{(n-1)+}(l-1)\otimes {}_{S_{l-1}}V_{\lambda}$ and thus to a module with a
series whose sections are isomorphic to the modules $\Delta_{n+}(\nu)$ where
$\lambda-\nu$ is of the form $(0,0,\ldots ,1,0,0,\ldots )$ (where $\lambda$
is a partition of $l$).
The remaining basis elements of $\Delta_n(\lambda)$
correspond to a basis for the quotient which
is isomorphic to $\Delta_{(n-1)+}(\lambda)$ over $P_{(n-1)+}$.

Next, restricting from $P_n^+$ to $P_n$, we see that the basis elements of
$\Delta_{n+}(\lambda+)$ with a propagating component of the form $\{n+1\}$
form a submodule (over $P_n$) isomorphic to $\Delta_n(\lambda)$.
The remaining basis elements of $\Delta_{n+}(\lambda+)$
form a basis for the quotient, isomorphic
to $M_{n}(l+)\otimes \mathbb{C}S_{l+1}\otimes_{\mathbb{C}S_l} V_{\lambda}$. Such a basis element $(\sigma,T)$ has a non-trivial part containing
$\{n+1\}$ in its half-partition and corresponds to $(\sigma,\pi_k\otimes T)$
where $\sigma$ is now interpreted as a $P_n$ half-partition by removing
$n+1$ and $k$ is $1$ plus the number of propagating components containing
parts between $n+1$ and the largest element not equal to $n+1$ of the part
containing $n+1$.
Thus the quotient is isomorphic to a module with a series with sections
isomorphic to the
modules $\Delta_n(\mu)$ where $\mu-\lambda$ is of the form $(0,0,\ldots ,1,
\ldots ,0,0)$. Axiom (N4) is shown, and (N9) also, with corresponding rooted
graph $(\YounGG,0)$, the double Young graph.

Axiom (N5) is trivial from the definitions, and axiom (N8) follows from a
standard diagram algebra argument. \Qed

\Subsection{Pair partitions}
A subset of the set of partitions of $\underline{n}$ is the set
$j(n)$ of partitions into pairs, known as \emph{pair partitions}.
Obviously this set is empty unless $n$ is even, and we have 
 $|j(2m)| =  \frac{ (2m)!}{ 2^m m!}= (2m-1)(2m-3)...(1)$. 
Let $\CC_{br}(n)$ denote the set of pair partitions of $\underline{n}\cup 
\underline{n'}$.
Define a half pair partition of $\underline{n}$ to be a partition of 
$\underline{n}$ into pairs and singletons (known as propagating
components).
For $\lambda$ an integer partition of $l$ where $n-l$ is even,
let $Y_{br}(n;\lambda)$ denote the set of pairs $(\sigma,T)$ where $\sigma$
is a half pair partition of $\underline{n}$ with $l$ singletons and $T$ is a
standard tableau filling of the Young diagram corresponding to $\lambda$.


(1) Bra-ket decomposition. \\
Given a partition of $\underline{n}\cup \underline{n'}$ we take the set of
parts contained entirely within $\underline{n}$ as the pairs of a
half pair partition; the remaining elements are designated as singletons.
This gives a half pair partition; the same procedure is followed for
the parts contained in $\underline{n'}$ and in this way we obtain a pair
of half pair partitions. The correspondence between the propagating
components of the two partitions determines an element of the symmetric
group $S_l$ and thus, via the Robinson-Schensted correspondence, a pair
of standard tableaux of the same shape. The first is associated with the
first half pair partition, and the second with the second half pair
partition. In this way we obtain a pair of elements of $Y_{br}(n;\lambda)$
for some integer partition $\lambda$.

(2) Edge maps. \\
If $\lambda$ is an integer partition of $l\leq n$, let $\nu$ be an integer
partition such that $\mu-\lambda$ is of the form
$(0,0,\ldots ,1,0,0,\ldots )$. We define
$\edge{(n;\lambda),(n+1,\mu)}$ to be the map that takes a pair
$(\sigma,T)$ to the half partition with an extra propagating component
$\{n+1\}$ and the tableau of shape $\mu$ obtained from $T$ by adding a new
box labelled $l+1$.
Let $\nu$ be an integer partition such that $\lambda-\nu$ is of the form
$(0,0,\ldots ,1,0,0,\ldots )$. The map $\edge{(n;\lambda),(n+1,\nu)}$ is
defined entirely analagously to the map $\edge{(2n;\lambda),(2n+1,\nu+)}$
for the partition algebra case considered above.

We have:

\prl(p:braueralgebraarray)
The sets $Y_{br}(n;\lambda)$ described above form a
$(\YounG,0)$-array of sets with corresponding Catalan sequence of
of sets given by $\CC_{br}(n)$. \Qed
\end{pr}

\begin{figure}
\[
\includegraphics{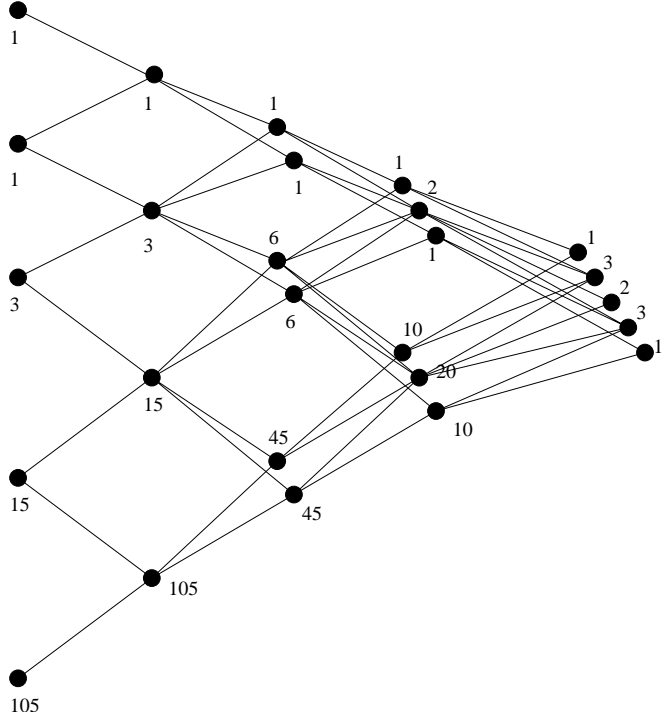}
\]
\caption{\label{BrauerBrattelii} Counting walks on the  Young graph.
That is, 
part of the `Bratteli' diagram
for the Brauer algebras $J_0(\delta) \subset J_{1}(\delta) \subset J_{2}(\delta) \subset J_{3}(\delta) \subset J_{4}(\delta) \cdots$
(to be precise, all restrictions
of the standard module $\SS({(0)}) (4)$),
in which each  $\SS(\lam)$  is represented by its dimension. }
\end{figure}

The start of the array is displayed in Figure~\ref{f:brauerarray}.
The sizes of some of the sets in the array are displayed in
Figure~\ref{BrauerBrattelii}. As a consequence of the above and
Proposition~\ref{p:equivalence} we have:


\begin{figure}
\[
\includegraphics{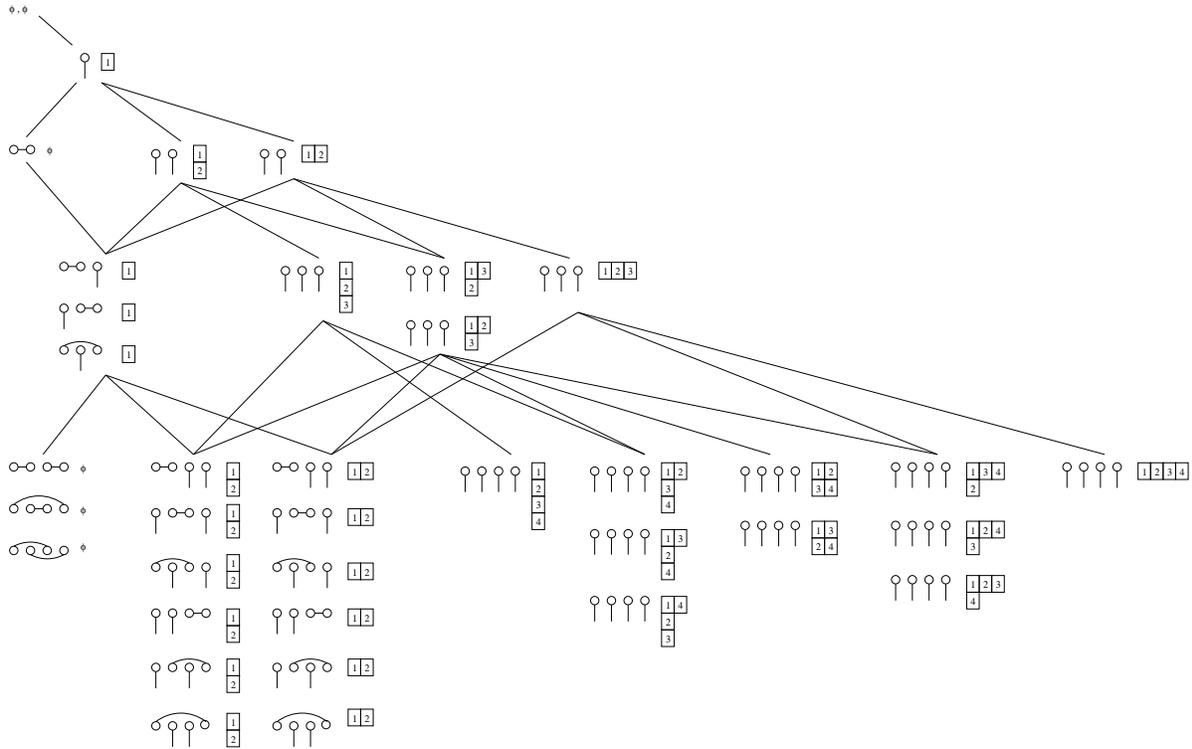}
\]
\caption{\label{f:brauerarray} The start of the array of half pair partitions.}
\end{figure}


\begin{prop}
Let $\YounG$ be the Young graph with distinguished vertex $0$ given by
the empty integer partition. Then the sequence $\{ N_{\YounG}(2n;0) \}_n$ is
the sequence 1, 1, 3, 15 , 105, ... , $(2m-1)(2m-3)...(1)$, ... 
\end{prop}
{\em Proof:} This is analogous to Proposition~\ref{dYg}
(cf. \cite{MartinRollet98}), using the Brauer algebra instead of the
partition algebra. 
\Qed

We remark that the above array also arises in an algebraic context.
The \emph{Brauer algebra} $J_n(\delta)$ has
$\mathbb{C}$-basis $\CC_{br}(n)$ with multiplication defined as for
$P_n$ (see also \cite{Weyl46}).
Each basis element is usually envisaged as a partition of an arrangement
of two rows of $m$ vertices into pairs. 
The following is implicit in~\cite{Brown55} (see
also~\cite{CoxDevisscherMartin05}).
The proof is analagous to that for the partition algebra case considered above.

\prl(Brauertower)
The tower of algebras $J_0(\delta)\subseteq J_1(\delta)\subseteq \cdots $
is a $(\YounG,0)$-tower of algebras, where
$\YounG$ is the Young graph with distinguished vertex $0$ given by the
empty partition. \Qed
\end{pr}


Remarks:
\newline
Both the Brauer and partition algebras have interesting non-semisimple
specialisations, associated to their roles in invariant theory. 
We briefly review this setting, and consider the combinatorics of
the non-semisimple cases, in Section~\ref{s:truncation}.


\Subsection{Weight lattice graphs}
Here we briefly touch on cases in which $N$ 
(from Section~\ref{new axioms}) is different from 2
(and the graph $G$ is directed).
The objective here is only to demonstrate that examples exist. 
We reserve a detailed exposition for a later work.

Our $(A_{\infty},0)$ graph may be considered as the dominant part of
the $\type = A_1$ weight lattice \cite{Bourbaki81}, with the root at the
 boundary of the dominant region, in the following sense.

For $N\in\mathbb{N}$, let $P(A_{N-1})$ denote the integral weight lattice of
type $A_{N-1}$, with the dominant integral weights denoted by $P^+(A_{N-1})$.
Let $e_1,e_2,\ldots ,e_N$ be the standard basis for $\mathbb{R}^N$, and let
$V$ be the subspace consisting of elements whose coordinates in this basis
sum to zero. Then $P_{A_{N-1}}$ can be realised in $V$ as the integral
linear combinations of the fundamental weights $\omega_i$, where
$\omega_i:=e_1+e_2+\cdots +e_i$, for $i=1,2,\ldots ,N-1$. The dominant
weights are the nonnegative integral combinations of the $\omega_i$.
Note also that the dominant weights are in one-to-one correspondance with the
$N$-row Young diagrams, which we denote by $\Lambda^N$.

Denote by $\lat(A_{N-1})$ the graph with vertices given by $P(A_{N-1})$
and an arrow from $\lambda$ to $\mu$ whenever $\lambda-\mu=e_i$ for some
$i$. Let $\lat^+(A_{N-1})$ denote the induced subgraph on the vertices
$P^+(A_{N-1})$. Each graph has distinguished vertex given by the zero
weight.

In the $\type=A_1$ case this simply says that each adjacent pair of vertices
has an edge between them in each direction (i.e. the graph is
undirected), and we thus see that $(\lat(A_1),0)$ is isomorphic
to $(A_{\infty}^{\infty},0)$ and $(\lat^+(A_1),0)$ is isomorphic to
$(A_{\infty},0)$.
The case $A_{N-1}$ is $N$--partite, rather than bipartite as all our
previous examples are. We shall return to this point shortly. 

These cases have an extensive literature
(see~\cite{MartinWoodcock03} for examples).

\begin{figure}
\includegraphics{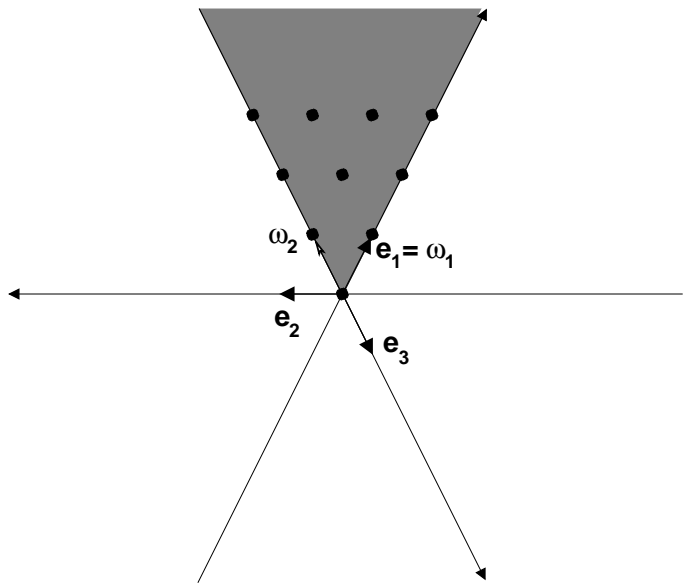}
\caption{\label{f:a2} The $A_2$ dominant integral weight lattice.}
\end{figure}
\begin{figure}
\includegraphics{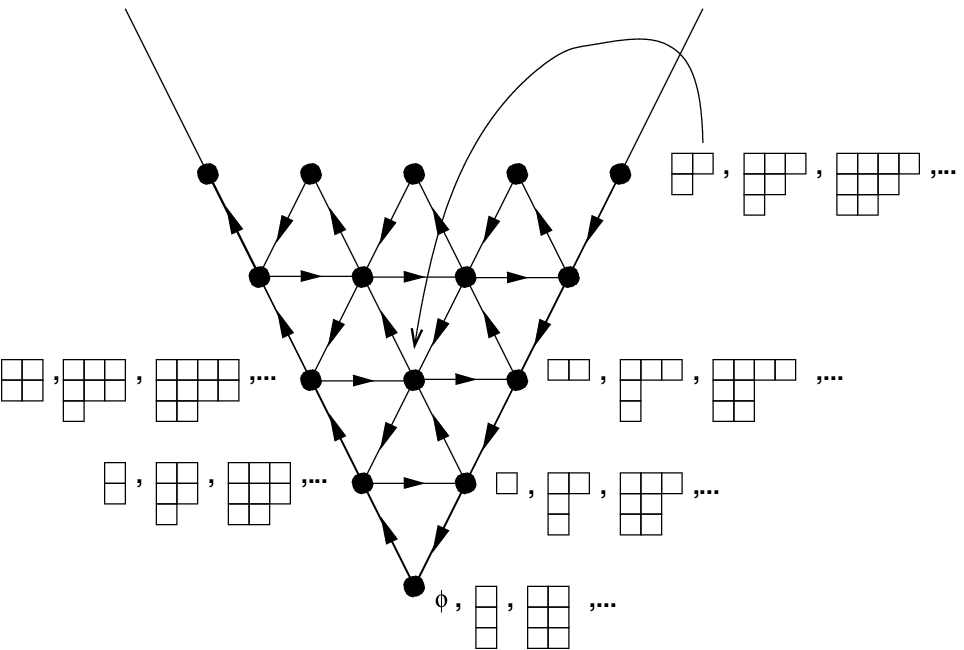}
\caption{\label{f:a2i} The graph $\lat^+(A_2)$.}
\end{figure}

The case for the $\type= A_2$ weight lattice is illustrated in
Figures~\ref{f:a2} and~\ref{f:a2i}.
In this case, it will be convenient to label the vertices as linear
combinations of fundamental weights, so that $d_1\omega_1+d_2\omega_2$
is denoted $(d_1,d_2)$. Thus the root is $(0,0)$;
the unique vertex one step from the root is $(0,1)$; the two vertices one
step on are $(1,0)$ and $(0,2)$.
In the third layer of the Pascal array there are sets labelled by the
vertices $(0,0)$, $(1,1)$ and $(0,3)$. 
The set $Y_{\lat^+(A_2)}(3;(1,1))$ has two paths in it, and the
others considered above have one. 
(Since the layers of the Pascal array are now effectively two-dimensional
we will not attempt to draw it.)

One way to put this in an algebraic context is to consider 
the action of the group algebra of the symmetric group $S_n$ on the tensor
space $(V_N)^{\otimes n}$, where $V_N$ is the $k$-vector space of dimension
$N$. Let us say for definiteness that $k=\C$, and we consider $N=3$. 
The former condition says that irreducible representations
$L_{\lambda}(n)$ of $S_n$ may be indexed by 
the set $\Lambda_n$ of integer partitions of $n$, or
equivalently by Young diagrams with $n$ boxes.  
For $n<N+1$ the action 
of $S_n$ on $(V_N)^{\otimes n}$ 
is faithful, and otherwise
the annihilator of this action is generated by the $S_{N+1}$
antisymmetriser $e_{(1^{N+1})} \in S_{N+1} \subseteq S_n$. 
This means that only 
irreducible representations indexed by 
the set $\Lambda^N_n$ of Young
diagrams with fewer than $N+1$ rows (and exactly $n$ boxes) appear in the
decomposition of tensor space. 
The dimensions of these representations are given by the hook law
(which counts the number of standard Young tableaux --- see
e.g. \cite{Hamermesh62}).
The point here is that they are also given, in case $N=3$, by 
\[
\dim L_{\lambda}(n) 
  = | Y_{\lat^+(A_2)}(n;(\lambda^f_1 -\lambda^f_2,\lambda^f_2)) |
\]
where $\lambda^f_i=\lambda_i - \lambda_3$.
For example, our $n=3$ layer as given above shows that 
$\dim L_{(2,1,0)}(3) 
  =  2  $. 
To see this let $H^N_n$ denote the quotient of $k S_n$ by its
annihilator on tensor space. Then there is an isomorphism
$H^N_{n} \cong e_{(1^N)} H^N_{n+N}  e_{(1^N)} $
\cite{Martin91}. This enables us to define a functor
\begin{eqnarray*}
F: H^N_{n+N}\modules  & \rightarrow &  H^N_{n}\modules
\\
M & \mapsto &  e_{(1^N)} M
\end{eqnarray*}
(where the action on $ e_{(1^N)} M$ is via the isomorphism).
This takes the simple module with label 
$(\lambda_1,\lambda_2,...,\lambda_N)\in \Lambda^N_n$ to one with label 
$(\lambda_1-1,\lambda_2-1,...,\lambda_N-1)$ 
(or to zero if $\lambda_N=0$). 
This says that the sequence of algebras $(H^N_n)_n$ has a global limit
with irreducible modules indexed by the set $\Lambda^{N-1}$ of $N-1$-row
Young diagrams --- which are in bijection with the set $P^+(A_{N-1})$  of 
dominant integral $A_{N-1}$-weights.
This brings us to an example of the setup in Section~\ref{new axioms}
with $N$ different from 2.
Let us write $\lambda^f$ for the image of $\lambda \in \Lambda^N_n$
in  $\Lambda^{N-1}$ (i.e. $\lambda^f_i=\lambda_i-\lambda_N$).
Thus for each $\mu \in \Lambda^{N-1}$ such that 
$|\mu| \leq n$ and $|\mu| \equiv n$ modulo $N$ we have an irreducible
representation
$
L[\mu](n) 
$
of $H^N_n$ 
given by $L[\lambda^f](n) = L_{\lambda}$. In Figure~\ref{f:a2i}, we have
labelled each $\mu\in \Lambda^{N-1}$ with the set of tableaux corresponding to
all $\lambda$ arising in this way.
We have (see, for example, \cite[7.3, Cor. 3]{Fultonyoungtableaux97})
\[
\res_{S_{n-1}} L[\mu](n) = \bigoplus_{\nu} L[\nu](n-1)
\]
where the sum is over weights which have an edge incoming from $\mu$
in $\lat^+(A_{N-1})$.
Now cf. Section~\ref{new axioms}. 

As with our previous examples there are a number of other
combinatorial sets related to sequences, Young diagrams and so on,
which are also counted in this Pascal array $Y_{\lat^+(A_2)}$. 
Sadly a diagram algebra style basis for our quotient algebra,
analogous to the \TL\ diagrams, is 
{\em not} known however. 
(We will return to this in another work.)

{\rem
In the manner of this section 
we may also embed the Young graph into $\lat^+(A_{\infty})$,
and hence give a conventional geometric setting for the index set 
for standard $P_n$--modules \cite{MartinWoodcock98}.
}

\section{Truncation, alternative roots, and non-semisimplicity}
\label{s:truncation}

The Stirling numbers of the second kind count set partitions with
certain extra properties \cite{Liu68,MartinRollet98}. 
It is interesting both from 
the combinatorial and representation theoretic point of
view to ask if walks with special properties on the double Young graph
can be put in correspondence with these special partitions, as a
refinement of the Bell sequence analysis above.
The answer is that the appropriate walks are those which are
restricted to a certain subgraph --- a {\em truncation}. 
This truncation then also has a role in the representation theory of
the partition algebra for non-semisimple values of the parameters.
This case was first treated in \cite{MartinRollet98}.
Here we begin with a much simpler example of such a truncation,
which also illustrates the  combinatorial and representation theoretic
effects of `moving' the root of the graph. 

\Subsection{\TL\ and blob cases}
In this section we reparameterise $TL_n(\delta)$ by $\delta=q+q^{-1}$. 
When $q$ is a (primitive $l$-th) root of unity the tower
  $(TL_n(q+q^{-1}))_{n\in\N}$ 
from Section~\ref{TLA}
is not semisimple, and the standard modules are not all simple.
It is interesting to ask what variant of our combinatoric can describe
  the simple modules in these cases. 
Here we recall the case over the complex field. In this case the only
  parameter is $l$. A convenient geometrical starting point is to
  embed $(A_{\infty},0)$ in  $(A_{\infty}^{\infty},0)$ in the obvious
  way. One then views $TL_{\bullet}$ as a quotient of the blob algebra
$b_{\bullet}$ (Section~\ref{s:blob}) by the relation 
\[
e = 0
\]
(For such a non-zero quotient we need $\delta' = 0$). 
Only one standard module in each layer respects this quotient, but the
consistency requirement $ \delta' = 0$ is not a semisimple
specialisation, and, appropriately arranged, the simple modules
corresponding to the vertices of 
$(A_{\infty},0)  \hookrightarrow (A_{\infty}^{\infty},0)$
do respect the quotient. 


\prl(goopiod){\rm \cite{MartinWoodcock2000}}
Fix $\delta'=0$. 
Walks on  $(A_{\infty}^{\infty},0)$ from 0 to $\lambda$ of length $n$
which never touch vertex $-1$ are a basis for the 
$\delta$-generically simple module 
$L^{b_n}(\lambda)=\Head(\Delta^{b_n}(\lambda) )$. \Qed
\end{pr}
Of course this restriction on walks precisely corresponds to the
embedding of  $(A_{\infty},0)$ in  $(A_{\infty}^{\infty},0)$ 
and we have that these modules are the simple modules of the 
$TL_{\bullet}(q+q^{-1})$ quotient. 
Indeed these bases make sense for any $q$, and give the standard
modules of $TL_{\bullet}(q+q^{-1})$.

Next we consider $l \in \N$.
\prl(goopiod2){\rm \cite{Martin91}}
Let $q$ be  a primitive $l$-th root of unity, and let $\lambda\in \mathbb{Z}$.
Let $m\in \mathbb{Z}$ be minimal such that $\lambda<ml-1$.
Walks on  $(A_{\infty}^{\infty},0)$ from 0 to $\lambda \geq 0$ of length $n$
which never touch vertex $-1$ 
and never touch vertex $ml-1$ unless they subsequently touch $(m-1)l-1$ 
are a basis for the simple module 
$L^{TL_n}(\lambda)$ (note that this last restriction is vacuous if
$\lambda=ml-1$ for some $m$). \Qed
\end{pr}

Example: Suppose $\lambda < l-1$. Since the walks in the proposition
never touch $-1$ they may never touch $l-1$, so we have the set of walks on
the truncated graph 
$(A_l,0)  \hookrightarrow (A_{\infty},0)$.

\begin{figure}
\includegraphics{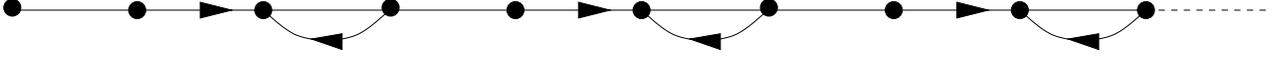}
\caption{\label{rollet l} Rollet for $l=3$.
Undirected edges have arrows in both directions.}
\end{figure}

Another natural question at this point, then, is:
Is there a Rollet graph {\em unrestricted} walks on which encode
the restriction of simple modules over the Temperley-Lieb algebra?
To answer this, note the following alternative statement of
Prop~\ref{goopiod2}: For $r=-1,-2,\ldots ,l-2$,
\[
Res L(ml+r) = \left\{ \begin{array}{lll}
L(ml+r+1) + L(ml+r-1) & r \neq -1,l-2  \\
L(ml+r-1)             & r=l-2   \\
2L(ml+r+1)+L(ml+r-1)  & r=-1
\end{array} \right. 
\] 
(where $L(s)$ is regarded as zero for negative $s$).
This gives us the Rollet graph in Figure~\ref{rollet l}. 

This in turn raises another natural question. 
The new Rollet graph is similar to that for standard modules ---
i.e. it is obtained from $A_{\infty}$ by borrowing one edge and giving
it to the next pair of vertices, at three step intervals. 
We have equation~(\ref{sumPL}), 
showing that there is a kind of combinatorial duality between the sets
of simple and projective modules. So the question is:
Can we construct a `dual' graph, walks on which give bases for
projectives? 
This is a question of a different nature, since unsurprisingly
\prl(nogo)
Not every \TL\ projective module restricts to a module with a
projective filtration. \Qed
\end{pr}
(And induction of \TL\ modules is a bit of an art.)
It is possible to give projective decompositions of restricted
projective modules up to some `boundary effects' (rather like \TL\ tilting
modules \cite{Donkin93,Martin92}), 
that is, to adding in one or two not-quite
projective modules in finitely many cases.
In the interests of brevity we will return to this interesting
question in a separate paper.  


Having given some examples of truncation and non-semisimplicity, now
now turn to the algebra of shifted roots. 


Set $[n]=q^{n-1}+q^{n-3}+\ldots+q^{1-n}$. 
Returning to the blob algebra, with $\delta=[2]$ as above
it is convenient to reparameterise $\delta'=[m+1]/[m]$.
The blob algebra $b_n$ is then semisimple (over $\C$) unless $m$ is
integral (see \cite{MartinSaleur94a,MartinWoodcock2000}). 
Note that $\delta'=0$ is the case $m=-1$. 
For each $b_n$ note that there are two standard modules of dimension
one, that is $\Delta^b(n)$ and $\Delta^b(-n)$. Since the algebras are
generically semisimple we may associate a primitive central idempotent
$E(n)$ (resp. $E(-n)$)
to each of these modules, which is well defined for $q$ and $m$
indeterminate (but which is `singular' in certain specialisations). 
For example, the idempotent for $\Delta(-1)$ is $E(-1)=e$ .
It is always possible to find a scalar multiple 
$E'(-n)$ of $E(-n)$ which is
well defined on specialising $m$ to an integer $l_0$ (say)
(leaving $q$ indeterminate). 
We may then define a quotient tower of $b_{\bullet}$ by $E'(-n) \equiv 0$
(for fixed $n$), analogous to the \TL\ quotient (the case $n=1$).
We have
\prl(goopiodx){\rm \cite{MartinWoodcock2000}}
Fix $\delta'=[l_0+1]/[l_0]$, where $l_0<0$ is integral.
Walks on  $(A_{\infty}^{\infty},0)$ from 0 to $\lambda$ of length $n$
which never touch vertex $l_0$ are a basis for the 
$\delta$-generically simple module 
$L^{b_n}(\lambda)=\Head(\Delta^{b_n}(\lambda) )$. 
These are the simple modules of the quotient by  $E'(l_0)=0$. \Qed
\end{pr}
Details of the associated representation theory are given in
\cite{MartinWoodcock2000}. The point for our purposes is simply this.
The tower of quotient algebras obtained by quotienting out by
$E'(l_0)$ (for $l_0<0$ integral) is a $G$-tower with $G= (A_{\infty},-l_0-1)$.
To see this we simply note that the truncation is at $l_0$ rather than
at -1.  This is a shift of $l_0+1$, which may be seen as 
a different embedding of $A_{\infty}$ in $A^{\infty}_{\infty}$, 
shifting the root by $l_0+1$. 


\Subsection{Brauer case}
For $\delta \in \N$ let $V_{\delta}$ be the (complex) defining module for the
action of the general linear group $GL(\delta)$. 
Then the natural action of $GL(\delta)$ on $V_{\delta}^{\otimes m}$ is
in Schur--Weyl duality with an action of the symmetric group $S_m$
permuting the tensor factors \cite{Weyl46}. 
Restricting to the action of the 
orthogonal group $O(\delta) \subset GL(\delta)$, 
then the Brauer algebra $J_m(\delta)$ has a natural action on
$V_{\delta}^{\otimes m}$ extending the $S_m$ action, 
which commutes with the  $O(\delta)$ action. 
This was the context in which the Brauer algebra was originally
introduced. 
Restricting further to the action of $S_{\delta} \subset O(\delta)$
{\em permuting} standard ordered basis elements in each $V_{\delta}$, 
the dual action extends to one of the partition algebra. 

These specialisations of $\delta$ are non-semisimple cases of the
Brauer algebra, so it is natural to ask if there are 
`truncations' of
the set of walks on the Young graph which enumerate the bases of
simple modules in these cases.
In general this seems to be a very hard problem (see
\cite{CoxDevisscherMartin05}). 
A simpler (but closely related) 
problem is to restrict attention to the simple modules which appear in
the tensor space module. 

The simplest case is $\delta=1$. Here only the (one-dimensional)
simple module in the
head of the spine standard module (with $\lambda=0$) survives. 
The truncation is to walks on the full subgraph containing the
vertices for $S_0$ and $S_1$. There is only one such walk in each
layer, so this truncation works almost trivially. Subsequent cases
are less straightforward to describe, and we reserve details on them for
a separate work.



\section{A
$(D_{\infty},0)$ tower of algebras}\label{s:dnalgebra}
\Subsection{A blob algebra in type $D$}
The defining Pascal $D$--array is sketched in Figure~\ref{D-seq walk}.
More precisely, this is the Pascal array for the rooted graph
$(D_{\infty},0)$ where $D_{\infty}$ is as in Figure~\ref{f:Dinfinity}.

To obtain such an array in an algebraic setting we consider the algebra
$d_n$ introduced by Richard Green in \cite[\S 7]{Green98}. A diagram basis of
$d_n$ is given by the Temperley-Lieb diagrams with $2n$ vertices
(divided into north and south vertices as in the Temperley-Lieb case).
As in the blob algebra case, there is the additional possibility of
decorating any arc exposed to the western end of the diagram. There is an
additional restriction: there must be an even number of blobs in total.

\begin{figure}
\setlength{\unitlength}{1cm}
\begin{picture}(8,2)
\put(1,1){\circle{1}}
\put(1.75,1){$\displaystyle =\delta,$}
\put(3.5,1){\circle{1}}
\put(4,1){\circle*{0.2}}
\put(4.25,1){$\displaystyle =0,$}

\put(5.7,0){\line(0,2){2}}
\put(5.7,0.7){\circle*{0.2}}
\put(5.7,1.3){\circle*{0.2}}
\put(6.2,1){$=$}
\put(7,0){\line(0,2){2}}
\end{picture}
\caption{\label{dnrules} The rules for simplifying diagrams for the
algebra $d_n$.}
\end{figure}
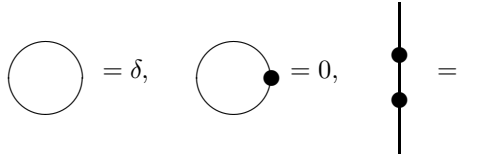

Multiplication is given by concatenation of diagrams in the usual way;
in addition, any undecorated loop is replaced with the scalar $\delta$,
any loop with an odd number of blobs is zero and any pair of
blobs on the same arc can be removed; see Figure \ref{dnrules}.
We assume in the sequel that $\delta\not=0$.

We denote the set of $d_n$ diagrams with $n$ vertices on the north edge
and $m$ on the south edge by $D^d(n,m)$. Thus a basis for $d_n$ is
$D^d(n,n)$.
We still have the bra-ket construction as in Section~\ref{s:blob}
above, with the
proviso that if the two diagrams produced have at least one undecorated
propagating line and an odd number of blobs in total, the westernmost
propagating line is decorated with a blob on each of the half diagrams
(and any pair of blobs on a single arc produced in this way is removed).
We denote the half diagrams obtained from $D^d(n,n)$ with $l$
propagating lines by $D^d_l(n,l)'$; note that if $l>0$ then these half
diagrams always have an even number of blobs.
The constructions in (2) and (3) in the previous section go through
unchanged.

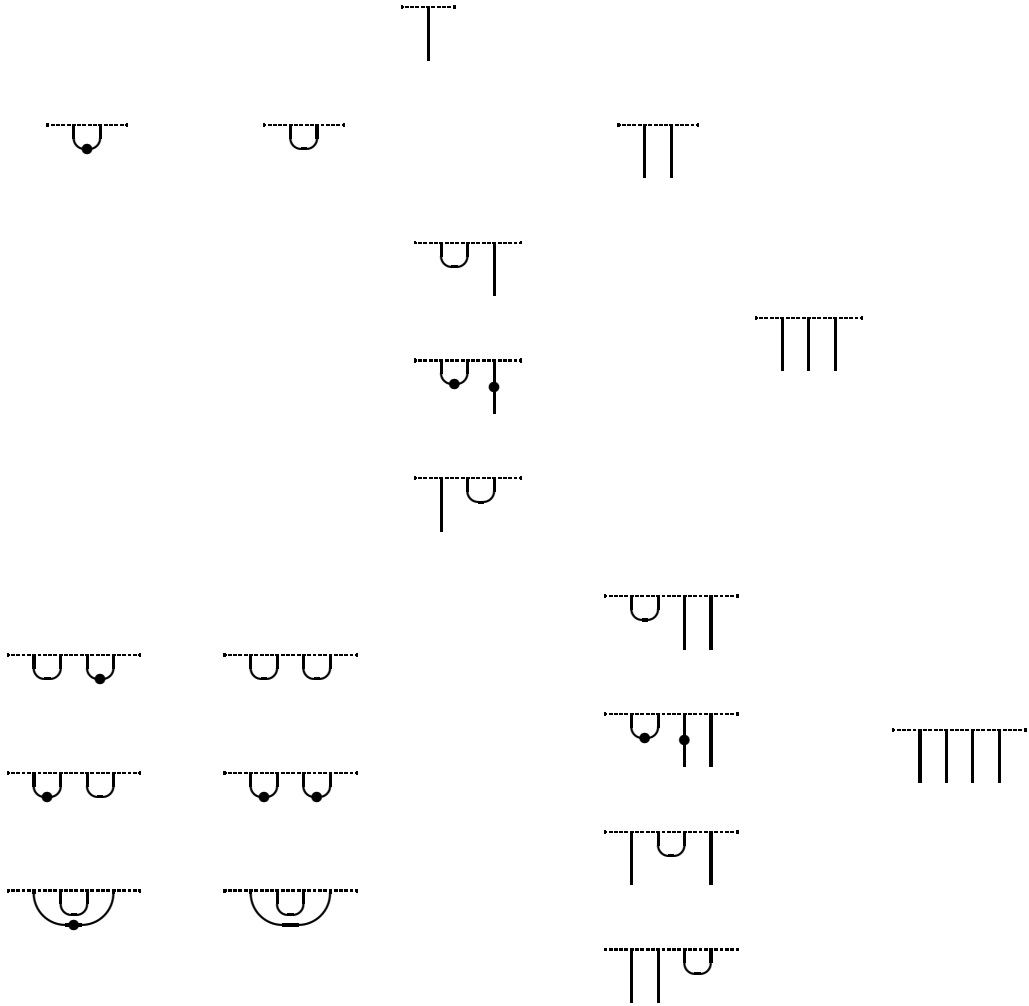
\begin{figure}
\[ \begin{array}{lllllll}
&&
\begin{picture}(22,41)(-0,-40)
{\put(20,0){
\put(0, 0){\dashbox{1.0}(20,0){}}
\thicklines
\put(10, 0){\line(0,-1){20}}
}}
\end{picture}
\\ 
\begin{picture}(22,41)(-0,-40)
{\put(20,0){
\put(0, 0){\dashbox{1.0}(30,0){}}
\thicklines
\put(15.0, 0){\oval(10,18.0)[b]}
\put(15.0,-9.0){\circle*{4}}
}}
\end{picture}
& 
\begin{picture}(82,41)(-40,-40)
{\put(20,0){
\put(0, 0){\dashbox{1.0}(30,0){}}
\thicklines
\put(15.0, 0){\oval(10,18.0)[b]}
}}
\end{picture}
& &
\begin{picture}(82,41)(-40,-40)
{\put(20,0){
\put(0, 0){\dashbox{1.0}(30,0){}}
\thicklines
\put(10, 0){\line(0,-1){20}}
\put(20, 0){\line(0,-1){20}}
}}
\end{picture}
\\
&&
\begin{array}{l}
\begin{picture}(22,41)(-0,-40)
{\put(20,0){
\put(0, 0){\dashbox{1.0}(40,0){}}
\thicklines
\put(15.0, 0){\oval(10,18.0)[b]}
\put(30, 0){\line(0,-1){20}}
}}
\end{picture}
\\
\begin{picture}(22,41)(-0,-40)
{\put(20,0){
\put(0, 0){\dashbox{1.0}(40,0){}}
\thicklines
\put(15.0, 0){\oval(10,18.0)[b]}
\put(15.0,-9.0){\circle*{4}}
\put(30, 0){\line(0,-1){20}}
\put(30,-10.0){\circle*{4}}
}}
\end{picture}
\\
\begin{picture}(22,41)(-0,-40)
{\put(20,0){
\put(0, 0){\dashbox{1.0}(40,0){}}
\thicklines
\put(10, 0){\line(0,-1){20}}
\put(25.0, 0){\oval(10,18.0)[b]}
}}
\end{picture}
\end{array}
&&
\begin{picture}(22,41)(-0,-40)
{\put(20,0){
\put(0, 0){\dashbox{1.0}(40,0){}}
\thicklines
\put(10, 0){\line(0,-1){20}}
\put(20, 0){\line(0,-1){20}}
\put(30, 0){\line(0,-1){20}}
}}
\end{picture}
\\
\begin{array}{l}
\begin{picture}(22,41)(20,-40)
{\put(20,0){
\put(0, 0){\dashbox{1.0}(50,0){}}
\thicklines
\put(15.0, 0){\oval(10,18.0)[b]}
\put(35.0, 0){\oval(10,18.0)[b]}
\put(35.0,-9.0){\circle*{4}}
}}
\end{picture}
\\
\begin{picture}(22,41)(20,-40)
{\put(20,0){
\put(0, 0){\dashbox{1.0}(50,0){}}
\thicklines
\put(15.0, 0){\oval(10,18.0)[b]}
\put(15.0,-9.0){\circle*{4}}
\put(35.0, 0){\oval(10,18.0)[b]}
}}
\end{picture}
\\
\begin{picture}(22,41)(20,-40)
{\put(20,0){
\put(0, 0){\dashbox{1.0}(50,0){}}
\thicklines
\put(25.0, 0){\oval(30,26.0)[b]}
\put(25.0,-13.0){\circle*{4}}
\put(25.0, 0){\oval(10,18.0)[b]}
}}
\end{picture}
\end{array}
&
\begin{array}{l}
\begin{picture}(22,41)(-20,-40)
{\put(20,0){
\put(0, 0){\dashbox{1.0}(50,0){}}
\thicklines
\put(15.0, 0){\oval(10,18.0)[b]}
\put(35.0, 0){\oval(10,18.0)[b]}
}}
\end{picture}
\\
\begin{picture}(22,41)(-20,-40)
{\put(20,0){
\put(0, 0){\dashbox{1.0}(50,0){}}
\thicklines
\put(15.0, 0){\oval(10,18.0)[b]}
\put(15.0,-9.0){\circle*{4}}
\put(35.0, 0){\oval(10,18.0)[b]}
\put(35.0,-9.0){\circle*{4}}
}}
\end{picture}
\\
\begin{picture}(22,41)(-20,-40)
{\put(20,0){
\put(0, 0){\dashbox{1.0}(50,0){}}
\thicklines
\put(25.0, 0){\oval(30,26.0)[b]}
\put(25.0, 0){\oval(10,18.0)[b]}
}}
\end{picture}
\end{array}
&&
\begin{array}{l}
\begin{picture}(22,41)(-30,-40)
{\put(20,0){
\put(0, 0){\dashbox{1.0}(50,0){}}
\thicklines
\put(15.0, 0){\oval(10,18.0)[b]}
\put(30, 0){\line(0,-1){20}}
\put(40, 0){\line(0,-1){20}}
}}
\end{picture}
\\
\begin{picture}(22,41)(-30,-40)
{\put(20,0){
\put(0, 0){\dashbox{1.0}(50,0){}}
\thicklines
\put(15.0, 0){\oval(10,18.0)[b]}
\put(15.0,-9.0){\circle*{4}}
\put(30, 0){\line(0,-1){20}}
\put(30,-10.0){\circle*{4}}
\put(40, 0){\line(0,-1){20}}
}}
\end{picture}
\\
\begin{picture}(22,41)(-30,-40)
{\put(20,0){
\put(0, 0){\dashbox{1.0}(50,0){}}
\thicklines
\put(10, 0){\line(0,-1){20}}
\put(25.0, 0){\oval(10,18.0)[b]}
\put(40, 0){\line(0,-1){20}}
}}
\end{picture}
\\
\begin{picture}(22,41)(-30,-40)
{\put(20,0){
\put(0, 0){\dashbox{1.0}(50,0){}}
\thicklines
\put(10, 0){\line(0,-1){20}}
\put(20, 0){\line(0,-1){20}}
\put(35.0, 0){\oval(10,18.0)[b]}
}}
\end{picture}
\end{array}
&&
\begin{picture}(82,41)(-20,-40)
{\put(20,0){
\put(0, 0){\dashbox{1.0}(50,0){}}
\thicklines
\put(10, 0){\line(0,-1){20}}
\put(20, 0){\line(0,-1){20}}
\put(30, 0){\line(0,-1){20}}
\put(40, 0){\line(0,-1){20}}
}}
\end{picture}
\end{array} \]
\caption{\label{D blob} $D$--type blob `bra' diagrams. 
These sets have the property that if combined with their 
corresponding `ket' set they produce diagrams with an even number
of blobs.}
\end{figure}

See Figure~\ref{D blob} for an illustration of the Pascal array
in this case.

For $l\in\{n,n-2,\ldots, \mbox{\ $2$\ or\ $1$}\}$, we define
$D^d_l(n,l)=D^d_l(n,l)'$. If $n$ is even, we define $D^d_0(n,0)$
to be the set of half diagrams in $D^d(n,0)'$ with an even number of blobs,
and $D^d_{0'}(n,0')$
to be the set of half diagrams in $D^d(n,0)'$ with an odd number of blobs.


\begin{prop} \label{Dblobarray}
The array of sets $((D^d_l(n,l))_{l\in D_{\infty}})_n$
is a Pascal $D^{\infty}$--sequence.
(a) If $l\not=0'$, the edge maps $\edge{i,i+1}$ are given
by $\phi^1$ and the edge maps $\edge{i+1,i}$ are given by $\phi^u$. \\
(b) The edge maps $\edge{0',1}$ are given by $\phi^{1'}$, i.e. $\phi^1$
modified by decorating the new propagating edge with a blob. \\
(c) The edge maps $\edge{1,0'}$ are given by $\phi^{u'}$, i.e. $\phi^u$
modified by decorating the propagating edge that is ``bent over''
(to form a new vertex) with an extra blob (so if a blob is already there, it is
removed by the relation involving two blobs).
\end{prop}

{\bf Proof:} It suffices to note that \\
(a) If $l\in \{n,n-2,\ldots \mbox{3\ or\ 2}\}$ then
$$
D^d_l(n,l)
\; = \; \phi^u( D^d_{l+1}(n-1,l+1) ) \; \bigcup \; \phi^1(
D^d_{l-1}(n-1,l-1) ).$$
(b) If $l=1$ then
$$D^d_l(n,1)
\; = \; \phi^u( D^d_{2}(n-1,2) ) \; \bigcup \; \phi^1(
D^d_{0}(n-1,0) \; \bigcup \; \phi^{1'}(D_{0'}(n-1,0') ).$$
(c) If $l=0$ then
$$
D^d_0(n,0)
\; = \; \phi^u( D^d_{1}(n-1,1) ).$$
(d) If $l=0'$ then
$$
D^d_0(n,0)
\; = \; \phi^{u'}( D^d_{1}(n-1,1) ).$$
\Qed

It can be shown that the axioms (A1) to (A4), (A4'), (A5)(a) and
(N5) hold. This has the consequence that the axioms (N1) to (N8) all
hold for the sequence of algebras
$d_1,d_2,\ldots $ (note that we start our numbering at $1$ rather than
$0$). We remark that a consequence of this is that the
$d_n$ are all quasihereditary.
In this case the standard module $\Delta_n(l)$ has a basis parametrized
by $D^d_l(n,l)$.

By Proposition~\ref{ToR-GT} it follows that the algebras $d_{\bullet}$
(with the standard modules) form a $(G,0)$-tower of algebras for some rooted
graph $(G,0)$. Analysis of the restriction of the standard modules
for $d_n$ to $d_{n-1}$ shows that $(G,0)=(D_{\infty},0)$
(see Figure~\ref{f:Dinfinity}),
and hence that the array associated to the tower
$d_{\bullet}$ is a Pascal $(D_{\infty},0)$-array of sets and thus equivalent
to the array of sets $((D^d_l(n,l))_{l\in D_{\infty}})_n$ shown in
Figure~\ref{D blob}. We thus have:

\begin{theorem}
The tower of algebras $d_{\bullet}$, together with their standard
modules, is a $(D_{\infty},0)$-tower of algebras. 
\end{theorem}

We see that the basis of $\Delta_n(\lambda)$ provided by
Proposition~\ref{towerbases} can be explicitly parametrized by the half
diagrams $D^d_l(n,l)$, recovering the above description of the modules
$\Delta_n(l)$.



\section{Final Remarks}

Recall that 
the Robinson--Schensted correspondence can be regarded as giving a
bijection between pairs of paths (of length $n$ from the empty diagram
to the same vertex in the Young graph) and permutations of $n$ elements.
The Young graph is a oriented graph with a natural $\mathbb{Z}$-grading.
Furthermore, if $g_1,g_2$ are vertices then the set of
vertices $g$ such that $g\rightarrow g_1$ and $g\rightarrow g_2$ has
the same cardinality (which must be $0$ or $1$) as the set of vertices
$g$ such that $g_1\rightarrow g$ and $g_2\rightarrow g$. In addition,
for any vertex $g$, the number of vertices $h$ such that
$g\rightarrow h$ is equal to $1$ more than the number of vertices $k$
such that $k\rightarrow g$. A graph with all of these properties
is known as a \emph{$Y$-graph}.

Fomin \cite{Fomin88} has shown that a generalised Robinson-Schensted
correspondence can be defined on such a graph $G$ provided it has a unique
source $0$, giving a bijection between pairs of paths of length $n$ from
$0$ to the same vertex and permutations of $n$ elements. The graph $G$
could be, for example, the Young graph or the Young-Fibonacci graph
considered by Fomin. In this case the Catalan sequence of sets given by
the symmetric groups has underlying Pascal $(G,0)$-array $Y_{G,0}$,
and \cite[K7]{Fomin88} can be regarded as giving a natural construction of
the bra-ket extraction and its inverse. 
The result \cite[Theorem C]{Fomin88}
does not obviously fit into the general framework considered here,
and suggests the consideration of more general (sequences of) paths.


\appendix
\section{Generating functions 
for Catalan numbers}
\label{s:generatingfunctions}

Let $H_0(x)=\sum_{n\geq 0} C(n) x^n$ be the ordinary generating
function for the Catalan numbers. 
Since the Catalan numbers count rooted planar trees with $n$ vertices,
we see that
the constant term is $1$, since there is a unique tree with no edges. 
Every other tree has at least one edge. 
Having drawn this edge the remainder of every tree which contains it as a
subtree can be decomposed as two trees: one growing from the end of
this edge, and one growing from the root (on the right, say). Thus
we have, solving for $H_0$
(e.g. \cite{BloteNightingale82}):
$$
H_0 = 1+xH_0^2=\frac{1-\sqrt{1-4x}}{2x}.
$$
We can think of each such tree with $n+1$ vertices 
as a forest of rooted trees with $n$ vertices in total embedded in
the upper half-plane, each with its root lying on the boundary.
We do this by simply moving the boundary up to the line of vertices
one-removed from the original root and then discarding the root.
Note that this does not work for the 1-vertex tree, since this procedure
would give us the empty tree, which is not rooted.

Now suppose in the forest picture that there are actually two
different types of tree in the forest. They look the same as
trees, but a forest with two pines is different from a forest with a
pine and an oak, and so on. Let $H_1$ denote the generating function for
such two-species forests. As the two-species forests are counted by the
$(A^{\infty}_{\infty},0)-$, or $(\groot(2,1),0)$-Catalan numbers,
$H_1$ is the generating function for these.

We note that such a forest can be regarded as a tree of one of the two
species together with a remaining forest. By an analogous argument to
the above we have that
\[
H_1 =  1 + 2x H_0 H_1 = \;\; \frac{1}{\sqrt{1-4x}}
\]
We have proved:
\begin{prop}
The sequence 
 $\{ N_{(A^{\infty}_{\infty},0)}(2n;v_0) \}_n$ 
is the sequence generated by $H_1$. \Qed
\end{prop}
Viewed in isolation, this is nothing but
a well known result for the Pascal triangle.
But we can continue. 

Similarly the  $(\groot(2,2,1),0)$ Catalan numbers are given by 
\[
H_2 = 1 + 2x H_1 H_2 \;\; = \;\; \frac{1}{1-\frac{2x}{\sqrt{1-4x}}}
\]
and so on. Thus  the  $(\groot(2^l,1),0)$, Catalan numbers are given by 
\[
H_l = 1 + 2x H_{l-1} H_l   
\]

Let us write $H^{\lambda}$ for the generating function for the 
$(\groot(\lambda),0)$-Catalan numbers. We have:
\begin{prop}
Let $G$ be the graph $\groot(\lambda)$. Then the generating
function $H^{\lambda}$ for the sequence $\{ N_G(2n;v_0) \}_n$
satisfies the relation
\[
H^{\lambda} = 1+ \lambda_1 x H^{\lambda'} H^{\lambda}
\]
where $\lambda' = (\lambda_2,\lambda_3,...)$ and $H^{\emptyset}=H_0$. \Qed
\end{prop}

Finally, 
recall that the exponential generating function for the Bell numbers
is 
\[
\sum_{n \geq 0} \frac{b(n) x^n}{n!} \; = \; \exp(\exp(x) -1)
\; = \; 1+x+\frac{2x^2}{2}+\frac{5x^3}{3!}+ \frac{15 x^4}{4!} +...
\]
(see for example \cite{Wilf94}), 
and that the Rollet graph for the Bell
sequence considered in Section~\ref{bellnumbers} was not a tree.
It is intriguing
to speculate on the relationship between the generating function of a
Catalan sequence of sets and the properties of the Rollet graph of
a Pascal array for the sequence. For example, it may be that ordinary
functions correspond to the case where the Rollet graph is a tree,
with exponential generating functions corresponding to the case where
the Rollet graph is not a tree.


\bibliographystyle{amsplain}
\bibliography{new31,catalan4}

\end{document}